\documentclass[pdflatex,sn-mathphys-num]{sn-jnl}% Math and Physical Sciences Numbered Reference Style
\usepackage{graphicx}%
\usepackage{multirow}%
\usepackage{amsmath,amssymb,amsfonts}%
\usepackage{amsthm}%
\usepackage{mathrsfs}%
\usepackage[title]{appendix}%
\usepackage{xcolor}%
\usepackage{textcomp}%
\usepackage{manyfoot}%
\usepackage{booktabs}%
\usepackage{algorithm}%
\usepackage{algorithmicx}%
\usepackage{algpseudocode}%
\usepackage{listings}%
\usepackage{indentfirst}
\newcommand{\intO}{\int_{\Omega}}

\DeclareMathOperator{\divv}{div}

\newcommand{\uu}{{\mathbf u}}

\newcommand{\del}{\partial}
\theoremstyle{thmstyleone}%
\newtheorem{theorem}{Theorem}%  meant for continuous numbers
\newtheorem{proposition}[theorem]{Proposition}%

\newtheorem{lemma}[theorem]{Lemma}

\newtheorem{remark}{Remark}%

\theoremstyle{thmstylethree}%
\newtheorem{definition}{Definition}%

\begin{document}
	\title{ Local-in-time existence of strong solutions to a class of compressible Power-Law flows }
	%%=============================================================%%
	%% GivenName	-> \fnm{Joergen W.}
	%% Particle	-> \spfx{van der} -> surname prefix
	%% FamilyName	-> \sur{Ploeg}
	%% Suffix	-> \sfx{IV}
	%% \author*[1,2]{\fnm{Joergen W.} \spfx{van der} \sur{Ploeg}
		%%  \sfx{IV}}\email{iauthor@gmail.com}
	%%=============================================================%%
	\author*[1]{\fnm{Fang Li}}\email{fangli@nwu.edu.cn}
	\author[1]{\fnm{Chang Mengge} }\email{changmg@stumail.nwu.edu.cn}
	\author[2]{\fnm{Guo Zhenhua} }\email{zhguo@gxu.edu.cn}
	\equalcont{These authors contributed equally to this work.}

	\affil*[1]{\orgdiv{School of Mathematics and CNS}, \orgname{Northwest University}, \city{Xi'an}, \postcode{710000},
		\country{People\textquoteright s Republic of China}}
	
	\affil*[2]{\orgdiv{School of Mathematics and Center for Applied Mathematics of Guangxi}, \orgname{Guangxi University}, \city{Nanning}, \postcode{530004},
		\country{People\textquoteright s Republic of China}}

	%%==================================%%i'an
	%% Sample for unstructured abstract %%
	%%==================================%%
	
\abstract{
We consider a model of the compressible non-Newtonian fluids for power-law flow fulfilling a periodic domain in ${\mathbb R}^3,$ in which the extra stress tensor is induced by a potential with $p(t,x)$-structure. The local-in-time existence and uniqueness of strong is proved for all
$\frac{7}{5} < \inf p(t,x) \leqslant \sup p(t,x) \leqslant 2.$ Further, two improved blow-up criterions for the obtained strong solution is given.}
	
	%%================================%%
	%% Sample for structured abstract %%
	%%================================%%
	
	
	%%\pacs[JEL Classification]{D8, H51}
	
	%%\pacs[MSC Classification]{35A01, 65L10, 65L12, 65L20, 65L70}
	\maketitle	
	
	\section{Introduction}
	
Let $\Omega\subset {\mathbb R}^3$ be a domain occupied by a fluid. Let $\rho$ and ${\mathbf u}$ be the density and velocity of the fluid, respectively.  Let $t$ and $x$ denote the time and spatial variables, respectively. Then the Navier-Stokes equations, describing the flow of a compressible non-Newtonian fluid, read as
\begin{eqnarray}
\begin{cases}\label{XT-0}
	\partial_{t}\rho+\divv(\rho{\mathbf u})=0\mbox{ in }(0,T)\times\Omega,\\
	\partial_{t}(\rho{\mathbf u})+\divv(\rho{\mathbf u}\otimes{\mathbf u} )-\divv{\mathcal A}({\mathbb D}{\mathbf u})+\nabla \rho^\gamma
     =\rho{\mathbf f}\mbox{ in }(0,T)\times\Omega.
\end{cases}
\end{eqnarray}
Here $T>0$ is the time of evolution, ${\mathbf f}$ is the external force, $\gamma>1,$
${\mathbb D}{\mathbf u}=\frac{1}{2}(\nabla{\mathbf u}+\nabla^{T}{\mathbf u}),$ and the non-Newtonian tensor ${\mathcal A}(\mathbb{ D}\uu)$ is given by the following power law
\begin{align}\label{S}
	{\mathcal A}({\mathbb D}{\mathbf u})=(1+|{\mathbb D}{\mathbf u}|^2)^{\frac{p-2}{2}}{\mathbb D}{\mathbf u}
\end{align}
with $1<p<+\infty.$
	
Recent decades have witnessed a revival of interest in the theoretical and numerical analysis of non-Newtonian fluids, driven by a broad range of applications ranging from biological systems to industrial processes. In particular, many industrial and natural processes involve materials that exhibit non-Newtonian behavior, including food processing, polymer manufacturing, tribology, foam injection molding, and rubber extrusion.
First existence results concerning incompressible non-Newtonian fluid with
	$${\mathcal A}({\mathbb D}{\mathbf u})=\nu_1{\mathbb D}{\mathbf u}+\nu_2|{\mathbb D}{\mathbf u}|^{p-2}{\mathbb D}{\mathbf u},$$
were obtained for constant $p\geqslant \frac{11}{5}$ by Lions and Ladyzhenskaya in \cite{Ladyzhenskaya-1969,Lions-1969}.
Since then, the field has witnessed substantial progress. Bellout, Bloom and Ne\v{c}as \cite{Bellout-1994} proved the existence of Young-measure-valued solutions to the incompressible non-Newtonian fluid equations for space-periodic problems under certain conditions. M\'{a}lek, Ne\v{c}as, Rokyta, and R\r{u}\v{z}i\v{c}ka \cite{Malek-1997} proved the existence and uniqueness of global strong solutions for incompressible non-Newtonian fluids of power-law type when $p>\frac{3d-4}{d}.$ Kaplick\'{y}, M\'{a}lek and Star\'{a} \cite{P.kaplicky-2002} proved the existence of $C^{1,\alpha}-$solutions to the nonlinear system describing the plane stationary flow of a class of non-Newtonian fluids. Specifically, a global $C^{1,\alpha}-$solution was constructed for $p >\frac{3}{2},$ and  a solution with interior $C^{1,\alpha}$ was obtained for $p>\frac{6}{5}.$
Later, Zhang, Guo, and Guo \cite{Z-G-G2010} proved
the existence of weak solutions for steady-state incompressible non-Newtonian fluids in spatial dimension $d\geqslant2$ under slip boundary conditions when $p\geqslant2$. Burczak, Modena, and Sz\'{e}kelyhidi \cite{B-M-S2021} proved the non-uniqueness of Leray-Hopf solutions for
three-dimensional incompressible non-Newtonian fluids of power-law type when $1<p <\frac{6}{5}.$
	
Motivated by the model for the motion of electrorheological fluids introduced in \cite{Rajagopal-1996,Rrzicka-2000,Rajagopal-2001}, where
$p$ is a function of space and time. Electrorheological fluids are a special type of smart fluids which change their material properties in response to the application of an electric field. For the model in \cite{Rajagopal-2001}, $p$ is not a constant but a function of the electric field $\bf E$, i.e., $p = p(|{\bf E}|^2).$ The electric field itself is a solution to the quasi-static Maxwell equations and is not influenced by the motion of the fluid. R\r{u}\v{z}i\v{c}ka  \cite{Rrzicka-1999} proved the existence of a strong solution to three-dimensional incompressible non-Newtonian fluids of power-law type for large times and large initial data, provided that $p(t,x)\in W^{1,\infty}((0,T)\times \Omega).$ Further, R\r{u}\v{z}i\v{c}ka \cite{Rrzicka-2000} established the existence under the condition $\frac{11}{5}<p^-\leqslant p^+<p^-+\frac43$, where
\begin{eqnarray}\label{S-1}
  p^-=\inf\limits_{(0,T)\times \Omega}p(t,x)\mbox{  and  }p^+=\sup\limits_{(0,T)\times \Omega}p(t,x).
\end{eqnarray}
Diening and R\r{u}\v{z}i\v{c}ka \cite{D-R2005} proved the local existence of strong solutions in a three-dimensional bounded periodic domain when
when $p(t,x)$ satisfies $\frac{7}{5}<p^-\leqslant p(t,x)\leqslant p^+\leqslant2.$
	
Research on compressible non-Newtonian fluids with constant power-law exponents can be traced back to the early 1990s. For instance, Ne\v{c}asov\'{a}, and Novotn\'{y} \cite{Necasova-1993} proved the existence of the measure-valued solutions in. Later, Mamontov \cite{Mamontov-1999,Mamontov-1999-1} established the global existence of a weak solution under the assumptions of an exponentially growing viscosity and isothermal pressure. Recently, Abbatiello, Feireisl, and Novotn\'{y} \cite{Abbatiello-2020} proved the existence of a dissipative solution. For strong solutions, the local-in-time existence was established in the absence of vacuum by Kalousek, M\'{a}cha, and Ne\v{c}asov\'{a} \cite{Kalousek-Macha-Necasova-2020}. Fang and Guo \cite{Fang-Guo-2012} studied the analytical solution to a class of non-Newtonian fluids with free boundaries. Xu and Yuan \cite{Yuan-Xu-2008} proved the local-in-time existence and uniqueness in one-space dimension with singularity and vacuum. When the initial energy was sufficiently small, Yuan, Si, and Feng \cite{Yuan-Si-Feng-2019} established the global well-posedness of strong solutions for the initial-boundary-value problem of the one-dimensional model (\ref{XT-0}). Recently, Bilal Al Taki \cite{BAT2023} proved the existence and uniqueness of strong solutions for compressible non-Newtonian fluids of power-law type in the local-in-time sense when $1<p<\infty.$ Al Baba, Al Taki, and Hussein \cite{AL2024} proved the existence, uniqueness, and continuous dependence on the data of local-in-time strong solutions to the Navier-Stokes equations describing compressible non-Newtonian fluids when vacuum is present initially. Fang and Zang \cite{F-Z2023} proved the global existence and uniqueness of strong solutions to the Cauchy problem for one-dimensional compressible non-Newtonian fluids with $p\in (1,2).$ Recently, Baba, Taki and Hussein \cite{AL2024} investigated the blow-up criterion for three-dimensional compressible non-Newtonian fluids allowing initial vacuum states and established the following result
\begin{equation}
\label{bbpp}
\lim\limits_{T\rightarrow T^{*}}\sup\limits_{0\leqslant t\leqslant	T}\left(\|\rho\|_{W^{1,q_{0}}(\Omega)}+\|\uu\|_{H^{1}(\Omega)}\right)=\infty,
\end{equation}
where $q_{0}=\min{\{6,q\}}$ with $q\in(3,\infty).$ Guo and Fang established the following two blow-up criteria for strong solutions to a three-dimensional compressible non-Newtonian fluid with vacuum and constant power law
\begin{eqnarray*}\label{1dlg-E9}
	\lim\limits_{T\rightarrow T^{*}}\int_{0}^{T}\|\nabla \uu\|_{L^{\infty}(\Omega)}^{4}dt=\infty\mbox{ or }
	\lim\limits_{T\rightarrow T^{*}}\left(\|\rho\|_{L^{\infty}(0,T;L^{\infty}(\Omega))}+\|\nabla \uu\|_{L^{\infty}(0,T;L^{3}(\Omega))}\right)=\infty.
\end{eqnarray*}
For further results on the mathematical theory of non-Newtonian fluids, we refer the reader to \cite{Guo-Zhu-2002, hanguo, Guo-Lin-Shang-2006, Diening-2010, Zang-2018, Kalousek-Macha-Necasova-2020, F-G}.
	
However, to our knowledge, few results exist for the compressible non-Newtonian fluids of power-law type $p(t,x)$ with
$p(t,x)\in W^{1,\infty}((0,T)\times \Omega).$ Motivated by the work of Diening and R\r{u}\v{z}i\v{c}ka \cite{D-R2005}, we aim to establish the existence of strong solutions to the following system
\begin{equation}\label{XT-1}
	\begin{cases}
	 \partial_{t}\rho+\divv(\rho{\mathbf u})=0\quad \mbox{in}\quad (0,T)\times\Omega,\\
	 \partial_{t}(\rho{\mathbf u})+\divv(\rho{\mathbf u}\otimes{\mathbf u} )
       -\divv((1+|{\mathbb D}{\mathbf u}|^2)^{\frac{p(t,x)-2}{2}}{\mathbb D}{\mathbf u})
	  +\nabla \rho^\gamma=\rho{\mathbf f}\quad \mbox{in}\quad (0,T)\times\Omega
	\end{cases}
\end{equation}
subjected to the initial condition
\begin{equation}\label{XT-2}
	\rho(x, 0) = \rho_0(x),\quad \rho {\mathbf u}(x, 0)={\bf m}_0(x),\quad x\in\Omega,
\end{equation}
where $\Omega$ is the torus $\Omega={\mathbb T}^3.$
	
Local strong solutions to the problem \eqref{XT-1}-\eqref{XT-2} are defined as follows.
	
\begin{definition}\label{definition}
Let $p(t,x)\in W^{1,\infty}((0,T)\times \Omega)$ satisfy $\frac{7}{5}<p^-\leqslant p(t,x)\leqslant p^+\leqslant 2,$ and let $T\in (0,+\infty)$ be a given positive time. A pair $( \rho, {\mathbf u})$ is called a strong solution to the problem \eqref{XT-1}-\eqref{XT-2} on $(0,T)\times\Omega,$ if $(\rho,{\mathbf u})$  possesses the following properties:
\begin{eqnarray}\label{1dlg-E6}
	\begin{cases}
		\rho\in L^{\infty}(0,T; W^{1,3p^-}(\Omega)),\\
		{\mathbf u} \in L^\infty(0,T; W^{1,\frac{12}{5}}(\Omega)),~~
		{\mathbf u}\in L^{\frac{p^-(5p^--6)}{2-p^-}}(0,T; W^{2,\frac{3p^-}{p^-+1}}(\Omega)) \cap L^{p^-}(0, T; W^{2,3p^-}(\Omega)),\\
		\sqrt{\rho}\partial_t{\mathbf u}\in L^{\infty}(0, T; L^2(\Omega)),~~\partial_t{\mathbf u}\in L^{p^-}(0, T; W^{1,\frac{3p^-}{p^-+1}}(\Omega)).
		\end{cases}
\end{eqnarray}
\end{definition}
	
\begin{remark}
Thanks to the regularities of the strong solution stated in Definition \ref{definition} and the system \eqref{XT-1}, one can show that
$$\partial_t\rho\in L^{\infty}(0,T; L^{3p^-}(\Omega)).$$
\end{remark}
	
We are now in a position to state our first theorem.
	
\begin{theorem}\label{theo1}
Let $T>0,$ $\Omega = \mathbb{T} ^3$ and $p(t,x)\in W^{1,\infty}((0,T)\times \Omega)$ satisfy $\frac{7}{5}<p^-\leqslant p(t,x)\leqslant p^+\leqslant 2.$ Assume that the initial data $(\rho_0,{\mathbf u}_0)$ and the external force ${\mathbf f}$ satisfy the following regularity
\begin{align}\label{CZ}\nonumber
  &0\leqslant\rho_0\in W^{1,3p^-}(\Omega),\quad \sqrt{\rho_0}\uu_0\in L^2(\Omega),\quad\nabla\uu_0\in L^{\min \{6(p^--1),2p^-\}}(\Omega),\\
  &{\mathbf f}\in L^\infty(0,T;W^{1,2}(\Omega)),\quad \partial_t{\mathbf f}\in{L^2(0,T;L^2(\Omega))},
\end{align}
and the compatibility condition
\begin{equation}\label{condition}
  -\divv((1+|{\mathbb D}{\mathbf u}_0|^2)^{\frac{p(t, x)-2}{2}}{\mathbb D}{\mathbf u}_0)+\nabla (\rho_0)^\gamma = \rho_0^{1/2}{\bf g}
\end{equation}
for some ${\bf g}\in L^{3p^-}(\Omega).$ Then, there exist a small time $T^{* }\in ( 0, T)$ and a unique strong solution $(\rho,{\mathbf u})$ to the problem \eqref{XT-1}-\eqref{XT-2} with the compatibility conditions \eqref{condition}.
\end{theorem}
	
\begin{remark}\label{remark-1}
Basically, the condition $\frac{7}{5}<p^-\leqslant p(t,x) \leqslant p^+\leqslant 2$ covers the constant $p\in (\frac{7}{5},2].$ Indeed, Yuan and Xu  proved the local-in-time existence and uniqueness of strong solutions to the problem \eqref{XT-1}-\eqref{XT-2} in \cite{Yuan-Xu-2008} for $1<p<2$ in the case of
\begin{align*}
  {\mathcal A}({\mathbb D}{\mathbf u})=|{\mathbb D}{\mathbf u}|^{p-2}{\mathbb D}{\mathbf u},
\end{align*}
which is limited to the one-dimensional space.
\end{remark}
	
Next, we establish a blow-up criterion for the problem \eqref{XT-1}-\eqref{XT-2} under the compatibility conditions \eqref{condition}.
	
\begin{theorem}\label{blow}
Let $T>0,$ $\Omega = \mathbb{T} ^3$ and $p(t,x)\in W^{1,\infty}((0,T)\times \Omega)$ satisfy $\frac{7}{5}<p^-\leqslant p(t,x)\leqslant p^+\leqslant 2.$  Assume that $(\rho_0,\uu_0)$ and extend force $\mathbf{f}$ satisfy \eqref{CZ} and \eqref{condition}. Let $(\rho,{\mathbf u})$ be a strong solution to the problem \eqref{XT-1}-\eqref{XT-2} satisfying the regularity \eqref{1dlg-E6}. If $T^*\in(0,+\infty)$	is the maximal times of existence, then
\begin{equation}\label{blow-1}
  \lim\limits_{T\rightarrow T^*}(\|\rho\|_{{L^\infty}(0,T;L^\infty(\Omega))}+\|\nabla{\mathbf u}\|_{{L^\infty}(0,T;L^3(\Omega))})=\infty.
\end{equation}
\end{theorem}
	
\begin{remark}
The blow up criterion \eqref{blow-1} involves both the density and velocity. It may be natural to expect the higher regularity of velocity if the blow up criterion without the density. We can take similar argument for \eqref{blow-1} to arrive at the following blow up criterion	
\begin{equation}\label{blow-2}
  \lim\limits_{T\rightarrow T^*}\int_{0}^{T}\|\nabla{\mathbf u}\|_{L^\infty(\Omega)}^4~dt=\infty.
\end{equation}
The main difficulty in the blow up criterions \eqref{blow-1} and \eqref{blow-2} is to control the gradient of density, which is not known and coupled with the second derivative of velocity.
\end{remark}
	
\begin{remark}
The blow up criteria \eqref{blow-1} and \eqref{blow-2} hold for the case that the power $p(t,x) $ is a constant with $p\geqslant1$ (see \cite{Guo-Fang-2026}). Further, compared with the blow-up criterion \eqref{bbpp} in \cite{AL2024}, the blow-up criterion \eqref{blow-1} seems to impose weaker regularity condition on the density $\rho.$
\end{remark}
	
{\bf Notation}

The functions discussed in this section is from $\Omega\times {\mathbb R}^{d\times d}$ to ${\mathbb R}.$ We distinguish the partial derivatives by $\partial_i$ and $\partial_{jk}:$ a single index denotes a partial derivative with respect to the $i$-th space coordinate, while a double index represents a partial derivative with respect to the $(j,k)$-component of the underlying $d \times d$-matrix space. Since we are dealing with functions defined on $\Omega \times \mathbb{R}^{d \times d}$ to $\mathbb{R},$ we denote the partial derivatives by $\partial_i$ (spatial derivatives) and $\partial_{jk}$ (derivatives with respect to the matrix entries). By $\nabla$ we denote the space gradient, while $\nabla_{d \times d}$ denotes the matrix consisting of the partial derivatives with respect to the space of matrices. In a few cases we use $d_i$ instead of $\partial_i$ to denote a total derivative. We denote by ${\mathbf B}^{\text{sym}}$ the symmetric part of a matrix ${\mathbf B} \in \mathbb{R}^{d \times d}$, i.e., ${\mathbf B}^{\text{sym}} = \frac{1}{2}({\mathbf B} + {\mathbf B}^\top).$ ${\mathbf B}_{jk}$ denotes the $(j, k)$-component of the matrix ${\mathbf B},$ and $|{\mathbf B}|$ denotes $(\mathop{\sum}\limits_{j,k=1}^3|{\mathbf B}_{jk}|^2)^\frac12.$ Further let $\mathbb{R}^{\text{sym}}_{d \times d}$ be the subspace of $\mathbb{R}^{d \times d}$ consisting of the symmetric matrices. Moreover we use $C$ as a positive constant which is generic but does not depend on the ellipticity constants.  $\mathbb{R}$ is the real numbers and \(\mathbb{R}^+\) denotes the positive real numbers.

%%%%%%%%%%%%%%%%%%%%%%%%%%%%%%%%%%%%%%%%%%%%%%%%%%%%%%%%%%%%%%%%%%%%%%%%%%%%%%%%%%%%%%%%%%%%%%%%%

\section{Preliminaries}\label{Sec-P}
	
Following a computational concept from the work in \cite{Diening-2010}, we let $I=(0,T)$ with $T>0$ and let $p(t,x): I\times \Omega \to ( 1, 2]$ be a $W^{1, \infty }( I\times \Omega )$ function with
	$$1< p^{-}: = \inf\limits_{I\times \Omega} p(t,x)\leqslant p^{+} := \sup \limits_{I\times \Omega}p(t,x)\leqslant 2.$$
Motivated by \cite{D-R2005}, we define a $p(t,x)$-potential and derive its basic properties as follows. Let
	$$F : I\times \Omega \times ({\mathbb R}^+\cup\{0\})\to{\mathbb R}^+\cup\{0\}$$
and let $\Phi: \mathbb{R} ^{d\times d}\to {\mathbb R}^+\cup\{0\}$ be defined as
$$\Phi(t, x, {\mathbf B}) = F(t, x, | {\mathbf B} ^{\mathrm{sym}}|)$$
with
\begin{eqnarray}
	&&\sum_{jklm}(\partial_{jk}\partial_{lm}\Phi)(t,x,{\mathbf B})\mathbf{C}_{jk}\mathbf{C}_{lm}\geqslant
	\gamma_{1}(1+|{\mathbf B}^{\mathrm{sym}}|^{2})^{\frac{p(t,x)-2}{2}}|\mathbf{C}^{\mathrm{sym}}|^{2},\label{2.1}\\
	&&\left|(\nabla_{d\times d}^{2}\Phi)(t,x,{\mathbf B})\right|\leqslant
	\gamma_{2}(1+|{\mathbf B}^{\mathrm{sym}}|^{2})^{\frac{p(t,x)-2}{2}}\label{2.2}
\end{eqnarray}
for all ${\mathbf B}, \mathbf{C} \in \mathbb{R} ^{d\times d}$ with constants $\gamma _{1}, \gamma _{2}> 0.$
Further, we assume that $F$ is continuously differentiable with respect to $t$ and $x.$ Moreover, for all $t\in I$ and $x\in \Omega,$ we have the following:
\begin{itemize}
	\item [$\bullet$] $(\partial_{t}F)(t,x,\cdot):{\mathbb R}^+\cup\{0\}\rightarrow\:{\mathbb R}^+\cup\{0\}$ is a $C^{1}$-function on  $[0,+\infty);$
	\item [$\bullet$] $(\partial_{j}F)(t,x,\cdot):{\mathbb R}^+\cup\{0\}\:\rightarrow\:{\mathbb R}^+\cup\{0\}$ is a $C^{2}$-function on $(0,+\infty)$ and (and continuous up to $0$).
\end{itemize}
Moreover, we assume that for $j=1,2,3$
\begin{align}\label{D}\nonumber
	&(\partial_{t}F)(t,x,0)=0,\nonumber\\
	&(\partial_{t}F)(t,x,R)>0\mbox{ for all }\:R>0,\nonumber\\
	&(\partial_{j}F)(t,x,0)=0,\nonumber\\
	&(\partial_{j}F)(t,x,R)>0\mbox{ for all } R>0,\nonumber\\
	&\left|(\partial_t\nabla_{d\times d}\Phi)(t, x, {\mathbf B})\right|\leqslant
		\gamma_{3}(1+|{\mathbf B}^{\mathrm{sym}}|^{2})^{\frac{p(t, x)-1}{2}}\ln(1+|{\mathbf B}^{\mathrm{sym}}|),\nonumber\\
	&\left|(\nabla\nabla_{d\times d}\Phi)(t, x, {\mathbf B})\right|\leqslant
		\gamma_{3}(1+|{\mathbf B}^{\mathrm{sym}}|^{2})^{\frac{p(t, x)-1}{2}}\ln(1+|{\mathbf B}^{\mathrm{sym}}|),
\end{align}
where $\gamma_{3}>0$. Such function $\Phi$ is called a $p(t,x)$-potential. We define the extra stress $\mathbf{S}$ induced by $F$ and $\Phi $ as follows:
$$\mathbf{S}({\mathbf B}):=\nabla_{d\times d}\Phi({\mathbf B})
=F^{\prime}(|{\mathbf B}^{\mathrm{sym}}|)\frac{{\mathbf B}^{\mathrm{sym}}}{|{\mathbf B}^{\mathrm{sym}}|}$$
for all ${\mathbf B}$ in $\mathbb{R} ^{d\times d}\setminus \{ \mathbf{0} \}.$

\begin{remark}
The standard example of such ${\mathbf S}$ is
$${\mathbf S}({\mathbb D}{\uu}):={\mathcal A}({\mathbb D}{\uu})=(1+|{\mathbb D}{\mathbf u}|^2)^{\frac{p(t,x)-2}{2}}{\mathbb D}{\mathbf u},$$
where $p=p(t,x)\in W^{1,\infty}(I\times \Omega).$
\end{remark}
	
Following \cite{Malek-1997} and \cite{Diening-2010},  \eqref{2.1} and \eqref{2.2} imply that $\mathbf{S}$ has the following properties.
	
\begin{remark}(\cite{Diening-2010})
Let $d\geqslant2$, then
\begin{eqnarray*}
&&(\partial_{jk}\Phi)({\mathbf B})=F^{\prime}(|{\mathbf B}^{\mathrm{sym}}|)\frac{{\mathbf B}_{jk}^{\mathrm{sym}}}{|{\mathbf B}^{\mathrm{sym}}|},\\
&&(\partial_{jk}\partial_{lm}\Phi)({\mathbf B})
=F^{\prime}(|{\mathbf B}^{\mathrm{sym}}|)\left(\frac{\delta_{jk,lm}^{\mathrm{sym}}}{|{\mathbf B}^{\mathrm{sym}}|}-\frac{{\mathbf B}_{jk}^{\mathrm{sym}}{\mathbf B}_{lm}^{\mathrm{sym}}}{|{\mathbf B}^{\mathrm{sym}}|^{3}}\right)+F^{\prime\prime}(|{\mathbf B}^{\mathrm{sym}}|)\frac{{\mathbf B}_{jk}^{\mathrm{sym}}}{|{\mathbf B}^{\mathrm{sym}}|}\frac{{\mathbf B}_{lm}^{\mathrm{sym}}}{|{\mathbf B}^{\mathrm{sym}}|},\\
&&\gamma_1(1+|{\mathbf B}^{\mathrm{sym}}|^2)^{\frac{p(t,x)-2}{2}}\leqslant F''(|{\mathbf B}^{\mathrm{sym}}|)\leqslant\gamma_2(1+|{\mathbf B}^{\mathrm{sym}}|^2)^{\frac{p(t,x)-2}{2}},\\
&&\lim\limits_{|{\mathbf B}|\to0}\mathbf{S}_{jk}({\mathbf B})=\lim\limits_{|{\mathbf B}|\to0}(\partial_{jk}\Phi)({\mathbf B})=0
\end{eqnarray*}
hold for all ${\mathbf B}\in \mathbb{R} ^{d\times d}\setminus \{ \mathbf{0} \},$
where $\delta _{jk, lm}^\mathrm{sym}: = \frac 12( \delta _{jl}\delta _{km}+ \delta _{jm}\delta _{kl}).$
\end{remark}
	
\begin{remark}(\cite{Diening-2010})
Let ${\mathbf B}$,$\mathbf{C}\in~\mathbb{R}^{d\times d}$. Due to $\Phi({\mathbf B})= F(|{\mathbf B}^{\mathrm{sym}}|)$,
then $\Phi({\mathbf B})=\Phi({\mathbf B}^{\mathrm{sym}})$ and
\begin{eqnarray*}
&&\sum_{jklm}(\partial_{jk}\partial_{lm}\Phi)({\mathbf B})\mathbf{C}_{jk}\mathbf{C}_{lm}=\sum_{jklm}(\partial_{jk}\partial_{lm}\Phi)({\mathbf B}^{\mathrm{sym}})\mathbf{C}_{jk}^{\mathrm{sym}}\mathbf{C}_{lm}^{\mathrm{sym}},\\
&&(\nabla_{d\times d}\Phi)({\mathbf B})=(\nabla_{d\times d}\Phi)({\mathbf B}^{\mathrm{sym}}),\\
&&(\nabla_{d\times d}^{2}\Phi)({\mathbf B})=(\nabla_{d\times d}^{2}\Phi)({\mathbf B}^{\mathrm{sym}}).
\end{eqnarray*}
Since we will mostly deal with symmetric matrices later, we use ${\mathbf B}$ instead of ${\mathbf B}^{\mathrm{sym}}.$
\end{remark}

\begin{lemma}(\cite{Diening-2010})\label{lem1}
There exist constants $c_1,c_2>0$ independent of $\gamma_1$ and $\gamma_2$ such that
\begin{eqnarray}
  &&\mathbf{S}(\mathbf{0})=0,\label{E-14}\\
  &&\sum_{ij}(\mathbf{S}_{ij}({\mathbf B})-\mathbf{S}_{ij}(\mathbf{C}))({\mathbf B}_{ij}-\mathbf{C}_{ij})
	 \geqslant c_{1}\gamma_{1}(1+|{\mathbf B}|^{2}+|\mathbf{C}|^{2})^{\frac{p(t,x)-2}{2}}|\mathbf{B-C}|^{2},\label{E-15}\\
  &&\sum_{ij}\mathbf{S}_{ij}({\mathbf B}){\mathbf B}_{ij}\geqslant c_{1}\:
	  \gamma_{1}(1+|{\mathbf B}|^{2})^{\frac{p(t,x)-2}{2}}|{\mathbf B}|^{2},\label{E-16}\\
  &&|\mathbf{S}({\mathbf B})-\mathbf{S}(\mathbf{C})|\leqslant c_{2}\:
	  \gamma_{2}(1+|{\mathbf B}|^{2}+|\mathbf{C}|^{2})^{\frac{p(t,x)-2}{2}}\:|{\mathbf B}-\mathbf{C}|,\label{E-17}\\
  &&|\mathbf{S}({\mathbf B})|\leqslant c_{2}\:\gamma_{2}(1+|{\mathbf B}|^{2})^{\frac{p(t,x)-2}{2}}\:|{\mathbf B}|\label{E-18}
\end{eqnarray}
hold for all ${\mathbf B},~\mathbf{C} \in \mathbb{R} _{\mathrm{sym}}^{d\times d}.$
\end{lemma}
	
Let $\Phi$ be a $p(t,x)$-potential and let ${\mathbf u}$ denote a sufficiently smooth function over the space-time cylinder. The brackets $\langle\cdot,\cdot\rangle$ denote integration over the space domain $\Omega.$ The following important expressions will be used.
\begin{eqnarray}
	&&{\mathcal I}_{\Phi}(t,{\mathbf u}):=\Big\langle\sum_{r}\sum_{jk\alpha\beta}(\partial_{\alpha\beta}\partial_{jk}\Phi)
		\big({\mathbb D}{\mathbf u}\big)\partial_{r}({\mathbb D}{\mathbf u})_{\alpha\beta},
		\partial_{r}({\mathbb D}{\mathbf u})_{jk}\Big\rangle(t),\label{ID}\\
	&&{\mathcal J}_{\Phi}(t,{\mathbf u}):=\Big\langle\sum_{jk\alpha\beta}(\partial_{\alpha\beta}\partial_{jk}\Phi)
		\big({\mathbb D}{\mathbf u}\big)\partial_{t}({\mathbb D}{\mathbf u})_{\alpha\beta},
		\partial_{t}({\mathbb D}{\mathbf u})_{jk}\Big\rangle(t),\label{JD}\\
	&&{\mathcal G}_{\Phi}(t,\mathbf{w},\mathbf{v}):=\Big\langle\sum_{jk\alpha\beta}(\partial_{\alpha\beta}\partial_{jk}\Phi)({\mathbb D}{\mathbf w})
		({\mathbb D}{\mathbf v})_{\alpha\beta},({\mathbb D}{\mathbf v})_{jk}\Big\rangle(t),\label{GD}
\end{eqnarray}
where $\mathbf{w}:I\times\Omega\to\mathbb{R}^3$ and $\mathbf{v}:I\times\Omega\to\mathbb{R}^3$ (or $\mathbf{v}:I\times\Omega\to\mathbb{R}^{3\times 3})$ are sufficiently smooth functions. For the sake of simplicity, we write ${\mathcal I}_\Phi({\mathbf u}),~{\mathcal J}_\Phi({\mathbf u})\mathrm{~and~}{\mathcal G}_\Phi(\mathbf{w},\mathbf{v})$ instead of ${\mathcal I}_{\Phi}(t,{\mathbf u}),~{\mathcal J}_{\Phi}(t,{\mathbf u})$ and ${\mathcal G}_{\Phi}(t,\mathbf{w},\mathbf{v}).$ It is easy to find that
\begin{equation*}
	{\mathcal I}_{\Phi}({\mathbf u})={\mathcal G}_{\Phi}({\mathbf u},{\mathbb D}{\mathbf u})\mbox{ and }
	{\mathcal J}_{\Phi}({\mathbf u})={\mathcal G}_{\Phi}({\mathbf u},\partial_{t}{\mathbf u}).
\end{equation*}
Due to the properties of $\Phi,$ there exits some $\gamma_1>0$ such that
\begin{equation}\label{Gxz}
	{\mathcal G}_\Phi(\mathbf{w},\mathbf{v})\geqslant
	\gamma_1\int_{\Omega}(1+|{\mathbb D}{\mathbf w}|^2)^{\frac{p(t,x)-2}{2}}|{\mathbb D}{\mathbf v}|^2\:dx.
\end{equation}
Since $(1+|{\mathbb D}{\mathbf w}|^2)^{\frac12}$ appears frequently, it is essential to use the following shortcut
\begin{equation}
  \widetilde{D}\mathbf{w}:=(1+|{\mathbb D}{\mathbf w}|^2)^{\frac{1}{2}}.
\end{equation}
Consequently,
\begin{eqnarray}
&&{\mathcal I}_{\Phi}({\mathbf u})\geqslant
          	C\gamma_{1}\:\int_{\Omega}(\widetilde{D}{\mathbf u})^{p(t,x)-2}|\nabla{\mathbb D}{\mathbf u}|^{2}\:dx,\label{I}\\
&&{\mathcal J}_{\Phi}({\mathbf u})\geqslant
		    C\gamma_{1}\:\int_{\Omega}(\widetilde{D}{\mathbf u})^{p(t,x)-2}|\partial_{t}{\mathbb D}{\mathbf u}|^{2}\:dx.\label{J}
\end{eqnarray}
Note that
	$$|\nabla^2{\mathbf u}|\leqslant3|\nabla{\mathbb D}{\mathbf u}|\leqslant6|\nabla^2{\mathbf u}|,$$
which shows that$|\nabla{\mathbb D}{\mathbf u}|$ and $|\nabla^{2}{\mathbf u}|$ are equivalent up to a multiplicative constant.
	
Next, we present the following two lemmas. The first is a special form of the Gagliardo-Nirenberg interpolation inequality (see \cite{G-N}), and the second is a local version of Gronwall's lemma (see \cite{D-R2005}).
	
\begin{lemma}({Gagliardo-Nirenberg Inequality}\cite{G-N})\label{GN}
Let $\Omega\subset {\mathbb R}^d(d\geqslant2)$ be a bounded domain and let $j, k$ be positive integers. For $1\leqslant p,~r\leqslant \infty$ and $0\leqslant j< k,$ $j/k\leqslant\theta\leqslant1$, there exists some constant $C=C(d,k,p,r,j,\theta,\Omega)$ such that
\begin{equation*}
  \|\nabla^{j}{\mathbf u}\|_{L^q(\Omega)}\leqslant C\|\nabla^{k}{\mathbf u}\|_{L^p(\Omega)}^{\theta} \|{\mathbf u}\|_{L^r(\Omega)}^{1-\theta}
\end{equation*}
holds for any ${\mathbf u}\in W_0^{k,p}(\Omega),$ where $\frac{1}{q}=\frac{j}{d}+\theta(\frac{1}{p}-\frac{k}{d})+\frac{1-\theta}{r}.$
\end{lemma}
	
\begin{lemma}(Local Version of Gronwall's Lemma \cite{D-R2005})\label{local-gronwall}
Let $T, \alpha, c_0>0$ be given constants, and let $h\in L^1(0,T)$ with $h\geqslant 0$ almost every on $[0,T].$ Let $f\in C^1[0,T]$ satisfy $f\geqslant 0$ and
\begin{equation}\label{eq:local-gronwall-ineq}
  f'(t)\leqslant h(t)+c_0f(t)^{1 + \alpha}.
\end{equation}
Suppose $t_0\in [0,T]$ is such that $\alpha c_0 H(t_0)^{\alpha}t_0<1$, where
	$$H(t):=f(0)+\int_{0}^{t}h(s)~ds.$$
Then
	$$f(t)\leqslant H(t)\left(1-\alpha c_0 H(t)^{\alpha}t\right)^{-\frac{1}{\alpha}} $$
holds for all $t\in [0,t_0].$
\end{lemma}
	
Next, we introduce the following auxiliary result. 
	
\begin{lemma}( \cite{Feireisl-2004})\label{lem2-2}
Let ${\bf v}\in W^{1,2}(\Omega)$, and let $\rho$ be a non-negative function such that
$$0<M\leqslant\int_{\Omega}\rho ~dx, \ \ \int_{\Omega}\rho^{\gamma}~dx\leqslant  E_0,$$
where $\Omega\subset \mathbb{R}^{3}$ is a bounded domain, $M, E_0$ are positive constants, and $\gamma>1$.
		
Then there exists a constant $c,$ depending only on $M$ and $E_0$ such that
$$\|{\bf v}\|_{L^{2}(\Omega)}^{2}\leqslant  c(E_0,M)\big(\|\nabla {\bf v}\|_{L^{2}(\Omega)}^{2}+\big(\int_{\Omega}\rho|{\bf v} |~ dx\big)^{2}\big).$$
\end{lemma}
	
At the end of this section, we follow the idea of \cite{BAT2023} to establish some regularity properties of a solution to the following nonlinear elliptic system
\begin{align}
\begin{cases}\label{elliptic-ope-GENERAL}
  -\divv {\mathcal A}({\mathbb D}{\uu})={\mathbf f},\\
  {\mathcal A}({\mathbb D}{\uu})=(1+|{\mathbb D}{\mathbf u}|^2 )^{\frac{p(t,x)-2}{2}} {\mathbb D}{\mathbf u},
\end{cases}
\end{align}
where $p(t,x)\in W^{1, \infty }( (0,T)\times \Omega )$ satisfies $1< p^{-}\leqslant p(t,x)\leqslant p^{+} \leqslant 2.$ Define the quasi-linear differential operator ${\mathcal A}({\mathbf u},{\mathbb D})$ as
	$$\mathcal{A}({\mathbf u},\mathbb{D})=\sum_{k,l=1}A^{k,l}\:\mathbb{D}_k\:\mathbb{D}_l,$$
where the matrix-valued coefficients
	$$A^{k,l}({\mathbf u})=(a_{i,j}^{k,l}).$$
If ${\mathbf u}\in C^1(\overline{\Omega}),$ the coefficients of the differential operator $\mathcal{A}(x,\mathbb{D})=\mathcal{A}({\mathbf u}(x),\mathbb{ D}$) are uniformly bounded continuous. Freezing $\mathcal{A}$ at a reference solution ${\mathbf u}^{*}$, we obtain the linear operator
\begin{eqnarray}\label{E-A-U*}
  \mathcal{A}({\mathbf u}^*,\mathbb{D})=\sum_{k,l=1}A_*^{k,l}\:\mathbb{D}_k\:\mathbb{D}_l\:.
\end{eqnarray}
	
\begin{lemma}\label{lemzzx}
Let $p(t,x)\in W^{1, \infty }( (0,T)\times \Omega )$ satisfy $1< p^{-}\leqslant p(t,x)\leqslant p^{+} \leqslant 2,$ and let $\Omega$ be a periodic domain in ${\mathbb R}^3.$  Given a function ${\mathbf f}\in L^s(\Omega)~(1<s<\infty)$ such that
$$\int_\Omega {\mathbf f}\,dx=0.$$
Assume that ${\mathbf u}$ is a unique solution to the system \eqref{elliptic-ope-GENERAL}. If $d=2$ (or $3$), then the linearized operator associated to equation \eqref{elliptic-ope-GENERAL} at a reference solution ${\mathbf u}^*\in W^{2-\frac{2}{s},s}(\Omega)~(2s>5),$ denoted by $\mathcal{A}({\mathbf u}^*, \mathbb{D}),$  still yield maximal $L^s-$regularity. Moreover, for $s\in(\frac{5}{2},+\infty),$ there exists a constant $C>0$ such that
\begin{equation}\label{E-28}
	\| {\mathbf u} \|_{W^{2,s}(\Omega)}\leqslant\|{\mathcal A}({\mathbf u}^*, {\mathbb D}){\mathbf u}\|_{L^s (\Omega)}
	\leqslant C (\|{\mathbf f}\|_{L^s (\Omega)}+\|{\mathbb D}{\mathbf u}\|_{L^{p^-}(\Omega)}^{p^-}+1).
\end{equation}
\end{lemma}
	
\begin{proof}
It is important to prove that the linearized operator $\mathcal{A}({\mathbf u}^{*},\mathbb{D})$ satisfies the strong ellipticity condition. To start, we rewrite the elliptic operator as follows
\begin{align}\label{DIVS}
	&\divv{\mathcal A}({\mathbb D}{\uu})\nonumber\\
	&= \divv\left((1+|{\mathbb D}{\mathbf u}|^2)^{\frac{p(t,x)-2}2}{\mathbb D}{\mathbf u}\right)\nonumber\\
	&= \left(\frac{1}{2}(1+|{\mathbb D}{\mathbf u}|^2)^{\frac{p(t,x)-2}{2}}(\Delta{\mathbf u}+\nabla\divv{\mathbf u})\right)
	+\frac{p(t,x)-2}{2}(1+|{\mathbb D}{\mathbf u}|^2)^\frac{p(t,x)-4}{2}\nabla(|{\mathbb D}{\mathbf u}|^2){\mathbb D}{\mathbf u}\nonumber\\
	&~~~+\frac{1}{2}(1+|{\mathbb D}{\mathbf u}|^2)^\frac{p(t,x)-2}{2}\ln(1+|{\mathbb D}{\mathbf u}|^2)\nabla p(t,x){\mathbb D}{\mathbf u}\nonumber\\
	&\triangleq L_1+L_2+L_3.
\end{align}
Let ${\divv{\mathcal B}({\mathbb D}{\uu})}:=L_1+L_2$ and ${\mathbf h}:=L_3+{\mathbf f}$, then we can rewrite $\eqref{elliptic-ope-GENERAL}_1$ as
\begin{equation*}
	-{\divv{\mathcal B}({\mathbb D}{\uu})}= {\mathbf h}.
\end{equation*}
Since ${\mathbb D}{\mathbf u}$ is symmetric, the $i$-th entry of ${\divv{\mathcal B}({\mathbb D}{\uu})}$ becomes
\begin{align*}
	[{\divv{\mathcal B}({\mathbb D}{\uu})}]_i=
	&\sum_{k=1}^{3}\frac{1}{2}(1+|{\mathbb D}{\uu}|^2)^\frac{p(t,x)-2}{2}(\partial_{k}^{2}{\mathbf u}_{i}+\partial_{i}\partial_{k}{\mathbf u}_{k})\\
	&~~~+2\sum_{j,k,l=1}^{3}\frac{p(t,x)-2}{2}(1+|{\mathbb D}{\mathbf u}|^2)^\frac{p(t,x)-4}{2}
      ({\mathbb D}{\uu})_{ij}({\mathbb D}{\uu})_{kl}\:\partial_{j}({\mathbb D}{\uu})_{kl}\\
	&= \sum_{k=1}^{3}\frac{1}{2}(1+|{\mathbb D}{\uu}|^2)^\frac{p(t,x)-2}{2}
      (\partial_{k}^{2}{\mathbf u}_{i}+\partial_{i}\partial_{k}{\mathbf u}_{k}\bigr)\\
	&~~~+\sum_{j,k,l=1}^{3}({p(t,x)-2})(1+|{\mathbb D}{\mathbf u}|^2)^\frac{p(t,x)-4}{2}
      ({\mathbb D}{\uu})_{ik}({\mathbb D}{\uu})_{jl}\:\partial_{k}\partial_{l}{\mathbf u}_{j}\\
	&= \sum_{j,k,l=1}^{3}a_{ij}^{kl}\partial_{k}\partial_{l}{\mathbf u}_{j}.
\end{align*}
Set
$$a_{i,j}^{k,l}=\frac{1}{2}(1+|{\mathbb D}{\uu}|^2)^\frac{p(t,x)-2}{2}(\delta_{kl}\delta_{ij}
		+\delta_{il}\delta_{jk})+(p(t,x)-2)(1+|{\mathbb D}{\uu}|^2)^\frac{p(t,x)-4}{2}({\mathbb D}{\uu})_{ik}\:({\mathbb D}{\uu})_{jl},$$
where $\delta_{kl}$ denotes the Kronecker symbol. Recalling \eqref{E-A-U*}, one finds that the coefficients of the differential operator $\mathcal{A}({\mathbf u}^*,\mathbb{ D})=\mathcal{A}(x,\mathbb{D})$ are uniformly bounded continuous for any given ${\mathbf u}^*\in W^{2-\frac{2}{s},s}(\Omega)$ with $2s>d+2.$ Further, it is deduced from the definition of $a^{kl}_{jl}$ and the Cauchy-Schwarz inequality that
\begin{align}
	&(\mathcal{A}(x,\xi)\eta,\eta)\nonumber\\
  &= \sum_{i,j,k,1=1}^na_{ij}^{kl}\xi_k\eta_l\xi_i\bar{\eta}_j\nonumber\\
  &=\frac{1}{2}(1+|{\mathbb D}{\mathbf u}|^2)^\frac{p(t,x)-2}{2}(\sum_{i,k=1}^{n}\xi_{k}^{2}\eta_{i}\bar{\eta}_{i}
	+\sum_{i,k=1}^{n}\xi_{i}\eta_{k}\xi_{k}\bar{\eta}_{i})\nonumber\\
  &~~~+(p(t,x)-2)(1+|{\mathbb D}{\mathbf u}|^2)^\frac{p(t,x)-4}{2}\sum_{i,k=1}^{n}d_{ik} \xi_{k}\eta_{i}\sum_{j,l=1}^{n}d_{jl}\xi_{l}\bar{\eta}_{j} \nonumber\\
  &= \frac{1}{2}(1+|{\mathbb D}{\mathbf u}|^2)^\frac{p(t,x)-2}{2}(|\xi|^{2}|\eta|^{2}+(\xi,\eta)\overline{(\eta,\xi)})
    +(p(t,x)-2)(1+|{\mathbb D}{\mathbf u}|^2)^\frac{p(t,x)-4}{2}(\mathbb{D}\xi,\eta)\overline{(\mathbb{D}\xi,\eta)}\nonumber\\
  &= \frac{1}{2}(1+|{\mathbb D}{\mathbf u}|^2)^\frac{p(t,x)-2}{2}(|\xi|^{2}|\eta|^{2}+|(\xi,\eta)|^{2})
	+(p(t,x)-2)(1+|{\mathbb D}{\mathbf u}|^2)^\frac{p(t,x)-4}{2}|(\mathbb{D}\xi,\eta)|^{2}\nonumber\\
&\geqslant  \frac{1}{2}(1+|{\mathbb D}{\mathbf u}|^2)^\frac{p(t,x)-2}{2}(|\xi|^{2}|\eta|^{2}+|(\xi,\eta)|^{2})
	+(p(t,x)-2)(1+|{\mathbb D}{\mathbf u}|^2)^\frac{p(t,x)-4}{2}|{\mathbb D}{\mathbf u}|^2|(\xi,\eta)|^{2}\nonumber\\
&\geqslant  (1+|{\mathbb D}{\mathbf u}|^2)^\frac{p(t,x)-4}{2}(p(t,x)-1)(1+|{\mathbb D}{\mathbf u}|^2)
\end{align}
holds for any $\xi\in\mathbb{R}^d$ and $\eta\in\mathbb{C}^d$ with $|\xi|=|\eta|=1.$ Obviously, the operator $\mathcal{A}({\mathbf u}^*,\mathbb{D})$ is strongly elliptic. Hence, we have
\begin{equation}\label{zzx1}
	\| {\mathbf u} \|_{W^{2,s}(\Omega)}
	\leqslant \|\mathcal{A}({\mathbf u}^*, \mathbb{D}) {\mathbf u}\|_{L^s (\Omega)}\leqslant C \|{\mathbf h}\|_{L^s (\Omega)}.
\end{equation}
		
Note that
\begin{eqnarray*}
	\|{\mathbf h}\|_{L^s(\Omega)}
     &&=\frac{1}{2}\|(1+|{\mathbb D}{\mathbf u}|^2)^\frac{p(t,x)-2}{2}\ln(1+|{\mathbb D}{\mathbf u}|^2)\nabla p(t,x){\mathbb D}{\mathbf u}+{\mathbf f}\|_{L^s (\Omega)}\\
	&&\leqslant \frac{1}{2}\|(1+|{\mathbb D}{\mathbf u}|^2)^\frac{p(t,x)-2}{2}\ln(1+|{\mathbb D}{\mathbf u}|^2)\nabla p(t,x){\mathbb D}{\mathbf u}\|_{L^s (\Omega)}+\|{\mathbf f}\|_{L^s (\Omega)}.
\end{eqnarray*}
Moreover, it follows from Lemma \ref{GN} that
\begin{align}\label{zzx2}
	&\|(1 + |{\mathbb D}{\mathbf u}|^2)^{\frac{p(t,x)-2}{2}}\ln(1 + |{\mathbb D}{\mathbf u}|^2)\nabla p(t,x)
			{\mathbb D}{\mathbf u}\|_{L^s(\Omega)}\nonumber\\
	&\leqslant C\|(1 + |{\mathbb D}{\mathbf u}|^2)^{\frac{p(t,x)-1}{2}}\ln(1 + |{\mathbb D}{\mathbf u}|^2)\|_{L^s(\Omega)}\nonumber\\
	&\leqslant C\|(1 + |{\mathbb D}{\mathbf u}|^2)^{\frac{p(t,x)-1+\varepsilon}{2}}\|_{L^s(\Omega)}\nonumber\\
	&\leqslant C(1 + \|{\mathbb D}{\mathbf u}\|_{L^{s(p^+-1+\varepsilon)}(\Omega)}^{(p^+-1+\varepsilon)})\nonumber\\
	&\leqslant  C(1 + \|{\mathbb D}{\mathbf u}\|_{L^{p^-}(\Omega)}^{(1-\theta)(p^+-1+\varepsilon)}
	\|\nabla^2{\mathbf u}\|_{L^{s}(\Omega)}^{\theta(p^+-1+\varepsilon)})\nonumber\\
	&\leqslant  \delta\|\nabla^2{\mathbf u}\|_{L^s(\Omega)}
     +C(\delta)(1+\|{\mathbb D}{\mathbf u}\|_{L^{p^-}}^\frac{(1-\theta) (p^+-1+\varepsilon)} {1-\theta(p^+-1+\varepsilon)})
\end{align}
for any fixed $\varepsilon\in(0,1),$ where
$$\frac{1}{(p^+-1+\varepsilon)s}=\theta(\frac{1}{s}-\frac{1}{3})+\frac{1-\theta}{p^-}\mbox{ and } \theta=\frac{3(s(p^+-1+\varepsilon)-p^-)}{(p^-s+3s-3p^-)(p^+-1+\varepsilon)}.$$
Thus, one can choose sufficiently small $\varepsilon\in(0,1)$ such that
$$\frac{(1-\theta)(p^+-1+\varepsilon)}{1-\theta(p^+-1+\varepsilon)}\leqslant p^-.$$
So, it follows from \eqref{zzx1} and \eqref{zzx2} that
$$\| {\mathbf u} \|_{W^{2,s}(\Omega)}\leqslant\|{\mathcal A}({\mathbf u}^*, {\mathbb D}){\mathbf u}\|_{L^s (\Omega)}
		\leqslant C (\|{\mathbf f}\|_{L^s (\Omega)}+\|{\mathbb D}{\mathbf u}\|_{L^{p^-}(\Omega)}^{p^-}+1).$$
\end{proof}
	
\section{A priori estimates }\label{Sec-A}
	
First, we establish the following lemma, which plays a key role in the regularity analysis of solutions.
	
\begin{lemma}\label{GIJ}
Let $T>0,$ $\Phi$ be a $p(t,x)$-potential and $1< q\leqslant2.$ Then for all sufficiently smooth ${\mathbf u}$ and $\mathbf{v},$
\begin{equation}
 \|{\mathbb D}{\mathbf v}\|_{L^q(\Omega)}
 \leqslant C{\mathcal G}^\frac{1}{2}_\Phi(\mathbf{w},\mathbf{v})\|(\widetilde{D}\mathbf{w})^\frac{2-p^-}{2}\|_{L^\frac{2q}{2-q}(\Omega)},
\end{equation}
holds for almost every $t\in (0,T),$ where $\frac{2q}{2-q}=\infty$ for $q=2.$
\end{lemma}
	
\begin{proof}
Observe that $1\leqslant\frac2q<\infty$ and $1<\frac2{2-q}\leqslant\infty.$ For $1\leqslant q<2,$ one finds that
\begin{eqnarray*}
  &\|{\mathbb D}{\mathbf v}\|_{L^q(\Omega)}^{q}
  &=\int_{\Omega}\left((\widetilde{D}\mathbf{w})^{p^--2}|{\mathbb D}{\mathbf v}|^{2}\right)^{\frac{q}{2}}
     (\widetilde{D}\mathbf{w})^{\frac{(2-p^-)q}{2}}~dx \\
   &&\leqslant\left(\int_{\Omega}(\widetilde{D}\mathbf{w})^{p^--2}|{\mathbb D}{\mathbf v}|^{2}dx\right)^{\frac{q}{2}}
			\left\|(\widetilde{D}\mathbf{w})^{\frac{(2-p^-)q}{2}}\right\|_{L^{\frac{q}{2-q}}(\Omega)}\\
   &&=\left(\int_{\Omega}(\widetilde{D}\mathbf{w})^{p^--2}|{\mathbb D}{\mathbf v}|^{2}dx\right)^{\frac{q}{2}}
			\left\|(\widetilde{D}\mathbf{w})^{\frac{2-p^-}{2}}\right\|_{L^{\frac{2q}{2-q}}(\Omega)}^{q}.
\end{eqnarray*}
Taking \eqref{Gxz} into account, one prove the lemma for $1<q<2.$ The case $q=2$ is similar.
\end{proof}
	
Applying Lemma \ref{GIJ} to ${\cal I}_\Phi(\uu)={\cal G}_\Phi(\uu, \nabla\uu)$ and ${\cal J}_\Phi(\uu)={\cal G}_\Phi(\uu, {\uu}_t),$ one deduces the following Lemma.
	
\begin{lemma}\label{L9}
Let $T>0,$ $\Phi$ be a $p(t,x)$-potential and $1< q\leqslant2.$ Then for all sufficiently smooth ${\mathbf u}$ and $\mathbf{v},$
\begin{eqnarray*}
 &&\left(\int_\Omega((1+|{\mathbb D}{\mathbf u}|^2)^{\frac{p(t,x)-2}{2}}{\mathbb D}{\mathbf u}
	-(1+|{\mathbb D}{\mathbf v}|^2)^{\frac{p(t,x)-2}{2}}{\mathbb D}{\mathbf v}):{\mathbb D}({\mathbf u}-{\mathbf v}))~dx\right) ^\frac{1}{2}\\
 &&\geqslant C\|(\widetilde{D}{\mathbf u})^{\frac{2-p(x,t)}{2}}
    +(\widetilde{D}\bar{\mathbf u})^{\frac{2-p(x,t)}{2}}\|_{L^{\frac{2q}{2-q}}(\Omega)}^{-1}\|{\mathbb D}({\mathbf u}-{\mathbf v})\|_{L^q(\Omega)}
\end{eqnarray*}
holds for almost every $t\in (0,T),$ where $\frac{2q}{2-q}=\infty$ for $q=2.$
\end{lemma}
	
\begin{proof}
Analogously to the proof of Lemma \ref{GIJ}, one uses Lemma \ref{lem1} to get that 		
\begin{align*}
	&\|{\mathbb D}({\mathbf u}-{\mathbf v})\|_{L^q(\Omega)}^{q}\\
   &=\int_{\Omega}\left((\widetilde{D}\mathbf{u}+\widetilde{D}\mathbf{v})^{p(x,t)-2}
	 |{\mathbb D}({\mathbf u}-{\mathbf v})|^2\right)^{\frac{q}{2}}(\widetilde{D}\mathbf{u}+\widetilde{D}\mathbf{v})^{\frac{(2-p(x,t))q}{2}}~dx\\
&\leqslant \int_{\Omega}\left((\widetilde{D}\mathbf{u}+\widetilde{D}\mathbf{v})^{p(x,t)-2}
     |{\mathbb D}({\mathbf u}-{\mathbf v})|^2~dx\right)^{\frac{q}{2}}
	 \left\|(\widetilde{D}\mathbf{u}+\widetilde{D}\mathbf{v})^{\frac{(2-p(x,t))q}{2}}\right\|_{L^{\frac{2}{2-q}}(\Omega)}\\
&\leqslant \int_{\Omega}\left((\widetilde{D}\mathbf{u}+\widetilde{D}\mathbf{v})^{p(x,t)-2}
	 |{\mathbb D}({\mathbf u}-{\mathbf v})|^2~dx\right)^{\frac{q}{2}}
     \left\|(\widetilde{D}\mathbf{u}+\widetilde{D}\mathbf{v})^{\frac{2-p(x,t)}{2}}\right\|^q_{L^{\frac{2}{2-q}}(\Omega)}\\
&\leqslant  C\left(\int_\Omega((1+|{\mathbb D}{\mathbf u}|^2)^{\frac{p(t,x)-2}{2}}{\mathbb D}{\mathbf u}
  -(1+|{\mathbb D}{\mathbf v}|^2)^{\frac{p(t,x)-2}{2}}{\mathbb D}{\mathbf v}):{\mathbb D}({\mathbf u}-{\mathbf v}))~dx\right)^{\frac{q}{2}}\nonumber\\
  &\cdot\left\|(\widetilde{D}\mathbf{u})^{\frac{2-p(x,t)}{2}}+(\widetilde{D}\mathbf{v})^{\frac{2-p(x,t)}{2}}\right\|^q_{L^{\frac{2}{2-q}}(\Omega)}.
\end{align*}
This proves the lemma for $1<q<2.$ The case $q=2$ is similar.
\end{proof}
	
Next, we derive some estimates for ${\mathcal I}_\Phi(\uu)$ and ${\mathcal J}_\Phi(\uu).$
	
\begin{lemma}\label{ine}
Let $T>0.$	For all sufficiently smooth ${\mathbf u},$
\begin{eqnarray}
   &&\|{\mathbf u}\|_{W^{2,\frac{3p^-}{p^-+1}}(\Omega)}^{p^-}\leqslant C\left({\mathcal I}_{\Phi}({\mathbf u})+1\right),\label{III}\\
   &&\|\partial_{t}{\mathbf u}\|_{W^{1,\frac{3p^-}{p^-+1}}(\Omega)}^{p^-}
		\leqslant C{\mathcal J}^{\frac{p^-}{2}}_{\Phi}({\mathbf u})\left({\mathcal I}_{\Phi}({\mathbf u})+1\right)^{\frac{2-p^-}{2}}
		\leqslant C({\mathcal J}_{\Phi}({\mathbf u})+{\mathcal I}_{\Phi}({\mathbf u})+1)\label{IV}
\end{eqnarray}
hold for almost every $t\in (0,T).$
\end{lemma}
	
\begin{proof}	
According to Lemma \ref{GIJ}, one finds that
\begin{eqnarray*}
	&\|\nabla{\mathbb D}{\mathbf u}\|_{L^{\frac{3p^-}{p^-+1}}(\Omega)}
    &\leqslant C\:{\mathcal I}^{\frac{1}{2}}_{\Phi}({\mathbf u})\|(\widetilde{D}{\mathbf u})^{\frac{2-p^-}{2}}\|_{L^{\frac{6p^-}{2-p^-}}(\Omega)}\\
  &&\leqslant C\:{\mathcal I}^{\frac{1}{2}}_{\Phi}({\mathbf u})\|\widetilde{D}{\mathbf u}\|_{L^{3p^-}(\Omega)}^{\frac{2-p^-}{2}}\\
  &&\leqslant C\:{\mathcal I}^{\frac{1}{2}}_{\Phi}({\mathbf u})\left(1+\|{\mathbb D}{\mathbf u}\|_{L^{3p^-}(\Omega)}\right)^{\frac{2-p^-}{2}}\\
  &&\leqslant C\:{\mathcal I}^{\frac{1}{2}}_{\Phi}({\mathbf u})
     \left(1+C\left\|\nabla{\mathbb D}{\mathbf u}\right\|_{L^{\frac{3p^-}{p^-+1}}(\Omega)}\right)^{\frac{2-p^-}{2}},
\end{eqnarray*}
which implies that
$$\|\nabla{\mathbb D}{\mathbf u}\|_{L^{\frac{3p^-}{p^-+1}}(\Omega)}^{p^-}\leqslant C\left({\mathcal I}_\Phi({\mathbf u})+1\right).$$
Since $|\nabla^{2}{\mathbf u}|\leqslant3\left|{\nabla}{\mathbb D}{\mathbf u}\right|$, one deduces from Sobolev embedding and Korn's inequality that
$$\|{\mathbf u}\|_{W^{2,\frac{3p^-}{p^-+1}}(\Omega)}^{p^-}\leqslant C\left({\mathcal I}_{\Phi}({\mathbf u})+1\right).$$
Analogously, it is deduced from Lemma \ref{GIJ} to get that
\begin{eqnarray*}
	&\|\partial_{t}{\mathbb D}{\mathbf u}\|_{L^{\frac{3p^-}{p^-+1}}(\Omega)}
	&\leqslant C\:{\mathcal J}^{\frac{1}{2}}_{\Phi}({\mathbf u})\|(\widetilde{D}{\mathbf u})^{\frac{2-p^-}{2}}\|_{L^{\frac{6p^-}{2-p^-}}(\Omega)}\\
	&&\leqslant C\:{\mathcal J}^{\frac{1}{2}}_{\Phi}({\mathbf u})
         \left(1+\|\nabla{\mathbb D}{\mathbf u}\|_{L^{\frac{3p^-}{p^-+1}}(\Omega)}\right)^{\frac{2-p^-}{2}}\\
	&&{\leqslant}C\:{\mathcal J}^{\frac{1}{2}}_{\Phi}({\mathbf u})
        \big(1+\:({\mathcal I}_{\Phi}({\mathbf u})+1)^{\frac{1}{p^-}}\big)^{\frac{2-p^-}{2}}\\
	&&\leqslant C{\mathcal J}^{\frac{1}{2}}_{\Phi}({\mathbf u})\left(1+{\mathcal I}_{\Phi}({\mathbf u})\right)^{\frac{2-p^-}{2p^-}}.
\end{eqnarray*}
According to Korn's inequality, one obtains that
$$\left\|\partial_{t}{\mathbf u}\right\|_{W^{1,\frac{3p^-}{p^-+1}}(\Omega)}
		\leqslant C\left\|\partial_{t}{\mathbb D}{\mathbf u}\right\|_{L^{\frac{3p^-}{p^-+1}}(\Omega)}
		\leqslant C{\mathcal J}^{\frac{1}{2}}_{\Phi}({\mathbf u})\left(1+{\mathcal I}_{\Phi}({\mathbf u})\right)^{\frac{2-p^-}{2p^-}}.$$
Thus, \eqref{IV} follows from Young's inequality.
\end{proof}
	
Next, we establish the lower-order estimates for $(\rho,{\mathbf u}).$
	
\begin{lemma}
Let $(\rho,{\mathbf u})$ be a smooth solution to the problem \eqref{XT-1}-\eqref{XT-2} on $[0,T)\times\Omega$ for some $T>0.$ Then there exists a positive constant $C>0$ such that
\begin{eqnarray}\label{enery}
	&&\frac{d}{dt}\int_{\Omega}(\rho|{\mathbf u}|^{2}+\frac{1}{\gamma-1}\rho^{\gamma})\:dx
	+\int_{\Omega}(1+|{\mathbb D}{\mathbf u}|^2)^{\frac{p(t,x)-2}{2}}|{\mathbb D}{\mathbf u}|^{2}\:dx\nonumber\\
	&&\leqslant \|\sqrt{\rho}{\mathbf u}\|_{L^{2}(\Omega)}^2+C\|{\mathbf f}\|_{L^{2}(\Omega)}^2\|\rho\|_{W^{1,3p^-}(\Omega)}
\end{eqnarray}
for almost every $t\in [0,T).$
\end{lemma}
	
\begin{proof}
It is easy to deduce from \eqref{XT-1}$_1$ that
$$\int_\Omega\rho(t)\:dx=\int_\Omega\rho_0\:dx$$
for almost every $t\in [0,T).$ Next, multiplying \eqref{XT-1}$_2$ by ${\mathbf u}$ and integrating by parts over $\Omega,$ one obtains that
$$\frac{1}{2}\frac{d}{dt}\int_{\Omega}(\rho|{\mathbf u}|^{2}+\frac{1}{\gamma-1}\rho^{\gamma})\:dx
+\int_{\Omega}(1+|{\mathbb D}{\mathbf u}|^2)^{\frac{p(t,x)-2}{2}}|{\mathbb D}{\mathbf u}|^{2}\:dx
=\int_{\Omega}\rho {\mathbf f}\cdot{\mathbf u}\:dx.
$$
Moreover, it follows from H\"{o}lder inequality that
\begin{eqnarray*}
	&\left|\int_{\Omega}\rho {\mathbf f}\cdot{\mathbf u}\:dx\right|
	&\leqslant\|{\mathbf f}\|_{L^{2}(\Omega)}\|\rho\|_{L^\infty(\Omega)}^\frac{1}{2}\|\sqrt{\rho}{\mathbf u}\|_{L^{2}(\Omega)}\\
	&&\leqslant C\|{\mathbf f}\|_{L^{2}(\Omega)}\|\rho\|_{W^{1,3p^-}(\Omega)}^\frac{1}{2}\|\sqrt{\rho}{\mathbf u}\|_{L^{2}(\Omega)}\\
	&&\leqslant \frac{1}{2}\|\sqrt{\rho}{\mathbf u}\|_{L^{2}(\Omega)}^2+C\|{\mathbf f}\|_{L^{2}(\Omega)}^2\|\rho\|_{W^{1,3p^-}(\Omega)}.
\end{eqnarray*}
One applies \eqref{CZ} and Gronwall's inequality to deduce that
\begin{eqnarray*}
	&&\frac{d}{dt}\int_{\Omega}(\rho|{\mathbf u}|^{2}+\frac{1}{\gamma-1}\rho^{\gamma})\:dx
		+\int_{\Omega}(1+|{\mathbb D}{\mathbf u}|^2)^{\frac{p(t,x)-2}{2}}|{\mathbb D}{\mathbf u}|^{2}\:dx\\
    &&\leqslant \|\sqrt{\rho}{\mathbf u}\|_{L^{2}(\Omega)}^2+C\|{\mathbf f}\|_{L^{2}(\Omega)}^2\|\rho\|_{W^{1,3p^-}(\Omega)}
\end{eqnarray*}
for almost every $t\in [0,T).$
\end{proof}
	
Further, we establish the estimates for velocity.
	
\begin{lemma}\label{LEM1}
Let $\frac{12}{5}\leqslant q\leqslant \min \{6(p^--1),2p^-\}.$ Suppose that $(\rho,{\mathbf u})$ is a smooth solution to the problem \eqref{XT-1}-\eqref{XT-2} on $[0,T)\times\Omega$ for some $T>0.$ Then there exist positive constants $m_1=m_1(p^-),~m_2=m_2(p^-),~m_3=m_3(p^-)$ and $m_4=m_4(p^-)$ such that
\begin{align}\label{ans1}
	&\frac{d}{dt}\|\widetilde{D}{\mathbf u}\|_{L^q(\Omega)}^{q}+{\mathcal I}^r_\Phi({\mathbf u})\nonumber\\
	&\leqslant \eta{\mathcal J}_\Phi({\mathbf u})
	  +C(1+\|\nabla\rho\|_{L^{3p^-}(\Omega)}^{m_1}+\|\widetilde{D}{\mathbf u}\|_{L^q(\Omega)}^{m_2}+\|\rho\|_{L^\infty(\Omega)}^{m_3} +\|\sqrt{\rho}\partial_t{\mathbf u}\|_{L^2(\Omega)}^{m_4})
\end{align}
holds for any fixed $\eta\in(0,1),$ where $1\leqslant r<\frac{5p^--6}{2-p^-}.$
\end{lemma}
	
\begin{proof}
Multiplying $\eqref{XT-1}_{2}$ by $-\Delta{\mathbf u}$ and integrating the resulting equation over $\Omega,$ one obtains that
\begin{align}\label{delu}\nonumber
  &\int_\Omega\divv{\mathcal A}({\mathbb D}{\uu})\cdot\Delta{\mathbf u}~dx\nonumber\\
  &=\int_\Omega\divv\left((1+|{\mathbb D}{\mathbf u}|^2)^{\frac{p(t,x)-2}2}{\mathbb D}{\mathbf u}\right)\cdot\Delta{\mathbf u}~dx\nonumber\\
  &=\int_\Omega\rho{\mathbf u}\cdot\nabla{\mathbf u}\cdot\Delta{\mathbf u}~dx+\int_\Omega\nabla \rho^\gamma\cdot\Delta{\mathbf u}~dx
	+\int_{\Omega}\rho\partial_t{\mathbf u}\cdot\Delta{\mathbf u}~dx-\int_\Omega\rho {\mathbf f} \cdot\Delta{\mathbf u}~dx\nonumber\\\nonumber
  &\triangleq\sum_{i=1}^{4}J_i.
\end{align}
The first term on the right-hand side of equation \eqref{delu} is estimated as follows
\begin{align}
J_1&=\sum_{i,j,k}\int_\Omega\rho{\mathbf u}^j\partial_j{\mathbf u}^i\partial_k\partial_k{\mathbf u}^i~ dx\nonumber \\
   &=-\sum_{i,j,k}\int_\Omega(\partial_k\rho{\mathbf u}^j\partial_j{\mathbf u}^i\partial_k{\mathbf u}^i
	 +\rho\partial_k{\mathbf u}^j\partial_j{\mathbf u}^i\partial_k{\mathbf u}^i
	 +\rho {\mathbf u}^j\partial_k\partial_j{\mathbf u}^i\partial_k{\mathbf u}^i)~dx\nonumber\\
   &\leqslant\int_\Omega|\nabla\rho\|{\mathbf u}\|\nabla{\mathbf u}|^2~dx+\int_\Omega\rho|\nabla{\mathbf u}|^3~dx
	 +\int_\Omega\rho|{\mathbf u}\|\nabla^2{\mathbf u}\|\nabla{\mathbf u}|~dx\nonumber\\
   &\triangleq J_{11}+J_{12}+J_{12}.
\end{align}
Note that
\begin{align}\label{3-3}
	\|\nabla{\mathbf u}\|_{L^3(\Omega)}^3
		&\leqslant\|\nabla{\mathbf u}\|_{L^q(\Omega)}^\frac{3q(p^--1)}{3p^--q}
			\|\nabla{\mathbf u}\|_{L^{3p^-}(\Omega)}^\frac{3p^-(3-q)}{3p^--q}\nonumber\\
		&\leqslant C \|\nabla{\mathbf u}\|_{L^q(\Omega)}^\frac{3q(p^--1)}{2q+3p^--q}
			+\varepsilon\|\nabla{\mathbf u}\|_{L^{3p^-}(\Omega)}^{p^-}\nonumber\\
		&\leqslant C \|\nabla{\mathbf u}\|_{L^q(\Omega)}^\frac{3q(p^--1)}{2q+3p^--q}
			+\varepsilon\|\nabla{\mathbf u}\|_{W^{1,\frac{3p^-}{p^-+1}}(\Omega)}^{p^-}\nonumber\\
		&\leqslant C \|\nabla{\mathbf u}\|_{L^q(\Omega)}^\frac{3q(p^--1)}{2q+3p^--q}
			+\varepsilon({\mathcal I}_\Phi({\mathbf u})+1)\nonumber\\
		&\leqslant C \|\widetilde{D}{\mathbf u}\|_{L^q(\Omega)}^\frac{3q(p^--1)}{2q+3p^--q}+ \varepsilon({\mathcal I}_\Phi({\mathbf u})+1)
\end{align}
holds for any fixed $\varepsilon\in(0,1),$ when $\frac{12}{5}\leqslant q<3;$ the estimate \eqref{3-3} is also true when
$3\leqslant q \leqslant \min \{6(p^--1),2p^-\}.$  Using Lemma \ref{ine} and Sobolev embedding inequality, one gets that
\begin{align}\label{J11}
J_{11}&\leqslant C\|\nabla\rho\|_{L^{3p^-}(\Omega)}\|{\mathbf u}\|_{L^{12}(\Omega)}\|\nabla{\mathbf u}\|_{L^3(\Omega)}^2\nonumber\\
	  &\leqslant C(\|\nabla\rho\|_{L^{3p^-}(\Omega)}^6+\|\nabla{\mathbf u}\|_{L^\frac{12}{5}(\Omega)}^6
		+\|\nabla{\mathbf u}\|_{L^3(\Omega)}^3)\nonumber\\
	  &\leqslant\varepsilon{\mathcal I}_\Phi+ C(1+\|\nabla\rho\|_{L^{3p^-}(\Omega)}^6
		+\|\nabla{\mathbf u}\|_{L^{\frac{12}{5}}(\Omega)}^6+C\|\widetilde{D}{\mathbf u}\|_{L^q(\Omega)}^{R_1})
\end{align}
holds for any $\varepsilon\in(0,1),$ where $R_1= \frac{3q(p^--1)}{2q+3p^--9}.$
For $J_{12}$ and $J_{13},$ one uses Lemma \ref{GN} and Lemma \ref{ine} to get that
\begin{align}
J_{12}&\leqslant\|\rho\|_{L^\infty(\Omega)}\|\nabla{\mathbf u}\|_{L^3(\Omega)}^3 \nonumber\\
	  &\leqslant\|\rho\|_{L^\infty(\Omega)}\|\nabla{\mathbf u}\|_{L^\frac{12}{5}(\Omega)}^\frac{12(p^--1)}{5p^--4}
			\|{\mathbf u}\|_{L^{3p^-}(\Omega)}^\frac{3p^-}{(5p^--4)} \nonumber\\
	  &\leqslant\|\rho\|_{L^\infty(\Omega)}\|\nabla{\mathbf u}\|_{L^\frac{12}{5}(\Omega)}^\frac{12(p^--1)}{5p^--4}
			\|{\mathbf u}\|_{W^{2,\frac{3p^-}{3p^-+1}}(\Omega)}^\frac{3p^-}{(5p^--4)}\nonumber \\
	  &\leqslant\|\rho\|_{L^\infty(\Omega)}\|\nabla{\mathbf u}\|_{L^\frac{12}{5}(\Omega)}^\frac{12(p^--1)}{5p^--4}
			({\mathcal I}_\Phi({\mathbf u})+1)^\frac{3}{5p^--4} \nonumber\\
	  &\leqslant\varepsilon{\mathcal I}_\Phi({\mathbf u})
           +C(1+\|\rho\|_{L^\infty(\Omega)}^{R_2}+\|\nabla{\mathbf u}\|_{L^{\frac{12}{5}}(\Omega)}^{R_3}) \label{J12}
\end{align}
and
\begin{align}
J_{13}&\leqslant \|\rho\|_{L^\infty(\Omega)}\|\nabla^2{\mathbf u}\|_{L^\frac{3p^-}{p^-+1}(\Omega)}
		 \|\nabla{\mathbf u}\|_{L^\frac{12}{5}(\Omega)}\|{\mathbf u}\|_{L^{84}(\Omega)} \nonumber\\
	  &\leqslant\|\rho\|_{L^{\infty}(\Omega)}\|{\mathbf u}\|_{W^{2,\frac{3p^-}{p^-+1}}(\Omega)}
		\|\nabla{\mathbf u}\|_{L^{\frac{12}{5}}(\Omega)}\|\nabla^2{\mathbf u}\|_{L^\frac{3p^-}{p^-+1}(\Omega)}^\frac{6p^-}{7(5p^--4)}
		\|{\mathbf u}\|_{L^{12}(\Omega)}^\frac{29p^--28}{7(5p^--4)}\nonumber \\
	  &\leqslant\|\rho\|_{L^\infty(\Omega)}\|\nabla{\mathbf u}\|_{L^\frac{12}{5}(\Omega)}^\frac{64p^--56}{7(5p^--4)}
		({\mathcal I}_\Phi({\mathbf u})+1)^\frac{41p^--28}{7p^-(5p^--4)} \nonumber\\
	 &\leqslant\varepsilon{\mathcal I}_\Phi({\mathbf u})
		+C(1+\|\nabla{\mathbf u}\|_{L^\frac{12}{5}(\Omega)}^{R_4}+\|\rho\|_{L^\infty(\Omega)}^{R_5})\label{J13}
\end{align}
hold for any fixed $\varepsilon\in(0,1),$ where
$R_2=\frac{2(5p^--4)}{5p^--7},$ $R_3=\frac{24(p^--1)}{5p^--7},$ $R_4=\frac{2p^-(64p^--56)}{35(p^-)^2-69p^-+28},$ $R_5=\frac{14p^-(5p^--4)}{35(p^-)^2-69p^-+28}.$
The second term and the fourth term on the right-hand side of equation \eqref{delu} are estimated as
\begin{align}
	&|J_2|+|J_4|\nonumber\\
	&\leqslant C\|\rho\|_{L^\infty(\Omega)}\|f\|_{L^6(\Omega)}\|{\mathbf u}\|_{W^{2,\frac{3p^-}{p^-+1}}(\Omega)}
	  +C\|{\mathbf u}\|_{W^{2,\frac{3p^-}{p^-+1}}(\Omega)}  \|\rho\|_{L^\infty(\Omega)}^{\gamma - 1}\|\nabla\rho\|_{L^{3p^-}(\Omega)} \nonumber\\
	&\leqslant C\|\rho\|_{L^\infty(\Omega)}({\mathcal I}_\Phi({\mathbf u})+1)^\frac{1}{p^-}
		+C\|\rho\|_{L^\infty(\Omega)}^{\gamma - 1}\|\nabla\rho\|_{L^{3p^-}(\Omega)}({\mathcal I}_\Phi({\mathbf u})+1)^\frac{1}{p^-} \nonumber\\
	&\leqslant\varepsilon{\mathcal I}_\Phi({\mathbf u})+C(1+\|\rho\|_{L^\infty(\Omega)}^{R_6}+\|\nabla\rho\|_{L^{3p^-}(\Omega)}^{2R_6}
		+\|\rho\|_{L^\infty(\Omega)}^{2(\gamma - 1)R_6})\label{J2J4}
\end{align}
holds for any fixed $\varepsilon\in(0,1),$ where $R_6=\frac{p^-}{p^--1}.$ The third term on the right-hand side of equation \eqref{delu} is estimated as
\begin{align}
|J_3|&\leqslant\|\rho\|_{L^\infty(\Omega)}^\frac{1}{2}\|\sqrt{\rho}\partial_t{\mathbf u}\|_{L^\frac{3p^-}{2p^--1}(\Omega)}
		\|\Delta{\mathbf u}\|_{L^\frac{3p^-}{p^-+1}(\Omega)} \nonumber\\
	&\leqslant\|\rho\|_{L^\infty(\Omega)}^\frac{1}{2}\|\sqrt{\rho}\partial_t{\mathbf u}\|_{L^2(\Omega)}^\frac{4p^--4}{3p^--2}
		\|\sqrt{\rho}\partial_t{\mathbf u}\|_{L^{3p^-}(\Omega)}^\frac{2 - p^-}{3p^--2}\|{\mathbf u}\|_{W^{2,\frac{3p^-}{p^-+1}}(\Omega)} \nonumber\\
	&\leqslant\|\rho\|_{L^\infty(\Omega)}^\frac{p^-}{3p^--2}\|\sqrt{\rho}\partial_t{\mathbf u}\|_{L^2(\Omega)}^\frac{4p^--4}{3p^--2}
	\|\partial_t{\mathbf u}\|_{W^{1,\frac{3p^-}{p^-+1}}(\Omega)}^\frac{2 - p^-}{3p^--2}({\mathcal I}_\Phi({\mathbf u})+1)^\frac{1}{p^-} \nonumber \\
	&\leqslant\|\rho\|_{L^\infty(\Omega)}^\frac{p^-}{3p^--2}\|\sqrt{\rho}\partial_t{\mathbf u}\|_{L^2(\Omega)}^\frac{4p^--4}{3p^--2}
		({\mathcal J}^\frac{1}{2}_\Phi({\mathbf u})({\mathcal I}_\Phi({\mathbf u})+1)^\frac{2 - p^-}{2p^-})^\frac{2 - p^-}{3p^--2}
		({\mathcal I}_\Phi({\mathbf u})+1)^\frac{1}{p^-} \nonumber\\
	&\leqslant\|\rho\|_{L^\infty(\Omega)}^\frac{p^-}{3p^--2}\|\sqrt{\rho}\partial_t{\mathbf u}\|_{L^2(\Omega)}^\frac{4p^--4}{3p^--2}
		({\mathcal I}_\Phi({\mathbf u})+1)^\frac{p^- + 2}{2(3p^--2)}{\mathcal J}^\frac{2 - p^-}{2(3p^--2)}_\Phi({\mathbf u}) \nonumber\\
	&\leqslant\varepsilon{\mathcal I}_\Phi({\mathbf u})+\varepsilon+\|\rho\|_{L^\infty(\Omega)}^{R_7}
			\|\sqrt{\rho}\partial_t{\mathbf u}\|_{L^2(\Omega)}^{R_8}{\mathcal J}_\Phi({\mathbf u})^{R_9}\label{J3}
\end{align}
holds for any fixed $\varepsilon\in(0,1),$ where
$R_7=\frac{2p^-}{5p^--6},$ $R_8=\frac{8(p^--1)}{5p^--6} $ and $R_9=\frac{2-p^-}{5p^--6}.$
In addition, the nonlinear main part is estimated as follows
\begin{align}\label{L1}
	&\int_\Omega\divv{\mathcal A}({\mathbb D}{\uu})\cdot\Delta{\mathbf u}~dx\nonumber\\
	&=\sum_{ijklm}\int_{\Omega}(\partial_{ij}\partial_{kl}\Phi)({\mathbb D}{\mathbf u})\partial_{m}({\mathbb D}{\mathbf u})_{ij}
       \partial_{m}({\mathbb D}{\mathbf u})_{kl}~dx
		+\sum_{klm}\int_{\Omega}(\partial_{m}\partial_{kl}\Phi)({\mathbb D}{\mathbf u})\partial_{m}({\mathbb D}{\mathbf u})_{kl}~dx\nonumber\\
		&\geqslant  {\mathcal I}_{\Phi}({\mathbf u})-C\int_{\Omega}(\widetilde{D}{\mathbf u})^{\frac{p(t,x)-1}{2}}
		\ln(1+|{\mathbb D}{\mathbf u}|)|\nabla^{2}{\mathbf u}|~dx\nonumber\\
		&\geqslant  {\mathcal I}_{\Phi}({\mathbf u})-\varepsilon\int_{\Omega}({\mathbb D}{\mathbf u})^{p(t,x)-2}|\nabla^{2}{\mathbf u}|^{2}~dx
		-C \int_{\Omega}(\widetilde{D}{\mathbf u})^{p(t,x)}\ln^2(1+|{\mathbb D}{\mathbf u}|)~dx\nonumber\\
		&\geqslant  \frac{C}{2}{\mathcal I}_{\Phi}({\mathbf u})-C_{\varepsilon}\int_{\Omega}(\widetilde{D}
		{\mathbf u})^{p(t,x)}\left(1+\ln^{2}(\widetilde{D}{\mathbf u})\right)~dx
\end{align}
due to \eqref{D}, \eqref{ID} and \eqref{I}. Note that
$$(\widetilde{D}{\mathbf u})^{p(t,x)}\left(1+\ln^2(\widetilde{D}{\mathbf u})\right)\leqslant C_\delta(\widetilde{D}{\mathbf u})^{p(t,x)+\sigma}$$
holds for every $\sigma\in (0,1).$  With the help of \eqref{delu}-\eqref{L1}, we choose suitable $\sigma$ such that $p(t,x)+\sigma\leqslant q$ and
obtain that
\begin{align}\nonumber
	{\mathcal I}_\Phi({\mathbf u})\leqslant
	& C(1+\|\rho\|_{L^\infty(\Omega)}^{R_{10}}+\|\nabla\rho\|_{L^{3p^-}(\Omega)}^{R_{11}}+\|\nabla{\mathbf u}\|_{L^\frac{12}{5}(\Omega)}^{R_{12}}\\
	&~~~+\|\widetilde{D}{\mathbf u}\|_{L^q(\Omega)}^{R_{13}}
	+\|\rho\|_{L^\infty(\Omega)}^{R_7}\|\sqrt{\rho}\partial_t{\mathbf u}\|_{L^2(\Omega)}^{R_8}{\mathcal J}_\Phi({\mathbf u})^{R_9}),
\end{align}
where $R_{10}=\max \{R_2,R_5,R_6,2(\gamma-1)R_6\},$ $R_{11}=\max \{2R_6,6\},$ $R_{12}=\max \{6,R_3,R_4\}$ and $R_{13}=\max \{R_1,q\}.$
Thus, for any fixed $r\in[1,\frac{5p^--6}{2-p^-}),$
\begin{align}\label{Ir}
	&{\mathcal I}^r_\Phi({\mathbf u})\nonumber\\
    &\leqslant C\left(1+\|\rho\|_{L^\infty(\Omega)}^{rR_{10}}+\|\nabla\rho\|_{L^{3p^-}(\Omega)}^{rR_{11}}
	 +\|\nabla{\mathbf u}\|_{L^\frac{12}{5}(\Omega)}^{rR_{12}}+\|\widetilde{D}{\mathbf u}\|_{L^q(\Omega)}^{rR_{13}}\right.\nonumber\\
	&~~~~~~\left. +\|\rho\|_{L^\infty(\Omega)}^{rR_7}\|\sqrt{\rho}\partial_t{\mathbf u}\|_{L^2(\Omega)}^{rR_8}{\mathcal J}_\Phi({\mathbf u})^{rR_9}\right)\nonumber\\
	&\leqslant\eta{\mathcal J}_\Phi({\mathbf u})+C(1+\|\nabla\rho\|_{L^{3p^-}(\Omega)}^{rR_{11}}
		+\|\nabla{\mathbf u}\|_{L^\frac{12}{5}(\Omega)}^{rR_{12}}
		+\|\widetilde{D}{\mathbf u}\|_{L^q(\Omega)}^{rR_{13}}+\|\rho\|_{L^\infty(\Omega)}^{R_{14}}
		+\|\sqrt{\rho}\partial_t{\mathbf u}\|_{L^2(\Omega)}^{R_{15}}),
\end{align}
where $R_{14}=\max \{R_{10},\frac{2rR_7(5p^--6)}{5p^--6-(2-p^-)r}\} $ and $ R_{15}=\frac{2rR_8(5p^--6)}{5p^--6-(2-p^-)r}.$
Note that
$$\partial_{t}\left((\widetilde{D}{\mathbf u})^{q}\right)
=q(\widetilde{D}{\mathbf u})^{q-2}\sum_{jk=1}^3({\mathbb D}{\mathbf u})_{jk}\partial_{t}(\mathbb{ D}{\mathbf u})_{jk}$$
holds for $q>1,$ and the assumption $\frac{12}{5}\leqslant q\leqslant \min \{6(p^--1),2p^-\}$ implies
$$\frac{12}{5}\leqslant q\leqslant2q-p^-\leqslant3p^-\mbox{ and }1<\frac{3(q-p^-)}{3p^--q}<r$$ for $r\in[1,\frac{5p^--6}{2-p^-}).$
So, one uses Lemma \ref{GN} and Lemma \ref{ine} to derive that
\begin{align}\label{Dq}
  \frac{d}{dt}(\|\widetilde{D}{\mathbf u}\|_{L^q(\Omega)}^{q})
	&\leqslant\int_{\Omega}(\widetilde{D}{\mathbf u})^{q-1}|\partial_{t}{\mathbb D}{\mathbf u}|~dx\nonumber\\
	&=q\int_{\Omega}(\widetilde{D}{\mathbf u})^{\frac{p^--2}{2}}|\partial_{t}{\mathbb D}{\mathbf u}
			|(\widetilde{D}{\mathbf u})^{q-\frac{p^-}{2}}~dx\nonumber\\
	&\leqslant C{\mathcal J}_{\Phi}^{\frac{1}{2}}({\mathbf u})
			\left(\|\widetilde{D}{\mathbf u}\|_{L^{2q-p^-}}^{2q-p^-}\right)^{\frac{1}{2}}\nonumber\\
			&\leqslant\eta{\mathcal J}_\Phi({\mathbf u})+C\|\widetilde{D}{\mathbf u}\|_{L^{2q-p^-}(\Omega)}^{2q-p^-}\nonumber\\
			&\leqslant\eta{\mathcal J}_\Phi({\mathbf u})
			+C\|\widetilde{D}{\mathbf u}\|_{L^{q}(\Omega)}^{\frac{2q(2p^--q)}{3p^--q}}
			\|\widetilde{D}{\mathbf u}\|_{L^{3p^-}(\Omega)}^\frac{3p^-(q-p^-)}{3p^--q}\nonumber\\
			&\leqslant\eta{\mathcal J}_\Phi({\mathbf u})+C\|\widetilde{D}{\mathbf u}\|_{L^{q}(\Omega)}^{\frac{2q(2p^--q)}{3p^--q}}
			({\mathcal I}_\Phi({\mathbf u})+1)^\frac{3(q-p^-)}{3p^--q}\nonumber\\
			&\leqslant\eta{\mathcal J}_\Phi({\mathbf u})+\varepsilon{\mathcal I}_\Phi^r({\mathbf u})
			+ C(\|\widetilde{D}{\mathbf u}\|_{L^{q}(\Omega)}^{R_{16}}+1)
\end{align}
holds for any fixed $\varepsilon\in(0,1),$ where $R_{16}={\frac{2qp^-r(2p^--q)}{3p^-r(3p-q)-9p^-(q-p^-)}}.$ Taking account of \eqref{Ir} and \eqref{Dq}, one obtains that
\begin{align*}
	&\frac{d}{dt}\|\widetilde{D}{\mathbf u}\|_{L^q(\Omega)}^{q}+{\mathcal I}^r_\Phi({\mathbf u})\nonumber\\
	&\leqslant\eta{\mathcal J}_\Phi({\mathbf u})+C(1+\|\nabla\rho\|_{L^{3p^-}(\Omega)}^{m_1}
		+\|\widetilde{D}{\mathbf u}\|_{L^q(\Omega)}^{m_2}+\|\rho\|_{L^\infty(\Omega)}^{m_3}+\|\sqrt{\rho}\partial_t{\mathbf u}\|_{L^2(\Omega)}^{m_4})
\end{align*}
holds for any fixed $\eta\in (0,1),$ where $m_1=rR_{11},$ $m_2=\max \{rR_{12},rR_{13},R_{16}\},$ $m_3=R_{14}\mbox{ and }m_4=R_{15}.$
\end{proof}
	
\begin{lemma}\label{LEM2}
Let $\max \{1,\frac{5}{5p^--4}\}\leqslant r<\frac{5p^--6}{2-p^-}$ and $\frac{12}{5}\leqslant q\leqslant \min \{6(p^--1),2p^-\}.$ Suppose that $(\rho,{\mathbf u})$ is a smooth solution to the problem \eqref{XT-1}-\eqref{XT-2} on $[0,T)\times\Omega$ for some $T>0.$ Then there exists positive constant $m_5=m_5(p^-),~m_6=m_6(p^-),~m_7=m_7(p^-)$ and $m_8=m_8(p^-)$ such that
\begin{align}\label{ans2}
	&\frac{d}{dt}\|\sqrt{\rho}\partial_t{\mathbf u}\|_{L^2(\Omega)}^2+{\mathcal J}_\Phi({\mathbf u})\nonumber\\
	&\leqslant\varepsilon{\mathcal I}^r_\Phi({\mathbf u})+C(1+\|\nabla\rho\|_{L^{3p^-}(\Omega)}^{m_5}
		+\|\widetilde{D}{\mathbf u}\|_{L^q(\Omega)}^{m_6}+\|\rho\|_{L^\infty(\Omega)}^{m_7}+\|\sqrt{\rho}\partial_t{\mathbf u}\|_{L^2(\Omega)}^{m_8})
\end{align}
hods for any fixed $\varepsilon\in(0,1).$
\end{lemma}
	
\begin{proof}
Differentiating \eqref{XT-1}$_2$ with respect to time, one gets that
\begin{align*}
	&\rho\partial_t(\partial_t{\mathbf u)}+\rho{\mathbf u}\cdot\nabla\partial_t{\mathbf u}
		-\divv\partial_t\left((1+|{\mathbb D}{\mathbf u}|^2)^{\frac{p(t,x)-2}2}{\mathbb D}{\mathbf u}\right)+\nabla \partial_t\rho^\gamma\\
	&=\partial_t(\rho {\mathbf f})-\partial_t\rho(\partial_t{\mathbf u}+{\mathbf u}\cdot\nabla{\mathbf u})
		-\rho\partial_t{\mathbf u}\cdot\nabla{\mathbf u}.
\end{align*}
Multiplying the above equation by $\partial_t{\mathbf u}$ and integrating the resulting equation over $\Omega$, one obtains that
\begin{align}{\label{pt}}
&\frac{1}{2}\frac{d}{dt}\int_\Omega\rho|\partial_t{\mathbf u}|^2~dx
	-\int_\Omega\divv\partial_t\left((1+|{\mathbb D}{\mathbf u}|^2)^{\frac{p(t,x)-2}2}{\mathbb D}{\mathbf u}\right)
         \cdot\partial_t{\mathbf u}~dx\nonumber\\
	&=\int_\Omega\partial_t(\rho {\mathbf f})\cdot{\partial_t\mathbf u}dx-\int_\Omega\partial_t\rho\partial_t{\mathbf u}\cdot\partial_t{\mathbf u}dx
			-\int_\Omega\partial_t\rho{\mathbf u}\cdot\nabla{\mathbf u}\cdot\partial_t{\mathbf u}dx\nonumber\\
	&~~~-\int_\Omega\rho(\partial_t{\mathbf u}\cdot\nabla{\mathbf u})\cdot{\partial_t\mathbf u}dx
     -\int_\Omega\nabla\partial_t(\rho^\gamma)\cdot\partial_t{\mathbf u}~dx\nonumber\\
	&=\sum_{i=1}^{5}I_i.
\end{align}
For $I_1$, one finds that
\begin{align}\label{I1}
I_1&=\int_\Omega\partial_t\rho {\mathbf f}\cdot\partial_t{\mathbf u}+\rho\partial_t{\mathbf f}\cdot\partial_t{\mathbf u}~dx\nonumber\\
   &= -\int_\Omega\divv(\rho{\mathbf u}){\mathbf f}\cdot\partial_t{\mathbf u}+\rho \partial_t{\mathbf f}\cdot\partial_t{\mathbf u}~ dx\nonumber\\
   &= \int_\Omega\rho{\mathbf u}\cdot\nabla {\mathbf f}\cdot\partial_t{\mathbf u}+(\rho{\mathbf u})
		\cdot ({\mathbf f}\cdot\nabla\partial_t{\mathbf u})+\rho \partial_t{\mathbf f}\cdot\partial_t{\mathbf u}~dx\nonumber\\
		&\leqslant \|\rho\|_{L^\infty(\Omega)}^\frac{1}{2}\|\nabla {\mathbf f}\|_{L^2(\Omega)}\|\sqrt{\rho}\partial_t{\mathbf u}\|_{L^2(\Omega)}
		\|{\mathbf u}\|_{L^\infty(\Omega)}\nonumber\\
	&~~~+\|\rho\|_{L^\infty(\Omega)}\|{\mathbf f}\|_{L^6(\Omega)}\|\nabla\partial_t{\mathbf u}\|_{L^\frac{3p^-}{p^+1}(\Omega)}
		\|{\mathbf u}\|_{L^{12}(\Omega)}+\|\rho\|_{L^{\infty}(\Omega)}^\frac{1}{2}\|\sqrt{\rho}\partial_t{\mathbf u}\|_{L^2(\Omega)}
        \|\partial_t{\mathbf f}\|_{L^2(\Omega)}\nonumber\\
		&\leqslant  C\|\rho\|_{L^\infty(\Omega)}^\frac{1}{2}\|\sqrt{\rho}\partial_t{\mathbf u}\|_{L^2(\Omega)}
		\|{\mathbf u}\|^{1-\theta}_{L^{12}(\Omega)}\|\nabla^2{\mathbf u}\|_{L^{\frac{3p^-}{3p^-+1}}(\Omega)}^\theta\nonumber\\
	&~~~+C\|\rho\|_{L^\infty(\Omega)}\|\partial_t{\mathbf u}\|_{W^{1,\frac{3p^-}{p^-+1}}(\Omega)}\|\nabla{\mathbf u}\|_{L^\frac{12}{5}(\Omega)}
			+C\|\rho\|_{L^\infty(\Omega)}+C\|\sqrt{\rho}\partial_t{\mathbf u}\|_{L^2(\Omega)}^2\nonumber\\
			&\leqslant  C\|\rho\|_{L^\infty(\Omega)}^\frac{1}{2}\|\sqrt{\rho}\partial_t{\mathbf u}\|_{L^2(\Omega)}
			\|\widetilde{D}{\mathbf u}\|_{L^q(\Omega)}^{1-\theta}({\mathcal I}_\Phi({\mathbf u})+1)^\frac{\theta}{p^-}\nonumber\\
	&~~~+C\|\rho\|_{L^\infty(\Omega)}({\mathcal J}_\Phi({\mathbf u})+{\mathcal I}_\Phi({\mathbf u})+1)^\frac{1}{p^-}
			\|\widetilde{D}{\mathbf u}\|_{L^q(\Omega)}+C\|\rho\|_{L^\infty(\Omega)}+C\|\sqrt{\rho}\partial_t{\mathbf u}\|_{L^2(\Omega)}^2\nonumber\\
			&\leqslant \varepsilon{\mathcal I}_\Phi({\mathbf u})+\eta{\mathcal J}_\Phi({\mathbf u})
			+C(1+\|\sqrt{\rho}\partial_t{\mathbf u}\|_{L^2(\Omega)}^{R_{17}}+\|\rho\|_{L^\infty(\Omega)}^{R_{18}}
			+\|\widetilde{D}{\mathbf u}\|_{L^q(\Omega)}^{R_{19}})
\end{align}
holds for any fixed $\varepsilon\in(0,1),$ where $\theta=\frac{p^-}{5p^--4}, R_{17}=\max \{\frac{3(5p^--4)}{5p^--5},2\},
R_{18}=\max \{\frac{3(5p^--4)}{10(p^--1)},2R_6\}$ and $R_{19}=\max \{\frac{12}{5},2R_6\}.$
For $I_2,$ one takes a similar argument to obtain that
\begin{eqnarray*}\label{I2}
&I_2 &= \int_\Omega\divv(\rho{\mathbf u})|\partial_t{\mathbf u}|^2~dx\nonumber\\
	 &&= -\int_\Omega\rho{\mathbf u}\partial_t{\mathbf u}\nabla\partial_t{\mathbf u}~ dx\nonumber\\
	 &&\leqslant \|\rho\|_{L^{\infty}(\Omega)}^\frac{1}{2}\|{\mathbf u}\|_{L^\infty(\Omega)}
		\|\sqrt{\rho}\partial_t{\mathbf u}\|_{L^\frac{3p^-}{2p^--1}(\Omega)}\|\nabla\partial_t{\mathbf u}\|_{L^\frac{3p^-}{p^-+1}(\Omega)}\nonumber\\
	&&\leqslant \|\rho\|_{L^\infty(\Omega)}^{\frac{1}{2}}\|{\mathbf u}\|_{L^{12}(\Omega)}^{1-\theta}
			\|\nabla^2{\mathbf u}\|_{L^{\frac{3p^-}{p^-+1}}(\Omega)}^\theta\|\sqrt{\rho}\partial_t{\mathbf u}\|_{L^2(\Omega)}^{1-\theta_1}
			\|\sqrt{\rho}\partial_t{\mathbf u}\|_{L^{3p^-}(\Omega)}^{\theta_1}\cdot\nonumber\\
    &&~~~~{\mathcal J}^\frac{1}{2}_\Phi({\mathbf u})
			({\mathcal I}_\Phi({\mathbf u})+1)^\frac{2-p^-}{2p^-}\nonumber\\
	&&\leqslant \|\rho\|_{L^\infty(\Omega)}^{\frac{1+\theta_1}{2}}\|\widetilde{D}{\mathbf u}\|_{L^q(\Omega)}^{1-\theta}
			\|\sqrt{\rho}\partial_t{\mathbf u}\|_{L^2(\Omega)}^{1-\theta_1}
			({\mathcal I}_\Phi({\mathbf u})+1)^\frac{\theta}{p^-}{\mathcal J}^\frac{1}{2}_\Phi({\mathbf u})
			({\mathcal I}_\Phi({\mathbf u})+1)^\frac{2-p^-}{2p^-})^{\theta_1+1}\nonumber\\
	&&\leqslant \|\rho\|_{L^\infty(\Omega)}^{\frac{1+\theta_1}{2}}\|\widetilde{D}{\mathbf u}\|_{L^q(\Omega)}^{1-\theta}
			\|\sqrt{\rho}\partial_t{\mathbf u}\|_{L^2(\Omega)}^{1-\theta_1}
			({\mathcal I}_\Phi({\mathbf u})+1)^\frac{2\theta+(2-p^-)(\theta_1+1)}{2p^-}{\mathcal J}^\frac{\theta_1+1}{2}_\Phi({\mathbf u}),\nonumber
\end{eqnarray*}
where $\theta_1=\frac{2-p^-}{3p^--2}.$ The fact
$$\max \{1,\frac{5}{5p^--4}\}\leqslant r<\frac{5p^--6}{2-p^-},$$
implies that
$$\frac{\theta_1+1}{2}+\frac{2\theta+(2-p^-)(\theta_1+1)}{2p^-r}<1.$$
So,
\begin{align}\label{I2}
I_2 &\leqslant \varepsilon{\mathcal I}^r_\Phi({\mathbf u})+\eta{\mathcal J}_\Phi({\mathbf u})
		+C(1+\|\rho\|_{L^\infty(\Omega)}^{R_{20}}+\|\widetilde{D}{\mathbf u}\|_{L^q(\Omega)}^{R_{21}}
		+\|\sqrt{\rho}\partial_t{\mathbf u}\|_{L^2(\Omega)}^{R_{22}})
\end{align}
holds for any fixed $\varepsilon,~\eta\in(0,1),$ where %$\frac{1}{\delta_1}=\frac{p^-r(5p^--4)-5(p^-)^2+17p^--10}{r(5p^--4)(3p^--2)}$,
$R_{20}=\frac{3\delta_1p^-}{(3p^--2)(\delta_1-1)},$ $R_{21}=\frac{12\delta_1(p^--1)}{(5p^--4)(\delta_1-1)}$
and $R_{22}=\frac{12\delta_1(p^--1)}{(3p^--2)(\delta_1-1)}.$
The third term on the righthand side of \eqref{pt} is estimated as follows
\begin{align}\label{I3}
I_3&=\int_\Omega\divv(\rho{\mathbf u})({\mathbf u}\cdot\nabla{\mathbf u}\cdot\partial_t{\mathbf u})~dx\nonumber\\
	&=-\int_\Omega\rho({\mathbf u}\cdot\nabla{\mathbf u})\cdot(\nabla{\mathbf u}\cdot\partial_t{\mathbf u})~dx
			-\int_\Omega\rho{\mathbf u}\otimes{\mathbf u}\cdot\nabla^2{\mathbf u}\cdot\partial_t{\mathbf u}~dx\nonumber\\
	&~~~~-\int_\Omega\rho({\mathbf u}\cdot\nabla{\mathbf u})\cdot({\mathbf u}\cdot\nabla\partial_t{\mathbf u})~dx\nonumber\\
			&\triangleq I_{31}+I_{32}+I_{33}.
\end{align}
For $I_{31},$ one finds that
\begin{align}\label{I_{31}}
	|I_{31}|&\leqslant\|\rho\|_{L^\infty(\Omega)}^\frac{1}{2}\|{\mathbf u}\|_{L^\infty(\Omega)}
			\|\sqrt{\rho}\partial_t{\mathbf u}\|_{L^2(\Omega)}\|\nabla{\mathbf u}\|_{L^4(\Omega)}^2\nonumber\\
			&\leqslant\|\rho\|_{L^\infty(\Omega)}^\frac{1}{2}\|{\mathbf u}\|_{L^{12}(\Omega)}^{1-\theta}
			\|{\mathbf u}\|_{W^{2,\frac{3p^-}{p^-+1}(\Omega)}}^\theta\|\sqrt{\rho}\partial_t{\mathbf u}\|_{L^2(\Omega)}\|\nabla
			{\mathbf u}\|_{L^\frac{12}{5}(\Omega)}^{2(1-\theta_2)}\|\nabla{\mathbf u}\|_{L^{3p^-}(\Omega)}^{2\theta_2}\nonumber\\
			&\leqslant\|\rho\|_{L^\infty(\Omega)}^\frac{1}{2}\|\sqrt{\rho}\partial_t{\mathbf u}\|_{L^2(\Omega)}
			\|\widetilde{D}{\mathbf u}\|_{L^q(\Omega)}^{3-2\theta_2-\theta}({\mathcal I}_\Phi({\mathbf u})+1)^\frac{\theta+2\theta_2}{p^-}\nonumber\\
			&\leqslant\varepsilon{\mathcal I}^r_\Phi({\mathbf u})+C(1+\|\rho\|_{L^\infty(\Omega)}^{R_{23}}
			+\|\sqrt{\rho}\partial_t{\mathbf u}\|_{L^2(\Omega)}^{R_{24}}+\|\widetilde{D}{\mathbf u}\|_{L^q(\Omega)}^{R_{25}})
\end{align}
holds for any fixed $\varepsilon\in(0,1),$ where
$\theta=\frac{p^-}{5p^--4},~\theta_2=\frac{2p^-}{5p^--4},$  %$\frac{\theta+2\theta_2}{p^-r}<1,$
$R_{23}=\frac{3r(5p^--4)}{2(r(5p^--4)-5)},$ $R_{24}=\frac{3r(5p^--4)}{r(5p^--4)-5}$ and $R_{25}=\frac{6r(5p^--6)}{r(5p^--4)-5}.$

For $I_{32},$ one gets that
\begin{align}\label{I_{32}}
|I_{32}|&\leqslant\|\rho\|_{L^\infty(\Omega)}^\frac{1}{2}\|\sqrt{\rho}\partial_t{\mathbf u}\|_{L^\frac{3p^-}{2p^--1}(\Omega)}
			\|\nabla^2{\mathbf u}\|_{L^\frac{3p^-}{p^-+1}}\|{\mathbf u}\|_{L^\infty(\Omega)}^2\nonumber\\
		&\leqslant\|\rho\|_{L^\infty(\Omega)}^\frac{1}{2}\|\sqrt{\rho}\partial_t{\mathbf u}\|_{L^2(\Omega)}^{1-\theta_1}
			\|\sqrt{\rho}\partial_t{\mathbf u}\|_{L^{3p^-}(\Omega)}^{\theta_1}\|{\mathbf u}\|_{W^{2,\frac{3p^-}{p^-+1}}}
			\|{\mathbf u}\|_{L^{12}(\Omega)}^{2(1-\theta)}\|\nabla^2{\mathbf u}\|_{L^{\frac{3p^-}{p^-+1}}(\Omega)}^{2\theta}\nonumber\\
			&\leqslant\|\rho\|_{L^\infty(\Omega)}^\frac{1+\theta_1}{2}\|\sqrt{\rho}\partial_t{\mathbf u}\|_{L^2(\Omega)}^{1-\theta_1}
			\|\widetilde{D}{\mathbf u}\|_{L^q(\Omega)}^{2(1-\theta)}\|{\mathbf u}\|_{W^{2,\frac{3p^-}{p^-+1}}(\Omega)}^{2\theta+1}
			\|\partial_t{\mathbf u}\|_{W^{1,\frac{3p^-}{p^-+1}}}^{\theta_1}\nonumber\\
			&\leqslant\|\rho\|_{L^\infty(\Omega)}^\frac{1+\theta_1}{2}\|\sqrt{\rho}\partial_t{\mathbf u}\|_{L^2(\Omega)}^{1-\theta_1}
			\|\widetilde{D}{\mathbf u}\|_{L^q(\Omega)}^{2(1-\theta)}({\mathcal I}_\Phi({\mathbf u})+1)^\frac{(2-p^-)\theta_1+4\theta+2}{2p^-}
			{\mathcal J}_\Phi({\mathbf u})^\frac{\theta}{2}\nonumber\\
			&\leqslant\varepsilon{\mathcal I}^r_\Phi({\mathbf u})+\eta{\mathcal J}_\Phi({\mathbf u})
			+C(1+\|\rho\|_{L^\infty(\Omega)}^{R_{26}}+\|\sqrt{\rho}\partial_t{\mathbf u}\|_{L^2(\Omega)}^{R_{27}}
			+\|\widetilde{D}{\mathbf u}\|_{L^q(\Omega)}^{R_{28}})
\end{align}
holds for any fixed $\varepsilon,\eta\in(0,1),$	where
$\frac{(2-p^-)\theta_1+4\theta+2+\theta p^-}{2p^-}<1,$ $\frac{1}{\delta_2}=\frac{5(p^-)^2-16+18p^-+(2-p^-)(5p^--4)r}{2r(3p^--2)(5p^--4)},$ $R_{26}=\frac{3p^-\delta_2}{(3p^--2)(\delta_2-1)},$ $R_{27}=\frac{12\delta_2(p^--1)}{(3p^--2)(\delta_2-1)}
$ and $R_{28}=\frac{12\delta_2(p^--1)}{(5p^--4)(\delta_2-1)}.$
For $I_{33},$ one has that
\begin{align}\label{I_{33}}
|I_{33}|	&\leqslant\|\rho\|_{L^\infty(\Omega)}\|{\mathbf u}\|_{L^\infty(\Omega)}\|{\mathbf u}\|_{L^{84}(\Omega)}
			\|\nabla{\mathbf u}\|_{L^{\frac{12}{5}}(\Omega)}\|\nabla\partial_t{\mathbf u}\|_{L^\frac{3p^-}{p^-+1}(\Omega)}\nonumber\\
			&\leqslant\|\rho\|_{L^\infty(\Omega)}\|{\mathbf u}\|_{L^{12}(\Omega)}^{2-\theta-\theta_3}
			\|\nabla{\mathbf u}\|_{L^\frac{12}{5}(\Omega)}\|\nabla^2{\mathbf u}\|_{L^\frac{3p^-}{p^-+1}(\Omega)}^{\theta+\theta_3}
			{\mathcal J}^\frac{1}{2}_\Phi({\mathbf u})({\mathcal I}_\Phi({\mathbf u})+1)^\frac{2-p^-}{2p^-})\nonumber\\
			&\leqslant\|\rho\|_{L^\infty(\Omega)}\|\widetilde{D}{\mathbf u}\|_{L^q(\Omega)}^{3-\theta-\theta_3}
			{\mathcal J}^\frac{1}{2}_\Phi({\mathbf u})({\mathcal I}_\Phi({\mathbf u})+1)^\frac{2-p^-+2\theta+2\theta_3}{2p^-}\nonumber\\
			&\leqslant\varepsilon{\mathcal I}^r_\Phi({\mathbf u})+\eta{\mathcal J}_\Phi({\mathbf u})
			+C(1+\|\rho\|_{L^\infty(\Omega)}^{R_{29}}+\|\widetilde{D}{\mathbf u}\|_{L^q(\Omega)}^{R_{30}})
\end{align}
holds for any fixed $\varepsilon,\eta\in(0,1),$ where
$\theta_3=\frac{6p^-}{7(5p^--4)},~\frac{1}{\delta_3}=\frac{-35(p^-)^2(1+r)-56+124p^--28p^-r}{14p^-r(5p^--4)},$ $R_{29}=2\frac{\delta_3}{\delta_3-1}$ and $ R_{30}=\frac{(92p^--84)\delta_3}{7(5p^--4)(\delta_3-1)}.$
Note that
\begin{align*}
	\|{\mathbb D}\partial_t{\mathbf u}\|_{L^\frac{2q}{2-p^-+q}}
		\leqslant C{\mathcal J}^\frac{1}{2}_\Phi({\mathbf u})\|(\widetilde{D}{\mathbf u})^\frac{2-p^-}{2}\|_{L^\frac{2q}{2-p^-}(\Omega)}
		\leqslant C{\mathcal J}^\frac{1}{2}_\Phi({\mathbf u})\|(\widetilde{D}{\mathbf u})\|_{L^q(\Omega)}^\frac{2-p^-}{2}.
\end{align*}
It is deduced from  Lemma \ref{GIJ} that
\begin{align}\label{I4}
I_{4}&\leqslant\|\rho\|_{L^\infty(\Omega)}^\frac{1}{2}\|\nabla{\mathbf u}\|_{L^q(\Omega)}
			\|\sqrt{\rho}\partial_t{\mathbf u}\|_{L^\frac{2q}{q-1}}^2\nonumber\\
	&\leqslant\|\rho\|_{L^\infty(\Omega)}^\frac{1}{2}\|\nabla{\mathbf u}\|_{L^q(\Omega)}
			\|\sqrt{\rho}\partial_t{\mathbf u}\|_{L^2(\Omega)}^{2(1-\theta_4)}
			\|\sqrt{\rho}\partial_t{\mathbf u}\|_{L^\frac{6q}{6-3p^-+q}}^{2\theta_4}\nonumber\\
			&\leqslant\|\rho\|_{L^\infty(\Omega)}^\frac{1+2\theta_4}{2}\|\nabla{\mathbf u}\|_{L^q(\Omega)}
			\|\sqrt{\rho}\partial_t{\mathbf u}\|_{L^2(\Omega)}^{2(1-\theta_4)}
			\|\nabla\partial_t{\mathbf u}\|_{L^\frac{2q}{2-p^-+q}}^{2\theta_4}\nonumber\\
			&\leqslant\|\rho\|_{L^\infty(\Omega)}^\frac{1+2\theta_4}{2}\|\nabla{\mathbf u}\|_{L^q(\Omega)}
			\|\sqrt{\rho}\partial_t{\mathbf u}\|_{L^2(\Omega)}^{2(1-\theta_4)}{\mathcal J}^\frac{1}{2}_\Phi({\mathbf u})
			\|\widetilde{D}{\mathbf u}\|_{L^q(\Omega)}^\frac{2-p^-}{2})^{2\theta_4}\nonumber\\
			&\leqslant\eta{\mathcal J}_\Phi({\mathbf u})+C(\|\rho\|_{L^\infty(\Omega)}^{R_{31}}
			+\|\sqrt{\rho}\partial_t{\mathbf u}\|_{L^2(\Omega)}^{R_{32}}+\|\widetilde{D}{\mathbf u}\|_{L^q(\Omega)}^{R_{33}})
\end{align}
holds for any fixed $\eta\in(0,1),$ where $\theta_4=\frac{3}{2q+3p^--6},~R_{31}=\frac{3(2q+3p^-)}{2(2q+3p^--9)},$ $R_{32}=6\text{ and }R_{33}=\frac{6q}{2q+3p^--9}.$
		
In addition, it follows from \eqref{XT-1}$_1$ 	that
\begin{equation}\label{transform}
	\partial_t\rho^\gamma+\operatorname{div}(\rho^\gamma{\mathbf u})+(\gamma-1)\rho^\gamma\operatorname{div}{\mathbf u}=0.
\end{equation}
One estimates $ I_{5}$ as follows
\begin{align}\label{I5}
|I_{5}|&=| \int_\Omega \partial_t \rho^\gamma\divv \partial_t{\mathbf u}\,dx| \nonumber\\
 	&=  -\intO (\nabla \rho^\gamma\cdot {\mathbf u} +\gamma \rho^\gamma\divv {\mathbf u})\divv \partial_t{\mathbf u}\,dx\nonumber\\
			&\leqslant  C\|\rho\|_{L^\infty(\Omega)}^{\gamma-1}\|\nabla\rho\|_{L^{3p^-}(\Omega)}
			\|\nabla\partial_t{\mathbf u}\|_{L^{\frac{3p^-}{p^-+1}}(\Omega)}\|{\mathbf u}\|_{L^{12}(\Omega)}\nonumber\\
         &~~+\|\rho\|_{L^\infty(\Omega)}^\gamma\|\nabla{\mathbf u}\|_{L^\frac{3p^-}{2p^--1}(\Omega)}
			\|\nabla\partial_t{\mathbf u}\|_{L^{\frac{3p^-}{p^-+1}}}\nonumber\\
			&\leqslant  C\|\rho\|_{L^\infty(\Omega)}^{\gamma-1}\|\nabla\rho\|_{L^{3p^-}(\Omega)}({\mathcal J}_\Phi({\mathbf u})
			+{\mathcal I}_\Phi({\mathbf u})+1)^\frac{1}{p^-}\|\widetilde{D}{\mathbf u}\|_{L^q(\Omega)}\nonumber\\
			&~~~+\|\rho\|_{L^{\infty}(\Omega)}^\gamma\|\nabla{\mathbf u}\|_{L^2(\Omega)}^{1-\theta_1}
			\|\nabla{\mathbf u}\|_{L^{3p^-}(\Omega)}^{\theta_1}({\mathcal J}_\Phi({\mathbf u})
			+{\mathcal I}_\Phi({\mathbf u})+1)^\frac{1}{p^-}\nonumber\\
			&\leqslant  C\|\rho\|_{L^\infty(\Omega)}^{\gamma-1}\|\nabla\rho\|_{L^{3p^-}(\Omega)}({\mathcal J}_\Phi({\mathbf u})
			+{\mathcal I}_\Phi({\mathbf u})+1)^\frac{1}{p^-}\|\widetilde{D}{\mathbf u}\|_{L^q(\Omega)}\nonumber\\
			&~~~+\|\rho\|_{L^{\infty(\Omega)}}^\gamma\|\widetilde{D}{\mathbf u}\|_{L^q(\Omega)}^{1-\theta_1}
			({\mathcal I}_\Phi({\mathbf u})+1)^\frac{\theta_1}{p^-}({\mathcal J}_\Phi({\mathbf u})
			+{\mathcal I}_\Phi({\mathbf u})+1)^\frac{1}{p^-}\nonumber\\
			&\leqslant  \varepsilon{\mathcal I}_\Phi({\mathbf u})+\eta{\mathcal J}_\Phi({\mathbf u})
       + C(1+\|\rho\|_{L^\infty(\Omega)}^{R_{34}}+\|\nabla\rho\|_{L^{3p^-}(\Omega)}^{R_{35}}+\|\widetilde{D}{\mathbf u}\|_{L^q(\Omega)}^{R_{36}})
\end{align}
holds for any fixed $\varepsilon,\eta\in(0,1),$ where $R_{34}=\max \{3(\gamma-1)R_6,\frac{2\gamma(3p^--2)}{3p^--4}\},$ $R_{36}=\max \{3R_6,\frac{8p^-}{3p^--4}\}$ and $R_{35}=3R_6.$  Moreover, by the same argument as in \eqref{L1}, one finds that
\begin{align}\label{L2}
	&-\int_\Omega\divv\partial_t\left((1+|{\mathbb D}{\mathbf u}|^2)^{\frac{p(t,x)-2}2}{\mathbb D}{\mathbf u}\right)\cdot\partial_t{\mathbf u}~dx
			\geqslant{\mathcal J}_\Phi({\mathbf u})-\|\widetilde{D}{\mathbf u}\|_{L^q(\Omega)}^q.
\end{align}
		
Therefore, it deduces from \eqref{pt}-\eqref{L2} that
\begin{align*}
	&\frac{d}{dt}\|\sqrt{\rho}\partial_t{\mathbf u}\|_{L^2(\Omega)}^2+{\mathcal J}_\Phi({\mathbf u})\nonumber\\
	&\leqslant \varepsilon{\mathcal I}^r_\Phi({\mathbf u})+C(1+\|\nabla\rho\|_{L^{3p^-}(\Omega)}^{m_5}
	+\|\widetilde{D}{\mathbf u}\|_{L^q(\Omega)}^{m_6}+\|\rho\|_{L^\infty(\Omega)}^{m_7}+\|\sqrt{\rho}\partial_t{\mathbf u}\|_{L^2(\Omega)}^{m_8}),
\end{align*}
holds for any fixed $\varepsilon\in(0,1),$ where $m_6=\max \{R_{19},R_{21},R_{25},R_{28},R_{30},R_{33},R_{36},q\},$ $m_5=R_{35},$
		$m_7=\max \{R_{18},R_{20},R_{23},R_{26},R_{29},R_{31},R_{34}\}$ and $m_8=\max \{R_{17},R_{22},R_{24},R_{27},R_{32}\}.$
	\end{proof}
	
Finally, an estimate for the higher-order derivatives of the density is given.
	
\begin{lemma}\label{LEM3}
Let $\max \{1,\frac{5}{5p^--4}\}\leqslant r<\frac{5p^--6}{2-p^-}$ and $\frac{12}{5}\leqslant q\leqslant \min \{6(p^--1),2p^-\}.$ Suppose that $(\rho,{\mathbf u})$ is a smooth solution to \eqref{XT-1}-\eqref{XT-2} on $[0,T)\times\Omega$ for some $T>0.$ Then there exists a positive constant $m_9=m_9(p^-),~m_{10}=m_{10}(p^-)$ and $m_{11}=m_{11}(p^-)$ such that
\begin{align}\label{ans3}
	\frac{d}{dt}\|\rho\|_{W^{1,{3p^-}}(\Omega)}	&
	\leqslant \varepsilon{\mathcal I}^r_\Phi({\mathbf u})+\eta{\mathcal J}_\Phi({\mathbf u})
			+C(1+\|\widetilde{D}{\mathbf u}\|_{L^q(\Omega)}^{m_9}+\|\rho\|_{W^{1,3p^-}}^{m_{10}})
		\end{align}
holds for any fixed $\varepsilon,\eta\in(0,1).$
\end{lemma}
	
\begin{proof}
The equation $\eqref{XT-1}_{1}$ implies that
\begin{equation}\label{cont-pressure}
	\partial_t\rho+\rho\divv{\mathbf u}+\nabla\rho\cdot{\mathbf u}=0.
\end{equation}
Differentiating \eqref{cont-pressure} with respect to $x_k$ yields
\begin{equation}\label{xk}
  \partial_t(\rho_{x_k})+{\mathbf u}\cdot \nabla \rho_{x_k}+\nabla \rho \cdot {\mathbf u}_{x_k}+\rho_{x_k}\divv {\mathbf u}
  +\rho\divv {\mathbf u}_{x_k} =0.
\end{equation}
Now, multiplying (\ref{xk}) by $\rho_{x_k}|\rho_{x_k}|^{3p^--2},$ integrating the resulting equation over $\Omega$ and summing  over $k$, we get that
\begin{align*}
  \frac{d}{dt}\|\nabla \rho\|^{3p^-}_{L^{3p^-}(\Omega)}
		&\leqslant  C\int_\Omega |\nabla {\mathbf u}\|\nabla \rho|^{3p^-}\,dx+\int_\Omega \rho|\nabla \divv {\mathbf u}\|\nabla \rho|^{3p^--1}\,dx\\
		&\leqslant  C \big( \|\nabla {\mathbf u}\|_{L^{\infty}(\Omega)}\|\nabla \rho\|_{L^{3p^-} (\Omega)}^{3p^-}
			+\|\rho\|_{L^\infty(\Omega)} \|\nabla \divv {\mathbf u}\|_{L^{3p^-}(\Omega)}\|\nabla \rho\|_{L^{3p^-} (\Omega)}^{3p^--1}\big),
\end{align*}
which implies that
\begin{align*}
	\frac{d}{dt}\|\nabla \rho\|_{L^{3p^-}(\Omega)}
      \leqslant C(\|\nabla{\mathbf u}\|_{L^\infty(\Omega)}\|\nabla \rho\|_{L^{3p^-}(\Omega)}
			+\|\rho\|_{L^\infty(\Omega)}\|\nabla^2{\mathbf u}\|_{L^{3p^-}(\Omega)}).
\end{align*}
Similarly, we can derive a bound on $\|\rho\|_{L^{3p^-}(\Omega)}$ by multiplying \eqref{cont-pressure} by $\rho^{3p^--2}$ and integrating the resultant over $\Omega.$ Thus,
\begin{align}\label{rho1}
	\frac{d}{dt}\|\rho\|_{W^{1,{3p^-}}(\Omega)}\leqslant C\|{\mathbf u}\|_{W^{2,{3p^-}}(\Omega)}\|\rho\|_{W^{1,{3p^-}}(\Omega)}.
\end{align}
		
In order to complete the estimate of \eqref{rho1}, one has to show the $L^s-$estimates on the elliptic operator associated to the system based on the results stated in Lemma \ref{lemzzx}. Indeed, we write the momentum equation as follows
\begin{equation}\label{lp-prob}
	-\mathcal{A}({\mathbf u}^*,\mathbb{D}){\mathbf u}=-\rho \partial_t{\mathbf u}-\rho ({\mathbf u}\cdot\nabla) {\mathbf u}
    -\nabla \rho^\gamma+\rho\mathbf{f}+(\mathcal{A}({\mathbf u},\mathbb{D}){\mathbf u}-\mathcal{A}({\mathbf u}^*,\mathbb{D}){\mathbf u}),
\end{equation}
where ${\mathbf u}^*$ is a reference solution. Regarding \eqref{lp-prob} as a quasi-linear elliptic equation, one finds that the linearized operator still yield maximal $L^s-$regularity according to Lemma \ref{lemzzx}. Since
\begin{align*}
&{\mathcal A}({\mathbf u},\mathbb{D}){\mathbf u}-\mathcal{A}({\mathbf u}^*,\mathbb{D}){\mathbf u}\\
		&= \frac{1}{2}((1+|{\mathbb D}{\mathbf u}|^2 )-(1+|{\mathbb D}{\mathbf u}^*|^2))(\Delta {\mathbf u}+\nabla \divv {\mathbf u}) +\\
		&{\small \underset{j,k,l=1}{\overset{3}{\sum}}(p(t,x)-2)\left((1+|{\mathbb D}{\mathbf u}|^2)^\frac{p(t,x)-4}{2}({\mathbb D}{\mathbf u})_{ik}({\mathbb D}{\mathbf u})_{jl}-(1+|{\mathbb D}{\mathbf u}^*|^2)^\frac{p(t,x)-4}{2}({\mathbb D}{\mathbf u}^*)_{ik}({\mathbb D}{\mathbf u}^*)_{jl}\right)\del_k \del_l{\mathbf u}_j},
\end{align*}
one gets that
\begin{align*}\label{difA-A*-est}
  &\|\mathcal{A}({\mathbf u},\mathbb{D}){\mathbf u}-\mathcal{A}({\mathbf u}^{*},\mathbb{D}){\mathbf u}\|_{L^{s}(\Omega)}\\
  &\leqslant  C\||{\mathbb D}{\mathbf u}|^{2}-|{\mathbb D}{\mathbf u}^*|^{2}\|_{L^{\infty}(\Omega)}\|{\mathbf u}\|_{W^{2,s}(\Omega)}\\
  &~~~ +C\sum_{k,l=1}^{3}
\||({\mathbb D}{\mathbf u})_{ik}({\mathbb D}{\mathbf u})_{jl}-({\mathbb D}{\mathbf u}^*)_{ik}({\mathbb D}{\mathbf u}^*)_{jl}|\|_{L^{\infty}(\Omega)}
			\|{\mathbf u}\|_{W^{2,s}(\Omega)}\\
	&~~~+C\sum_{k,l=1}^{3}\|(|{\mathbb D}{\mathbf u}|^{2}-|{\mathbb D}{\mathbf u}^*|^{2})({\mathbb D}{\mathbf u}^*)_{ik}({\mathbb D}{\mathbf u}^*)_{jl}\|_{L^{\infty}(\Omega)}\|{\mathbf u}\|_{W^{2,s}(\Omega)}\\
			&\leqslant  C\||{\mathbb D}{\mathbf u}|-|{\mathbb D}{\mathbf u}^*|\|_{L^{\infty}(\Omega)}
			\||{\mathbb D}{\mathbf u}|+|{\mathbb D}{\mathbf u}^*|\|_{L^{\infty}(\Omega)}\|{\mathbf u}\|_{W^{2,s}(\Omega)}\\
			&~~~+C\sum_{k,l=1}^{3}(\|({\mathbb D}{\mathbf u})_{ik}|_{L^{\infty}(\Omega)}
			\||({\mathbb D}{\mathbf u})_{jl}|-|({\mathbb D}{\mathbf u}^*)_{jl}|\|_{L^{\infty}(\Omega)}\\
			&~~~+C\sum_{k,l=1}^{3}\|({\mathbb D}{\mathbf u}^*)_{jl}\|_{L^{\infty}(\Omega)}
			\||({\mathbb D}{\mathbf u})_{ik}|-|({\mathbb D}{\mathbf u}^*)_{ik}|\|_{L^{\infty}(\Omega)})\|{\mathbf u}\|_{W^{2,s}(\Omega)}\\
			&~~~+C\sum_{k,l=1}^{3}\||{\mathbb D}{\mathbf u}|-|{\mathbb D}{\mathbf u}^*|\|_{L^{\infty}(\Omega)}
			\|({\mathbb D}{\mathbf u}^*)_{ik}({\mathbb D}{\mathbf u}^*)_{jl}\|_{L^{\infty}(\Omega)}
			\||{\mathbb D}{\mathbf u} |+|{\mathbb D}{\mathbf u}^*|\|_{L^{\infty}(\Omega)}\|{\mathbf u}\|_{W^{2,s}(\Omega)}\\
			&\leqslant  C\|\nabla({\mathbf u}-{\mathbf u}^{*})\|_{L^{\infty}(\Omega)}\|{\mathbf u}\|_{W^{2,s}(\Omega)}.	
\end{align*}
So, by choosing the reference solution ${\mathbf u}^*$ of \eqref{lp-prob} such that ${\mathbf u}$ is a small perturbation of ${\mathbf u}^*,$ one deduces from \eqref{lp-prob} that
\begin{eqnarray}\label{3dlg-E19}
	&&\|{\mathbf u} \|_{W^{2,s}(\Omega)}\\
	&&\leqslant C \left(\|-\rho \partial_t{\mathbf u}-\rho ({\mathbf u}\cdot\nabla) {\mathbf u}-\nabla \rho^\gamma+\rho\mathbf{f}\|_{L^{s}(\Omega)}
			+\|\nabla{\mathbf u}\|_{L^{p^-}(\Omega)}^{p^-}+1\right)+\frac12\|{\mathbf u}\|_{W^{2,s}(\Omega)}.\nonumber
\end{eqnarray}
Thus,
\begin{eqnarray}\label{u2}
  &&\|{\mathbf u}\|_{W^{2,{3p^-}}(\Omega)}\nonumber\\
  &&\leqslant  C(\|\rho \partial_t{\mathbf u}\|_{L^{3p^-}(\Omega)}+\|\rho({\mathbf u}\cdot\nabla) {\mathbf u}\|_{L^{3p^-}(\Omega)}
	  +\|\nabla \rho^\gamma\|_{L^{3p^-}(\Omega)}+\|\rho\mathbf{f}\|_{L^{3p^-}(\Omega)}\|\nabla{\mathbf u}\|_{L^{p^-}(\Omega)}^{p^-}+1)\nonumber\\
  &&\leqslant  C(\|\rho\|_{L^\infty(\Omega)}\|\partial_t{\mathbf u}\|_{L^{3p^-}(\Omega)}
			+\|\rho\|_{L^\infty(\Omega)}\|{\mathbf u}\|_{L^\infty(\Omega)}\|\nabla{\mathbf u}\|_{L^{3p^-}(\Omega)}
			+\|\rho\|_{L^\infty(\Omega)}^{\gamma-1}\|\nabla\rho\|_{L^{3p^-}(\Omega)}\nonumber\\
  &&~~~+\|\rho\|_{L^\infty(\Omega)}\|\mathbf{f}\|_{L^{3p^-}}+\|\nabla{\mathbf u}\|_{L^{p^-}(\Omega)}^{p^-}+1)  \nonumber\\
			&&\leqslant  C(\|\rho\|_{L^\infty(\Omega)}({\mathcal I}_\Phi({\mathbf u})+{\mathcal J}_\Phi({\mathbf u})+1)^\frac{1}{p^-}
			+\|\rho\|_{L^\infty(\Omega)}\|{\mathbf u}\|_{L^{12}(\Omega)}^{1-\theta}
			\|\nabla^2{\mathbf u}\|_{L^{\frac{3p^-}{p^-+1}}(\Omega)}^{1+\theta}\nonumber\\
			&&~~~+\|\rho\|_{L^\infty(\Omega)}^{\gamma-1}\|\nabla\rho\|_{L^{3p^-}(\Omega)}+\|\rho\|_{L^\infty(\Omega)}\|\mathbf{f}\|_{L^{3p^-}}
			+\|\nabla{\mathbf u}\|_{L^{p^-}(\Omega)}^{p^-}+1)\nonumber\\
			&&\leqslant  C(\|\rho\|_{L^\infty(\Omega)}({\mathcal I}_\Phi({\mathbf u})+{\mathcal J}_\Phi({\mathbf u})+1)^\frac{1}{p^-}
			+\|\rho\|_{L^\infty(\Omega)}\|\widetilde{D}{\mathbf u}\|_{L^{q}(\Omega)}^{1-\theta}
			({\mathcal I}_\Phi({\mathbf u})+1)^\frac{1+\theta}{p^-}\nonumber\\
			&&~~~+\|\rho\|_{L^\infty(\Omega)}^{\gamma-1}\|\nabla\rho\|_{L^{3p^-}(\Omega)}+\|\rho\|_{L^\infty(\Omega)}\|\mathbf{f}\|_{L^{3p^-}}
			+\|\nabla{\mathbf u}\|_{L^{p^-}(\Omega)}^{p^-}+1)\nonumber\\
			&&\leqslant  \varepsilon{\mathcal I}^r_\Phi({\mathbf u})+\eta{\mathcal J}_\Phi({\mathbf u})
			+C(1+\|\rho\|_{L^\infty(\Omega)}^{R_{37}}+\|\widetilde{D}{\mathbf u}\|_{L^q(\Omega)}^{R_{38}}+\|\nabla\rho\|_{L^{3p^-}(\Omega)}^2)
\end{eqnarray}
holds for any fixed $\varepsilon, \eta\in (0,1),$ where $R_{37}=\max \{\frac{2p^-(5p^--4)r}{p^-(5p^--4)r-6p^-+4},2(\gamma-1)\}$ and
$R_{38}=\max \{R_6,\frac{8p^-r(p^--1)}{p^-(5p^--4)r-6p^-+4}\}.$
Combining  \eqref{rho1} and \eqref{u2}, one can conclude that
\begin{align}
	&\frac{d}{dt}\|\rho\|_{W^{1,{3p^-}}(\Omega)}\nonumber\\
	&\leqslant C(\|\rho\|_{L^\infty(\Omega)}({\mathcal I}_\Phi({\mathbf u})
			+{\mathcal J}_\Phi({\mathbf u})+1)^\frac{1}{p^-}\|\rho\|_{W^{1,3p^-}(\Omega)}\nonumber\\
	&~~~+\|\rho\|_{L^\infty(\Omega)}^{\gamma-1}\|\rho\|_{W^{1,3p^-}(\Omega)}^2+\|\rho\|_{L^\infty(\Omega)}
			\|\mathbf{f}\|_{L^{3p^-}}\|\rho\|_{W^{1,3p^-}(\Omega)}\nonumber\\
	&~~~+\|\nabla{\mathbf u}\|_{L^{p^-}(\Omega)}^{p^-}\|\rho\|_{W^{1,3p^-}(\Omega)}
			+\|\rho\|_{L^\infty(\Omega)}\|\widetilde{D}{\mathbf u}\|_{L^{q}(\Omega)}^\theta
			({\mathcal I}_\Phi({\mathbf u})+1)^\frac{1+\theta}{p^-}\|\rho\|_{W^{1,3p^-}(\Omega)})\nonumber\\
			&\leqslant  \varepsilon{\mathcal I}^r_\Phi({\mathbf u})+\eta{\mathcal J}_\Phi({\mathbf u})
			+C(1+\|\widetilde{D}{\mathbf u}\|_{L^q(\Omega)}^{m_9}+\|\rho\|_{W^{1,3p^-}(\Omega)}^{m_{10}})
\end{align}
holds for any fixed $\varepsilon,\eta\in (0,1),$ where $m_9=\max \{\frac{12p^-r(p^--1)}{p^-(5p^--4)r-6p^-+4},2p^-\}$ and  $m_{10}=\max \{\frac{6p^-(5p^--4)r}{p^-(5p^--4)r-6p^-+4},2R_6,4,2(\gamma-1)\}.$
\end{proof}
	
Consequently, the following proposition follows immediately from Lemma \ref{LEM1}-Lemma \ref{LEM3}.
	
\begin{proposition}\label{prop-13}
Let $\max \{1,\frac{5}{5p^--4}\}\leqslant r<\frac{5p^--6}{2-p^-}$ and $\frac{12}{5}\leqslant q\leqslant \min \{6(p^--1),2p^-\}.$ Suppose that $(\rho,{\mathbf u})$ is a smooth solution to \eqref{XT-1}-\eqref{XT-2} on $[0,T)\times\Omega$ for some $T>0.$ Then there exists a positive constant $m_{11},~m_{12}$ and $m_{13}$ such that
\begin{align}\label{E-71}
	&\frac{d}{dt}(\|\sqrt{\rho}{\mathbf u}\|_{L^2(\Omega)}^2+\|\sqrt{\rho}\partial_t{\mathbf u}\|_{L^2(\Omega)}^2+\|\rho\|_{L^\gamma(\Omega)}
		+\|\widetilde{D}{\mathbf u}\|_{L^q(\Omega)}^{q}+\|\rho\|_{W^{1,3p^-}(\Omega)})\nonumber\\
	&~~~+\|{\mathbf u}\|_{W^{2,3p^-}(\Omega)}+{\mathcal I}^r_\Phi({\mathbf u})+{\mathcal J}_\Phi({\mathbf u})\nonumber\\
	&\leqslant C(1+\|\sqrt{\rho}{\mathbf u}\|_{L^2(\Omega)}^2+(1+\|f\|_{W^{1,2}(\Omega)}^2)\|\rho\|_{W^{1,3p^-}(\Omega)}^{m_{11}}
		+\|\widetilde{D}{\mathbf u}\|_{L^q(\Omega)}^{m_{12}}+\|\sqrt{\rho}\partial_t{\mathbf u}\|_{L^2(\Omega)}^{m_{13}}).
\end{align}
\end{proposition}

\section{Proof of Theorem 1 }\label{Sec-P}
	
\subsection{The Faedo-Galerkin approximations}
	
Based on the necessary a priori estimates derived in the previous section, the existence of strong solutions can be established by a standard argument: construct approximate solutions via a semi-discrete Galerkin scheme, derive uniform bounds, and obtain solutions by passing to the limit (see \cite{BAT2023}, \cite{hanguo}). Precisely, using a fixed-point argument combined with the characteristics method and the Faedo-Galerkin method, we construct approximate solutions in an appropriate finite-dimensional space, whose existence is then extended to the full space via the a priori estimates established in the previous section. Take $X=H^2(\Omega)$ as the basic function space as $X=H_0^1(\Omega)\cap H^2(\Omega)$ and its finite-dimensional subspaces as
	$$X^N=\operatorname{span}\{\phi^1,...,\phi^N\}\subset X\cap C^2(\overline{\Omega})\quad (N=1,2,\cdots),$$
where $\{\phi^k\}$ denotes the set consisting of the eigenvectors of the Stokes operator. Thus, $\{\phi^k\}$ is an orthogonal basis of $L^2(\Omega)$ (is also an orthogonal in $H^2(\Omega)$).
	
Suppose $\rho_0,{\mathbf u}_0$ and ${\mathbf f}$ satisfy the hypotheses of Theorem \ref{theo1}. We assume that $\rho_0\in C^1(\overline{\Omega})$ and  $\rho_0\geqslant \delta$ for some $\delta>0.$ Define ${\mathbf u}^N_\delta =\sum\limits_{r=1}^Nc_{r,\delta}^N(t)\phi^r(x).$ We look for a solution $(\rho^N_\delta ,{\mathbf u}^N_\delta )$ to the following initial problem
\begin{eqnarray}\label{XT1}
\begin{cases}
	\rho_t^N+\:\mathrm{div}\:(\rho^N_\delta {\mathbf u}^N_\delta )=0,\\
	 \int_\Omega(\rho^N_\delta \partial_t{\mathbf u}_\delta^N
			+\rho^N_\delta ({\mathbf u}^N_\delta \cdot\nabla ){\mathbf u}^N_\delta
			+\mathrm{div}((1+|{\mathbb D}{\mathbf u}^N_\delta |^2 )^{\frac{p(t,x)-2}{2}} {\mathbb D}{\mathbf u}^N_\delta )
			+\nabla (\rho^N_\delta )^\gamma)\cdot {\mathbf v}\:dx\\
			~~~~~~~~~=\int_\Omega\rho^N_\delta {\mathbf f}\cdot {\mathbf v}\:dx~~({\mathbf v}\in X^N ),\\
			\rho^N_\delta (0)=\rho_0,\quad{\mathbf u}^N_\delta (0)=\sum\limits_{r=1}^N\langle{\mathbf u}_0,\phi^r\rangle\phi^r.
\end{cases}
\end{eqnarray}
For any given ${\mathbf u}^N_\delta \in C^1([0,T);X^N),$ one can deduce from the classical theory of transport equation that there exists a classical solution $\rho^N_\delta$ to \eqref{XT1}$_1$ such that
$$\rho^N_\delta \in C^1([0,T);X^N),~~~~\rho^N_\delta (0)\geqslant \delta >0\mbox{ for all }x\in\Omega.$$
Moreover,
\begin{eqnarray}\label{E-78}
  &&\rho^N_\delta (t,x) \geqslant\inf\limits_{x\in\Omega}\rho_0(x)\exp\Big(-\int_0^t\|\operatorname{div}{\mathbf u}^N_\delta \|_{L^\infty(\Omega)}\:ds\Big),\nonumber\\
  &&\rho^N_\delta (t,x) \leqslant\sup\limits_{x\in\Omega}\rho_0(x)\exp\Big(\int_0^t\|\operatorname{div}{\mathbf u}^N_\delta \|_{L^\infty(\Omega)}\:ds\Big)
\end{eqnarray}
for all $(t,x)\in[0,T)\times\Omega,$ implies that there exists a positive constant $C$ such that
	$$\rho^N_\delta (t,x)>C\delta>0\quad\text{for all}\quad(t,x)\in[0,T)\times\Omega.$$
Introduce the operator
\begin{eqnarray*}
  &&{\cal T}:C^{1}([0,T];X^N)\mapsto  C^{1}([0,T];C^{1}(\Omega)),\\
  &&{\mathbf u}^N_\delta \mapsto \rho^N_\delta [{\mathbf u}^N_\delta ].
\end{eqnarray*}
That is, $\rho^N_\delta [{\mathbf u}^N_\delta ]={\cal T}({\mathbf u}^N_\delta ).$ Taking argument in [Chapter 7,  \cite{Feireisl-2004}], one finds that the mapping ${\cal T}$ is continuous.
	
Next, we turn our attention to the approximated problem represented by the integral equation \eqref{XT1}$_2.$ Suppose that $\rho^N_\delta \in C^{1}([0,T];C^{1}(\Omega))$ is given. Following the argument in [Chapter 7, \cite{Feireisl-2004}], we consider a family of linear operators
$${\cal M}[\rho]:X^N\mapsto (X^N)^*,\quad\langle {\cal M}[\rho]{\mathbf v},{\mathbf w}\rangle
	=\int_\Omega\varrho {\mathbf v}\cdot {\mathbf w}\:dx~(\forall {\mathbf v},{\mathbf w}\in X^N),$$
where $(X^N)^*$ stands for the dual space of $X^N.$ It is easy to see that
$$\inf\limits_{\|w\|_{L^2(\Omega)}=1}\langle {\cal M}[\rho]w,w\rangle=\inf_{\|w\|_{L^2(\Omega)}=1}\int\limits_{\Omega}\varrho|w|^2\:dx
	\geqslant\inf\limits_{(t,x)\in[0,T]\times\Omega}\varrho(t,x)\geqslant C\delta>0.$$
So, the operator ${\cal M}$ is invertible and
	$$\|{\cal M}^{-1}[\rho]\|_{\mathcal{L}((X^N)^*,X^N)}\leqslant(C\delta)^{-1},$$
where $\mathcal{L}((X^N)^*,X^N)$ denotes the set of bounded linear mappings from $(X^N)^*$ to $X^N.$  Further, ${\cal M}^{-1}$ is Lipschitz continuous in the sense
\begin{align}\label{2}
  \|{\cal M}^{-1}[\varrho_1]-{\cal M}^{-1}[\varrho_2]\|_{\mathcal{L}((X^N)^*,X^N)}\leqslant C(N,\delta)\|\varrho_1-\varrho_2\|_{L^1(\Omega)}
\end{align}
for $\varrho_1,\varrho_2\in L^1(\Omega)$ such that $\varrho_1,\varrho_2\geqslant C\delta>0.$ Thus, the integral equation \eqref{XT1}$_2$ can be rephrased as
$${\mathbf u}^N_\delta (t)={\cal M}^{-1}[{\cal T}({\mathbf u}^N_\delta )(t)]\Big({\cal M}[\rho_0]{\mathbf u}_0^N+\int_0^t{\cal N}[T({\mathbf u}^N_\delta ),{\mathbf u}^N_\delta (s)]\:ds\Big)\quad\text{in}\quad X^N,$$
where
\begin{eqnarray*}
		&&{\cal N}: X^N\mapsto (X^N)^*\mbox{ with }\\
		&&\langle {\cal N}[{\cal T}({\mathbf u}^N_\delta ),{\mathbf u}^N_\delta ],\eta\rangle
		\equiv \int_\Omega{\cal T}({\mathbf u}^N_\delta ){\mathbf f}\cdot \eta dx
		+\int_\Omega ({\cal T}({\mathbf u}^N_\delta ))^\gamma\mathrm{div}\eta~dx\\
		&&~~~~~~~~~~~~~~~~~~~~~~~~~+\int_\Omega[{\cal T}({\mathbf u}^N_\delta ){\mathbf u}^N_\delta \otimes {\mathbf u}^N_\delta
		-((1+|{\mathbb D}{\mathbf u}^N_\delta |^2 )^{\frac{p(t,x)-2}{2}} {\mathbb D}{\mathbf u}^N_\delta )]:\nabla \eta~dx.
\end{eqnarray*}
Therefore, the integral equation $\eqref{XT1}_2$ can be rewritten as a system of ordinary differential equations. This system in turn fulfills the Carath\'{e}odory conditions and is therefore solvable locally in time, i.e., on a small time interval $[0,T^*)$ for some $T^*<T.$ To ensure solvability for large times, at least for this finite-dimensional problem, we need to establish a first a priori estimate.
	
Since ${\mathbf f}\in L^\infty(0,T;W^{1,2}(\Omega))$ and $\partial_t{\mathbf f}\in{L^2(0,T;L^2(\Omega))},$ one deduces that
\begin{align}\label{E-80}
 \sup\limits_{0\leqslant t\leqslant T}\int_{\Omega}(\rho^N_\delta |{\mathbf u}^N_\delta |^2+(\rho^N_\delta )^\gamma)\:dx
		+\int_{0}^T\int_{\Omega}|{\mathbb D}{\mathbf u}^N_\delta |^{p(t,x)}\:dxdt\leqslant C(T, \rho_0,{\mathbf u}_0,{\mathbf f}).
\end{align}
Moreover, all norms on $X^N$ are equivalent. One finds that the velocity ${\mathbf u}^N_\delta $ is bounded in $L^1(0,T(n);W^{1,\infty}(\Omega)).$ Thus, it follows from \eqref{E-78} that the density $\rho^N_\delta $ is bounded above and below by a constant independent of $T(n).$
As a consequence, ${\mathbf u}^N_\delta $ and $\rho^N_\delta $ remain bounded in $X^N$ independently of $T(n).$  Thus, one can iterate the previous local existence result to construct a solution $(\rho^N_\delta ,{\mathbf u}^N_\delta )\in C^{1}([0,T];X^N)\times C^{1}([0,T];X^N)$ to \eqref{XT1}.
	
\subsection{Estimates independent of $N$ and $\delta$}
	
Following the same arguments as in the previous section, one obtains that for $\frac{12}{5}\leqslant q\leqslant \min\{6(p^--1),2p^-\}$ there exist positive constants $m_{11},~m_{12}, ~m_{13},$ independent of $N$ and $\delta$, such that
\begin{align}
	&\frac{d}{dt}(\|\sqrt{\rho^N_\delta}{\mathbf u}^N_\delta\|_{L^2(\Omega)}^2
     +\|\sqrt{\rho^N_\delta }\partial_t{\mathbf u}_\delta^N\|_{L^2(\Omega)}^2
     +\|\widetilde{D}{\mathbf u}^N_\delta\|_{L^q(\Omega)}^{q}+\|\rho^N_\delta\|_{L^\gamma(\Omega)}
	 +\|\rho^N_\delta \|_{W^{1,3p^-}(\Omega)})\nonumber\\
	&~~~~~+\|{\mathbf u}^N_\delta \|_{W^{2,3p^-}(\Omega)}+{\mathcal I}_\Phi({\mathbf u}^N_\delta )^r+{\mathcal J}_\Phi({\mathbf u}^N_\delta )\nonumber\\
	&\leqslant C(1+\|\sqrt{\rho^N_\delta}{\mathbf u}^N_\delta\|_{L^2(\Omega)}^2+(1+\|f\|_{W^{1,2}(\Omega)}^2)
		\|\rho^N_\delta \|_{W^{1,3p^-}(\Omega)}^{m_{11}}+\|\widetilde{D}{\mathbf u}^N_\delta\|_{L^q(\Omega)}^{m_{12}}\nonumber\\
    &~~~~+\|\sqrt{\rho}\partial_t{\mathbf u}_\delta^N\|_{L^2(\Omega)}^{m_{13}}).
\end{align}
Thus, one deduces from Lemma \ref{local-gronwall} that there exists a small time $T^*\in(0,T)$ such that
\begin{align}\label{ans4}
  &\sup_{0\leqslant t\leqslant T^{*}}\left(\|\sqrt{\rho^N_\delta}{\mathbf u}^N_\delta\|_{L^2(\Omega)}^2
     +\|\sqrt{\rho^N_\delta }\partial_t{\mathbf u}_\delta^N\|_{L^2(\Omega)}^2
		+\|\widetilde{D}{\mathbf u}^N_\delta\|_{L^q(\Omega)}^{q}+\|\rho^N_\delta\|_{L^\gamma(\Omega)}\right.\nonumber\\
  &\left.~~~~~~~~~~+\|\rho^N_\delta \|_{W^{1,3p^-}(\Omega)}\right)
    +\int_{0}^{T^*}(\|{\mathbf u}^N_\delta \|_{W^{2,3p^-}(\Omega)}+{\mathcal I}_\Phi({\mathbf u}^N_\delta )^r
		+{\mathcal J}_\Phi({\mathbf u}^N_\delta ))\:dt\nonumber\\
  &\leqslant C,
\end{align}
where $C$ is independent of $N$ and $\delta.$ So, one integrates \eqref{u2} over $(0,T^*)$ to obtain that
\begin{align}\label{ans5}
\int_{0}^{T^*}\|{\mathbf u}^N_\delta \|_{W^{2,3p^-}(\Omega)}^{p^-}\:dt\leqslant C,
\end{align}
where $C$ is independent of $N$ and $\delta.$ Recalling the definition of ${\mathcal I}_\Phi({\mathbf u}^N_\delta )$ and
${\mathcal J}_\Phi({\mathbf u}^N_\delta ),$ one deduces from Lemma \ref{ine} that the approximation solutions
$(\rho^N_\delta ,{\mathbf u}^N_\delta )$ have the following properties
\begin{eqnarray}\label{N}
\begin{cases}
			\|\sqrt{\rho^N_\delta }{\mathbf u}^N_\delta \|_{L^\infty (0,T;L^2(\Omega))}\leqslant C,
			&\|\rho^N_\delta \|_{L^\infty (0,T;L^\gamma(\Omega))}\leqslant C,\\
			\|\rho^N_\delta\|_{L^\infty(0,T^*;W^{1,3p^-}(\Omega))}\leqslant C,
			&\|\uu^N_\delta\|_{L^\infty(0,T^*;W^{1,\frac{12}{5}}(\Omega))}\leqslant C,\\
			\|{\mathbf u}^N_\delta \|_{L^{p^-} (0,T^*;W^{2,3p^-}(\Omega))}\leqslant C,
			&\|{\mathbf u}^N_\delta \|_{L^{\frac{p^-(5p^--6)}{2-p^-}}(0,T^*; W^{2,\frac{3p^-}{p^-+1}}(\Omega))}\leqslant C,\\
			\|\sqrt{\rho^N_\delta }\partial_t{\mathbf u}_\delta^N\|_{L^{\infty}(0, T^*; L^2(\Omega))}\leqslant C,
			&\|\partial_t{\mathbf u}_\delta^N\|_{L^{p^-}(0, T^*; W^{1,\frac{3p^-}{p^-+1}}(\Omega))}\leqslant C,
\end{cases}
\end{eqnarray}
where $C$ is independent of $N$ and $\delta$.
	
\subsection{Passage to the limit}
	
In this subsection, we are ready to pass the limit for $N\rightarrow\infty$ and $\delta\rightarrow 0^+$ in the sequence of approximate solution
$(\rho^N_\delta ,{\mathbf u}^N_\delta ).$ We first consider the limit for $N\rightarrow \infty.$ We can pick a subsequence (still denoted by $(\rho^N_\delta ,{\mathbf u}^N_\delta )$) with the help of the Aubin-Lions lemma such that
\begin{equation*}
\begin{array}{ll}
\rho_\delta^N\overset{*}\rightharpoonup\rho_\delta\mbox{ in }L^{\infty}(0,T^*; W^{1,3p^-}(\Omega)),
&{\mathbf u}^N_\delta\overset{*}\rightharpoonup{\mathbf u}_\delta\mbox{ in } L^\infty(0,T^*; W^{1,\frac{12}{5}}(\Omega)),\\
{\mathbf u}^N_\delta\rightharpoonup{\mathbf u}_\delta\mbox{ in } L^{\frac{p^-(5p^--6)}{2-p^-}}(0,T^*; W^{2,\frac{3p^-}{p^-+1}}(\Omega)),
&{\mathbf u}^N_\delta\rightharpoonup{\mathbf u}_\delta\mbox{ in } L^{p^-}(0, T^*; W^{2,3p^-}(\Omega)),\\
\sqrt{\rho_\delta^N}\partial_t{\mathbf u}_\delta^N\rightharpoonup\sqrt{\rho}\partial_t{\mathbf u}_\delta\mbox{ in } L^{\infty}(0, T^*; L^2(\Omega)),&\\
\partial_t{\mathbf u}_\delta^N\rightharpoonup\partial_t{\mathbf u}_\delta\mbox{ in } L^{p^-}(0, T^*; W^{1,\frac{3p^-}{p^-+1}}(\Omega)).&
\end{array}
\end{equation*}
These convergence properties are sufficient for passing to the limit in $N.$ Let us detail the convergence of the extra term that arises from the nonlinear part
$${\mathbf S}({\mathbb D}{\uu})=(1+|{\mathbb D}{\mathbf u}|^2)^{\frac{p(t,x)-2}{2}}{\mathbb D}{\mathbf u}.$$
On one hand, it is deduced from \eqref{N} that
$$|{\mathbf S}({\mathbb D}{\uu_\delta^N}(x))|\leqslant C\quad \mbox{for all } x\in\Omega.$$
On the other hand,
$$\mathbb{D}\uu_\delta^N\to\mathbb{ D}\uu_\delta\quad \mbox{a.e. in }(0,T)\times\Omega.$$
Thus, one obtains from the continuity of $\mathbf{S}$ that
$${\mathbf S}({\mathbb D}{\uu}_\delta^N)\to{\mathbf S}({\mathbb D}{\uu}_\delta)\quad \mbox{a.e. in }(0,T)\times\Omega.$$
In addition, it follows from the lower semi-continuity of norm that $(\rho_\delta,{\mathbf u}_\delta)$ satisfies the following estimate
\begin{align}\label{E-87}
	&\sup_{0\leqslant t\leqslant T^{*}}\left(\|\sqrt{\rho_\delta}{\mathbf u}_\delta\|_{L^2(\Omega)}^2
		+\|\sqrt{\rho_\delta}\partial_{t}{\mathbf u}_\delta\|_{L^{2}(\Omega)}^2+\|{\mathbf u}_\delta\|_{W^{1,\frac{12}{5}}(\Omega)}
		+\|\rho_\delta\|_{W^{1,3p^-}(\Omega)}\right)\nonumber\\
	&  +\int_{0}^{T^{*}}\left(\|{\mathbf u}_\delta\|_{W^{2,3p^-}(\Omega)}^{p^-}
		+\|\partial_{t}{\mathbf u}_\delta\|_{W^{1,\frac{3p^-}{p^-+1}}(\Omega)}^{p^-}
		+\|{\mathbf u}_\delta\|_{W^{2,\frac{3p^-}{p^-+1}}(\Omega)}^\frac{p^-(5p^--6)}{2-p^-}\right)\:dt\nonumber\\
		&\leqslant C,
\end{align}
where $C$ is dependent of $T^*, \rho_0,{\mathbf u}_0,{\mathbf f}$ and independent of $\delta.$
	
Our next goal is to carry out the limit process when $\delta\rightarrow 0^+.$ Note that
$$\rho_{\delta,0}\rightarrow \rho_{0}\mbox{ in } W^{1,3p^-}(\Omega)~\mbox{as } \delta\rightarrow 0^+.$$
Further, $T^*$ is independent of $\delta$ and
${\mathbf u}_\delta(0,x)\rightarrow {\mathbf u}(0,x)\in W^{2,p^-}(\Omega)~\mbox{as } \delta\rightarrow 0^+.$
Applying the same argument used in the passage for $N\rightarrow\infty,$ yields a subsequence of $(\rho_\delta,{\mathbf u}_\delta)$ convergent to a strong solution $(\rho,{\mathbf u})$ to the problem \eqref{XT-1}-\eqref{XT-2} satisfying compatibility conditions \eqref{condition}.
	
\subsection{Uniqueness }	
	
Let $(\rho,{\bf u})$ and $(\bar{\rho},\bar{\bf u})$ be two strong solutions to the problem \eqref{XT-1}-\eqref{XT-2} on $(0,T^*)\times \Omega$ satisfying compatibility conditions \eqref{condition}. It is deduced from the regularity of strong solutions stated in Definition \ref{definition} that
\begin{eqnarray}\label{r-1}
\rho, \bar{\rho}\in C([0,T^*),L^r(\Omega))~~(1<r<\infty)\mbox{  and }\uu, \bar{\uu}\in C([0,T^*),L^s(\Omega))~~(1<s<12).\nonumber\\
\end{eqnarray}
One finds that
\begin{eqnarray}\label{E4-4-1}
&&\begin{cases}
	\partial_t(\rho-\bar{\rho})+\divv((\rho-\bar{\rho})\uu+\bar{\rho}(\uu-\bar{\bf u}))=0\mbox{ in  }(0,T^*)\times\Omega,\\
	\partial_t\rho(\uu-\bar{\bf u})+\rho(\uu\cdot \nabla)(\uu-\bar{\bf u})\\
			~~~~~~-\divv((1+|{\mathbb D}{\mathbf u}|^2)^{\frac{p(t,x)-2}{2}}{\mathbb D}{\mathbf u}
		-(1+|{\mathbb D}\bar{\mathbf u}|^2)^{\frac{p(t,x)-2}{2}}{\mathbb D}\bar{\mathbf u})+\nabla (\rho^\gamma-(\bar{\rho})^\gamma)\\
			~~~~~~=(\rho-\bar{\rho})\left({\mathbf f}-\partial_t\bar{\uu}-(\bar{\uu}\cdot\nabla)\bar{\uu}\right)
		-\rho((\uu-\bar{\bf u})\cdot \nabla)\bar{\uu}~~~~\mbox{ in  }(0,T^*)\times\Omega,\\
		\left(\rho-\bar{\rho}\right)(x, 0) = 0\mbox{ and } \left(\rho {\mathbf u}-\bar{\rho}\bar {\mathbf u}\right)(x, 0)=0\mbox{ in  } x\in\Omega.
\end{cases}\nonumber\\
\end{eqnarray}
Taking $\uu-\bar{\bf u}$ as a test function of $\eqref{E4-4-1}_2$ and integrating over $\Omega,$ one obtains that
\begin{align}\label{E4-4-2}
&\frac{1}{2}\frac{d}{dt}\|\sqrt{\rho}(\uu-\bar{\bf u})\|_{L^2(\Omega)}^2\nonumber\\
& +\int_\Omega((1+|{\mathbb D}{\mathbf u}|^2)^{\frac{p(t,x)-2}{2}}{\mathbb D}{\mathbf u}
		-(1+|{\mathbb D}\bar{\mathbf u}|^2)^{\frac{p(t,x)-2}{2}}{\mathbb D}\bar{\mathbf u}):{\mathbb D}(\uu-\bar{\bf u})~dx\nonumber\\
&=\int_\Omega(\rho-\bar{\rho}){\mathbf f}\cdot(\uu-\bar{\bf u})~dx
  -\int_\Omega(\rho-\bar{\rho})\partial_t\bar{\uu}\cdot(\uu-\bar{\bf u})~dx\nonumber\\
&-\int_\Omega(\rho-\bar{\rho})(\bar{\uu}\cdot\nabla)\bar{\uu}\cdot(\uu-\bar{\bf u})~dx-\int_\Omega\rho((\uu-\bar{\bf u})\cdot \nabla)\bar{\uu}\cdot(\uu-\bar{\bf u})~dx\nonumber\\
&-\int_\Omega\nabla (\rho^\gamma-(\bar{\rho})^\gamma)\cdot(\uu-\bar{\bf u})~dx\nonumber\\
&\triangleq\sum_{i=1}^{5}K_i.
\end{align}
Since $p(x,t)>\frac{7}{5},$ there exists a constant $q\in (\frac{7}{4},p^{-}]$ such that
	$$\frac{2-p(x,t)}{2}\cdot \frac{2q}{2-q}<\frac{12}{5}.$$
So, $\|(\widetilde{D}{\mathbf u})^{\frac{2-p(x,t)}{2}}+(\widetilde{D}\bar{\mathbf u})^{\frac{2-p(x,t)}{2}}\|_{L^{\frac{2q}{2-q}}(\Omega)}>0$  and $$\|(\widetilde{D}{\mathbf u})^{\frac{2-p(x,t)}{2}}+(\widetilde{D}\bar{\mathbf u})^{\frac{2-p(x,t)}{2}}\|_{L^{\frac{2q}{2-q}}(\Omega)}\in L^\infty(0,T^*).$$
Using Lemma \ref{L9}, one gets that
\begin{eqnarray*}
&&\|{\mathbb D}(\uu-\bar{\bf u})\|_{L^q(\Omega)}\\
&&\leqslant C\left(\int_\Omega((1+|{\mathbb D}{\mathbf u}|^2)^{\frac{p(t,x)-2}{2}}{\mathbb D}{\mathbf u}
	-(1+|{\mathbb D}\bar{\mathbf u}|^2)^{\frac{p(t,x)-2}{2}}{\mathbb D}\bar{\mathbf u}):{\mathbb D}(\uu-\bar{\bf u})~dx\right)^\frac{1}{2}.
\end{eqnarray*}
The Korn's inequality implies that
$$\|\uu-\bar{\bf u}\|_{L^\frac{3q}{3-q}(\Omega)}\leqslant C\|\uu-\bar{\bf u}\|_{W^{1,q}(\Omega)}\leqslant C\|{\mathbb D}(\uu-\bar{\bf u})\|_{L^q(\Omega)}.$$
Hence,
\begin{eqnarray}\label{E4-4-3}
&& \int_\Omega((1+|{\mathbb D}{\mathbf u}|^2)^{\frac{p(t,x)-2}{2}}{\mathbb D}{\mathbf u}
		-(1+|{\mathbb D}\bar{\mathbf u}|^2)^{\frac{p(t,x)-2}{2}}{\mathbb D}\bar{\mathbf u}):{\mathbb D}(\uu-\bar{\bf u})~dx\nonumber\\
&&\geqslant C\|\uu-\bar{\bf u}\|_{L^\frac{3q}{3-q}(\Omega)}^2
\end{eqnarray}
	and
\begin{eqnarray}\label{E4-4-3-1}
&&\int_\Omega((1+|{\mathbb D}{\mathbf u}|^2)^{\frac{p(t,x)-2}{2}}{\mathbb D}{\mathbf u}
		-(1+|{\mathbb D}\bar{\mathbf u}|^2)^{\frac{p(t,x)-2}{2}}{\mathbb D}\bar{\mathbf u}):{\mathbb D}(\uu-\bar{\bf u})~dx\nonumber\\
&&\geqslant C\|\nabla(\uu-\bar{\bf u})\|_{L^q(\Omega)}^2.
\end{eqnarray}
Next, each term on right hand of \eqref{E4-4-2} is estimated as follows
\begin{align}
|K_1|&=|\int_\Omega(\rho-\bar{\rho}){\mathbf f}\cdot(\uu-\bar{\bf u})~dx|\nonumber\\
	&\leqslant \|\rho-\bar{\rho}\|_{L^\frac{6q}{5q-6}(\Omega)}\|f\|_{L^2(\Omega)}\|\uu-\bar{\uu}\|_{L^\frac{3q}{3-q}(\Omega)}\nonumber\\
	&\leqslant \varepsilon\|\uu-\bar{\uu}\|_{L^\frac{3q}{3-q}(\Omega)}^2+
       C(\varepsilon)\|\rho-\bar{\rho}\|_{L^\frac{6q}{5q-6}(\Omega)}^2\|f\|_{L^2(\Omega)}^2,\label{E4-4-4}\\
|K_2|&=|-\int_\Omega(\rho-\bar{\rho})\partial_t\bar{\uu}\cdot(\uu-\bar{\bf u})~dx|\nonumber\\
	&\leqslant \|\rho-\bar{\rho}\|_{L^\frac{6q}{5q-6}(\Omega)}\|\partial_t{\uu}\|_{L^2(\Omega)}
       \|\uu-\bar{\uu}\|_{L^\frac{3q}{3-q}(\Omega)}\nonumber\\
	&\leqslant \varepsilon\|\uu-\bar{\uu}\|_{L^\frac{3q}{3-q}(\Omega)}^2+
         C(\varepsilon)\|\rho-\bar{\rho}\|_{L^\frac{6q}{5q-6}(\Omega)}^{p^-}\|\partial_t{\uu}\|_{L^2(\Omega)}^{p^-}
          \|\uu-\bar{\uu}\|_{L^\frac{3q}{3-q}(\Omega)}^{2-p^{-}}\nonumber\\
	&\leqslant \varepsilon\|\uu-\bar{\uu}\|_{L^\frac{3q}{3-q}(\Omega)}^2
       + C(\varepsilon)\|\rho-\bar{\rho}\|_{L^\frac{6q}{5q-6}(\Omega)}^{p^-}\|\partial_t{\uu}\|_{L^2(\Omega)}^{p^-}
   \left(\|\uu\|_{L^{12}(\Omega)}^{2-p^{-}}+\|\bar{\uu}\|_{L^{12}(\Omega)}^{2-p^{-}}\right),\label{E4-4-5}\\
|K_3|&=|-\int_\Omega(\rho-\bar{\rho})(\bar{\uu}\cdot\nabla)\bar{\uu}\cdot(\uu-\bar{\bf u})dx|\nonumber\\
	 &\leqslant \|\rho-\bar{\rho}\|_{L^\frac{6q}{5q-6}(\Omega)}\|\bar{\uu}\|_{L^{12}(\Omega)}\|\nabla\bar{\uu}\|_{L^\frac{12}{5}(\Omega)}
		\|\uu-\bar{\uu}\|_{L^\frac{3q}{3-q}(\Omega)}\nonumber\\
	 &\leqslant \|\rho-\bar{\rho}\|_{L^\frac{6q}{5q-6}(\Omega)}\|(\bar{\uu}\cdot\nabla)\bar{\uu}\|_{L^2(\Omega)}
		\|\uu-\bar{\uu}\|_{L^\frac{3q}{3-q}(\Omega)}\label{E4-4-6}\\
	&\leqslant\varepsilon\|\uu-\bar{\uu}\|_{L^\frac{3q}{3-q}(\Omega)}^2+
        C(\varepsilon)\|\rho-\bar{\rho}\|_{L^\frac{6q}{5q-6}(\Omega)}^2\|\bar{\uu}\|_{L^{12}(\Omega)}^2
        \|\nabla\bar{\uu}\|_{L^\frac{12}{5}(\Omega)}^2\nonumber\\
	&\leqslant\varepsilon\|\uu-\bar{\uu}\|_{L^\frac{3q}{3-q}(\Omega)}^2
       + C(\varepsilon)\|\rho-\bar{\rho}\|_{L^\frac{6q}{5q-6}(\Omega)}^2\|\bar{\uu}\|_{W^{1,\frac{12}{5}}(\Omega)}^4,\nonumber\\
|K_4|&=|-\int_\Omega\rho((\uu-\bar{\bf u})\cdot \nabla)\bar{\uu}\cdot(\uu-\bar{\bf u})dx|\leqslant
		\|\sqrt{\rho}(\uu-\bar{\bf u})\|_{L^2(\Omega)}^2\|\nabla\bar{\uu}\|_{L^\infty(\Omega)},\label{E4-4-6}\\
|K_5|&=|-\int_\Omega\nabla (\rho^\gamma-(\bar{\rho})^\gamma)\cdot(\uu-\bar{\bf u})dx|\nonumber\\
		&=|\int_\Omega(\rho^\gamma-(\bar{\rho})^\gamma)\divv(\uu-\bar{\bf u})dx|\nonumber\\
		&\leqslant C\|\rho^\gamma-\bar{\rho}^\gamma\|_{L^\frac{q}{q-1}(\Omega)}\|\nabla(\uu-\bar{\uu})\|_{L^q(\Omega)}\nonumber\\
		&\leqslant\varepsilon\|\nabla(\uu-\bar{\uu})\|_{L^q(\Omega)}^2
        + C(\varepsilon)\|\rho^\gamma-\bar{\rho}^\gamma\|_{L^\frac{q}{q-1}(\Omega)}^{2},\label{E4-4-7}
\end{align}
where $\varepsilon\in(0,1)$ is any fixed. Taking suitable small $\varepsilon\in(0,1),$ one can deduce from \eqref{E4-4-2}-\eqref{E4-4-7} that
\begin{align}\label{E4-4-8}
	&\frac{d}{dt}\|\sqrt{\rho}(\uu-\bar{\bf u})\|_{L^2(\Omega)}^2
        +\frac{1}{C}\left(\|\uu-\bar{\bf u}\|_{L^\frac{3q}{3-q}(\Omega)}^2+\|\nabla(\uu-\bar{\bf u})\|_{L^q(\Omega)}^2\right)\\
	&\leqslant  \|\sqrt{\rho}(\uu-\bar{\bf u})\|_{L^2(\Omega)}^2\|\nabla\bar{\uu}\|_{L^\infty(\Omega)}
        +C\|\rho^\gamma-\bar{\rho}^\gamma\|_{L^\frac{q}{q-1}(\Omega)}^{2}\nonumber\\
	&~+C\|\rho-\bar{\rho}\|_{L^\frac{6q}{5q-6}(\Omega)}^{p^-}\left(\|f\|_{L^2(\Omega)}^2+\|\bar{\uu}\|_{W^{1,\frac{12}{5}}(\Omega)}^4\right)\nonumber\\
&~+C\|\rho-\bar{\rho}\|_{L^\frac{6q}{5q-6}(\Omega)}^{p^-}\|\partial_t{\uu}\|_{L^2(\Omega)}^{p^-}
   \left(\|\uu\|_{L^{12}(\Omega)}^{2-p^{-}}+\|\bar{\uu}\|_{L^{12}(\Omega)}^{2-p^{-}}\right).\nonumber
\end{align}

On one hand, it is deduced from $\eqref{E4-4-1}_1$ that
\begin{align}\label{E4-4-9}
	&\frac{d}{dt}\int_\Omega|\rho-\bar{\rho}|^\frac{6q}{5q-6}~dx\nonumber\\
    &=-\frac{6q}{5q-6}\int_\Omega\divv((\rho-\bar{\rho})\uu+\bar{\rho}(\uu-\bar{\bf u}))
        (|\rho-\bar{\rho}|^\frac{q+6}{5q-6}\mbox{sign}(\rho-\bar{\rho}))~dx\nonumber\\
	&\leqslant  C\left|\int_\Omega|\rho-\bar{\rho}|^\frac{6q}{5q-6}\divv\uu ~dx\right|
      +C\left|\int_\Omega\nabla\bar{\rho}\cdot(\uu-\bar{\bf u}) (|\rho-\bar{\rho}|^\frac{q+6}{5q-6})\mbox{sign}(\rho-\bar{\rho})~dx\right|\nonumber\\
	&~~~+C\left|\int_\Omega\bar{\rho}\divv(\uu-\bar{\bf u}) ((|\rho-\bar{\rho}|^\frac{q+6}{5q-6}\mbox{sign}(\rho-\bar{\rho}))~dx\right|\nonumber\\
	&\triangleq  L_{1}+L_{2}+L_{3}.
\end{align}
Each term on the right hand of \eqref{E4-4-9} is estimated as follows
\begin{align*}
|L_1|&=\left|\int_\Omega|\rho-\bar{\rho}|^\frac{6q}{5q-6}\divv\uu ~dx\right|
      \leqslant C\|\rho-\bar{\rho}\|_{L^\frac{6q}{5q-6}(\Omega)}^\frac{6q}{5q-6}\|\nabla\uu\|_{L^\infty(\Omega)},\\
|L_2|&=\left|\int_\Omega\nabla\bar{\rho}\cdot(\uu-\bar{\bf u}) ((|\rho-\bar{\rho}|^\frac{q+6}{5q-6})\mbox{sign}(\rho-\bar{\rho}))~dx\right|\\
	&\leqslant \|\nabla\bar{\rho}\|_{L^{3q}(\Omega)}\|\uu-\bar{\bf u}\|_{L^{\frac{3q}{3-q}}(\Omega)}
      \left(\int_\Omega|\rho-\bar{\rho}|^{\frac{q+6}{5q-6}\cdot\frac{3q}{4(q-1)}}~dx\right)^\frac{4(q-1)}{3q}\\
	&\leqslant C\|\nabla\bar{\rho}\|_{L^{3q}(\Omega)}\|\uu-\bar{\bf u}\|_{L^{\frac{3q}{3-q}}(\Omega)}
      \left(\int_\Omega|\rho-\bar{\rho}|^{\frac{6q}{5q-6}}~dx\right)^\frac{4(q-1)}{3q},\\
|L_3|&=\left|\int_\Omega\bar{\rho}\divv(\uu-\bar{\bf u}) ((|\rho-\bar{\rho}|^\frac{q+6}{5q-6}\mbox{sign}(\rho-\bar{\rho}))~dx\right|\\
	&\leqslant \|\bar{\rho}\|_{L^{\infty}(\Omega)}\|\divv(\uu-\bar{\bf u})\|_{L^{q}(\Omega)}
      \left(\int_\Omega|\rho-\bar{\rho}|^{\frac{q+6}{5q-6}\cdot\frac{q}{q-1}}~dx\right)^\frac{q-1}{q}\\
     &\leqslant C\|\bar{\rho}\|_{L^{\infty}(\Omega)}
		\|\nabla(\uu-\bar{\bf u})\|_{L^{q}(\Omega)}\left(\int_\Omega|\rho-\bar{\rho}|^{\frac{6q}{5q-6}}~dx\right)^\frac{q-1}{q}.
\end{align*}
Combining \eqref{E4-4-9}, one obtains that
\begin{align}\label{E4-4-10}
	\frac{d}{dt}\int_\Omega|\rho-\bar{\rho}|^\frac{6q}{5q-6}~dx
     &\leqslant   C\|\rho-\bar{\rho}\|_{L^\frac{6q}{5q-6}(\Omega)}^\frac{6q}{5q-6}\|\nabla\uu\|_{L^\infty(\Omega)}\nonumber\\
	&~~~+C\left(\int_\Omega|\rho-\bar{\rho}|^{\frac{6q}{5q-6}}~dx\right)^\frac{q-1}{q}\|\bar{\rho}\|_{L^{\infty}(\Omega)}
		\|\nabla(\uu-\bar{\bf u})\|_{L^{q}(\Omega)}\nonumber\\
	&~~~+C\|\nabla\bar{\rho}\|_{L^{3q}(\Omega)}\|\uu-\bar{\bf u}\|_{L^{\frac{3q}{3-q}}(\Omega)}
        \left(\int_\Omega|\rho-\bar{\rho}|^{\frac{6q}{5q-6}}~dx\right)^\frac{4(q-1)}{3q}.
\end{align}

On the other hand, it is deduced from $\eqref{E4-4-1}_1$ that
\begin{align}\label{E4-4-11}
	&\frac{d}{dt}\int_\Omega|\rho^\gamma-\bar{\rho}^\gamma|^\frac{q}{q-1}dx\nonumber\\
	&\leqslant C\left|\int_\Omega|\rho^\gamma-\bar{\rho}^\gamma|^\frac{q}{q-1}\divv\uu dx\right|
      +C\left|\int_\Omega\nabla\bar{\rho}^\gamma\cdot(\uu-\bar{\bf u}) |\rho^\gamma-\bar{\rho}^\gamma|^\frac{1}{q-1}\mbox{sign}(\rho-\bar{\rho})dx\right|\nonumber\\
	&~~~+C\left|\int_\Omega\bar{\rho}^\gamma\divv(\uu-\bar{\bf u}|\rho^\gamma-\bar{\rho}^\gamma|^\frac{1}{q-1}\mbox{sign}(\rho-\bar{\rho})dx\right|\nonumber\\
	&=Y_1+Y_2+Y_3.
\end{align}
Taking the similar argument for each term on the right hand of \eqref{E4-4-9}, one estimates each term on the right hand of \eqref{E4-4-11} as follows
\begin{align*}
|Y_1|&=\left|\int_\Omega|\rho^\gamma-\bar{\rho}^\gamma|^\frac{q}{q-1}\divv\uu dx\right|
        \leqslant C\|\rho^\gamma-\bar{\rho}^\gamma\|_{L^\frac{q}{q-1}(\Omega)}^\frac{q}{q-1}\|\nabla\uu\|_{L^\infty(\Omega)},\\
|Y_2|&=\left|\int_\Omega\nabla\bar{\rho}^\gamma\cdot(\uu-\bar{\bf u}) |\rho^\gamma-\bar{\rho}^\gamma|^\frac{1}{q-1}\mbox{sign}(\rho-\bar{\rho}))dx\right|\\
	 &\leqslant C\|\nabla\bar{\rho}^\gamma\|_{L^{3q}(\Omega)}\|\uu-\bar{\bf u}\|_{L^{\frac{3q}{3-q}}(\Omega)}
           \left(\int_\Omega|\rho^\gamma-\bar{\rho}^\gamma|^\frac{q}{q-1}dx\right)^\frac{1}{q},\\
|Y_3|&=\left|\int_\Omega\bar{\rho}^\gamma\divv(\uu-\bar{\bf u})|\rho^\gamma-\bar{\rho}^\gamma|^\frac{1}{q-1}\mbox{sign}(\rho-\bar{\rho})dx\right|\\
	 &\leqslant \|\bar{\rho}\|_{L^{\infty}(\Omega)}\|\divv(\uu-\bar{\bf u})\|_{L^{q}(\Omega)}
       \left(\int_\Omega|\rho^\gamma-\bar{\rho}^\gamma|^qdx\right)^\frac{q-1}{q}\\
	 &\leqslant C\|\bar{\rho}\|_{L^{\infty}(\Omega)}\|\nabla(\uu-\bar{\bf u})\|_{L^{q}(\Omega)}
        \left(\left(\int_\Omega|\rho^\gamma-\bar{\rho}^\gamma|^{q\cdot\frac{1}{q-1}}dx\right)^{q-1}\cdot
		\left(\int_\Omega1^{\frac{1}{2-q}}dx\right)^{2-q}\right)^\frac{q-1}{q}\\
	 &\leqslant C\|\bar{\rho}\|_{L^{\infty}(\Omega)}
     \|\nabla(\uu-\bar{\bf u})\|_{L^{q}(\Omega)}\left(\int_\Omega|\rho^\gamma-\bar{\rho}^\gamma|^{\frac{q}{q-1}}dx\right)^\frac{(q-1)^2}{q}.
\end{align*}
Combining \eqref{E4-4-11}, one obtains that
\begin{eqnarray}\label{E4-4-11}
  &\frac{d}{dt}\int_\Omega|\rho^\gamma-\bar{\rho}^\gamma|^\frac{q}{q-1}dx&\leqslant C\|\rho^\gamma-\bar{\rho}^\gamma\|_{L^\frac{q}{q-1}(\Omega)}^\frac{q}{q-1}\|\nabla\uu\|_{L^\infty(\Omega)}\\
  &&~~~+C\|\nabla\bar{\rho}^\gamma\|_{L^{3q}(\Omega)}\|\uu-\bar{\bf u}\|_{L^{\frac{3q}{3-q}}(\Omega)}
        \left(\int_\Omega|\rho^\gamma-\bar{\rho}^\gamma|^\frac{q}{q-1}dx\right)^\frac{1}{q}\nonumber\\
  &&~~~+C\|\bar{\rho}\|_{L^{\infty}(\Omega)}\|\nabla(\uu-\bar{\bf u})\|_{L^{q}(\Omega)}
        \left(\int_\Omega|\rho^\gamma-\bar{\rho}^\gamma|^{\frac{q}{q-1}}dx\right)^\frac{(q-1)^2}{q}.\nonumber
\end{eqnarray}
	
%%%%%%%%%%%%%%%%%%%%%%%%%%%%%%%%%%%%%%%%%%%%%%%%%%%%%%%%%%%%%%%%%%%%%%%%%%%%%%%%%%%%%%%%%%%%%%%%%%%%%%%%%%%%%%%%%%%%%%%%%%%%%%%%%%%%%%%%%%%%%%%%%%	

{\bf Claim.} $\rho=\bar{\rho}$ a.e. on $(0,T^*)\times\Omega.$
	
Assume that $\rho\neq\bar{\rho}$ a.e. on $(0,T^*)\times\Omega.$ Then
\begin{eqnarray}\label{E4-4-12}
	\int_\Omega|\rho-\bar{\rho}|^{\frac{6p^{-}}{5p^{-}-6}}dx=c_1(t)>0\mbox{   and  }
	\int_\Omega|\rho^\gamma-\bar{\rho}^\gamma|^{\frac{q}{q-1}}dx=c_2(t)>0
\end{eqnarray}
with $c_1(t)\in C[0,T^*)$ and $c_2(t)\in C[0,T^*).$ Note that $q\in (\frac{7}{4},p^-].$ Then one takes $q=p^-$ and deduces from \eqref{E4-4-10} and \eqref{E4-4-12} that
\begin{align}\label{E4-4-13}
  &\frac{d}{dt}\|\rho-\bar{\rho}\|_{L^\frac{6p^{-}}{5p^{-}-6}(\Omega)}^{p^{-}}\nonumber\\
	          &\leqslant  C\|\rho-\bar{\rho}\|_{L^\frac{6p^{-}}{5p^{-}-6}(\Omega)}^{p^{-}}\|\nabla\uu\|_{L^\infty(\Omega)}\nonumber\\
		       &~~~ +C\left(\int_\Omega|\rho-\bar{\rho}|^{\frac{6p^{-}}{5p^{-}-6}}dx\right)^\frac{5p^{-}-6}{6}\|\bar{\rho}\|_{L^{\infty}(\Omega)}
		         \|\nabla(\uu-\bar{\bf u})\|_{L^{p^{-}}(\Omega)}\nonumber\\
		     &~~~+C\|\nabla\bar{\rho}\|_{L^{3p^{-}}(\Omega)}\|\uu-\bar{\bf u}\|_{L^{\frac{3p^{-}}{3-p^{-}}}(\Omega)}
                \left(\int_\Omega|\rho-\bar{\rho}|^{\frac{6p^{-}}{5p^{-}-6}}dx\right)^\frac{2(p^--1)(p^-+1)}{3p^-}\nonumber\\
	         &\leqslant  \varepsilon \|\uu-\bar{\bf u}\|_{L^{\frac{3p^{-}}{3-p^{-}}}(\Omega)}^2
               +\varepsilon\|\nabla(\uu-\bar{\bf u})\|_{L^{p^{-}}(\Omega)}^2+
		       C\|\rho-\bar{\rho}\|_{L^\frac{6p^{-}}{5p^{-}-6}(\Omega)}^{p^{-}}\|\nabla\uu\|_{L^\infty(\Omega)}\nonumber\\
		     &~~~+C\|\nabla\bar{\rho}\|_{L^{3p^{-}}(\Omega)}^2
               \left(\int_\Omega|\rho-\bar{\rho}|^{\frac{6p^{-}}{5p^{-}-6}}dx\right)^\frac{2(p^--1)(3p^-+4)}{3p^-}\nonumber\\
             & +C\left(\int_\Omega|\rho-\bar{\rho}|^{\frac{6p^{-}}{5p^{-}-6}}dx\right)^\frac{2(p^--1)(p^-+1)}{3p^-}
                  \|\bar{\rho}\|_{L^{\infty}(\Omega)}^2\nonumber\\
	         &\leqslant  \varepsilon \|\uu-\bar{\bf u}\|_{L^{\frac{3p^{-}}{3-p^{-}}}(\Omega)}^2
                +\varepsilon\|\nabla(\uu-\bar{\bf u})\|_{L^{p^{-}}(\Omega)}^2+
                C\|\rho-\bar{\rho}\|_{L^\frac{6p^{-}}{5p^{-}-6}(\Omega)}^{p^{-}}\|\nabla\uu\|_{L^\infty(\Omega)}\nonumber\\
		    &~~~+C\|\nabla\bar{\rho}\|_{L^{3p^{-}}(\Omega)}^2\left(\int_\Omega|\rho-\bar{\rho}|^{\frac{6p^{-}}{5p^{-}-6}}dx\right)^\frac{6p^-6}{6}
		         \left(\int_\Omega|\rho-\bar{\rho}|^{\frac{6p^{-}}{5p^{-}-6}}dx\right)^{-\frac{5(p^{-})^2-3p^{-}-16}{6p^-}}\nonumber\\
		    &~~~+ C\left(\int_\Omega|\rho-\bar{\rho}|^{\frac{6p^{-}}{5p^{-}-6}}dx\right)^\frac{6p^-6}{6}		
                 \left(\int_\Omega|\rho-\bar{\rho}|^{\frac{6p^{-}}{5p^{-}-6}}dx\right)^{-\frac{5(p^{-})^2-6p^{-}-12}{6p^-}}
                  \|\bar{\rho}\|_{L^{\infty}(\Omega)}^2\nonumber\\
             &\leqslant  \varepsilon \|\uu-\bar{\bf u}\|_{L^{\frac{3p^{-}}{3-p^{-}}}(\Omega)}^2
                +\varepsilon\|\nabla(\uu-\bar{\bf u})\|_{L^{p^{-}}(\Omega)}^2
               +C\|\rho-\bar{\rho}\|_{L^\frac{6p^{-}}{5p^{-}-6}(\Omega)}^{p^{-}}\|\nabla\uu\|_{L^\infty(\Omega)}\nonumber\\
		     &~~~+	C\|\rho-\bar{\rho}\|_{L^\frac{6p^{-}}{5p^{-}-6}(\Omega)}^{p^{-}}
\|\nabla\bar{\rho}\|_{L^{3p^{-}}(\Omega)}^2 c_1^{-\frac{5(p^{-})^2-3p^{-}-16}{6p^-}}(t)\nonumber\\
&~~~+C\|\rho-\bar{\rho}\|_{L^\frac{6p^{-}}{5p^{-}-6}(\Omega)}^{p^{-}}c_1^{-\frac{5(p^{-})^2-6p^{-}-12}{6p^-}}(t)\|\bar{\rho}\|_{L^{\infty}(\Omega)}^2
\end{align}
holds for any fixed $\varepsilon\in(0,1).$ Also, one uses \eqref{E4-4-11} and \eqref{E4-4-12} to arrive that
\begin{align}\label{E4-4-14}
 &\frac{d}{dt}\|\rho^\gamma-\bar{\rho}^\gamma\|_{L^\frac{p^-}{p^--1}(\Omega)}^{2}\nonumber\\
 &\leqslant\varepsilon \|\uu-\bar{\bf u}\|_{L^{\frac{3p^{-}}{3-p^{-}}}(\Omega)}^2
           +\varepsilon\|\nabla(\uu-\bar{\bf u})\|_{L^{p^{-}}(\Omega)}^2\nonumber\\
 &~~~+ C\|\rho^\gamma-\bar{\rho}^\gamma\|_{L^\frac{p^-}{p^--1}(\Omega)}^{2}
  \left(\|\nabla\uu\|_{L^\infty(\Omega)}+\|\nabla\bar{\uu}\|_{L^\infty(\Omega)} +\|\nabla\bar{\rho}^\gamma\|_{L^{3p^{-}}(\Omega)}^2\right)\nonumber\\
 &~~~+C\|\rho^\gamma-\bar{\rho}^\gamma\|_{L^\frac{p^-}{p^--1}(\Omega)}^{2}\|{\rho}\|_{L^{\infty}(\Omega)}^2 c_2^{-\frac{2-p^-}{p^-}}(t)
\end{align}
holds for any fixed $\varepsilon\in(0,1).$ Now, one adds \eqref{E4-4-8}, \eqref{E4-4-13} and \eqref{E4-4-14} to get the following inequality
\begin{align}\label{E4-4-15}
    &\frac{d}{dt}\left(\|\sqrt{\rho}(\uu-\bar{\bf u})\|_{L^2(\Omega)}^2+\|\rho-\bar{\rho}\|_{L^\frac{6p^{-}}{5p^{-}-6}(\Omega)}^{p^{-}}
		+\|\rho^\gamma-\bar{\rho}^\gamma\|_{L^\frac{p^-}{p^--1}(\Omega)}^{2}\right)\nonumber\\
	&~~~+\frac{1}{C}\left(\|\uu-\bar{\bf u}\|_{L^\frac{3p^-}{3-p^-}(\Omega)}^2+\|\nabla(\uu-\bar{\bf u})\|_{L^{p^-}(\Omega)}^2\right)\nonumber\\
    &\leqslant  \|\sqrt{\rho}(\uu-\bar{\bf u})\|_{L^2(\Omega)}^2\|\nabla\bar{\uu}\|_{L^\infty(\Omega)}
		+C\|\rho^\gamma-\bar{\rho}^\gamma\|_{L^\frac{p^-}{p^--1}(\Omega)}^{2}\nonumber\\
	&~~~+C\|\rho-\bar{\rho}\|_{L^\frac{6p^{-}}{5p^{-}-6}(\Omega)}^{p^{-}}
		  \left(\|f\|_{L^2(\Omega)}^2+ \|\partial_t{\uu}\|_{L^2(\Omega)}^{p^-}(\|\uu\|_{L^{12}(\Omega)}^{2-p^{-}}
          +\|\bar{\uu}\|_{L^{12}(\Omega)}^{2-p^{-}})+\|\bar{\uu}\|_{W^{1,\frac{12}{5}}(\Omega)}^4\right.\nonumber\\
    &~~~\left. +\|\nabla\uu\|_{L^\infty(\Omega)}+
		c_1^{-\frac{5(p^{-})^2-3p^{-}-16}{6p^-}}(t)\|\nabla\bar{\rho}\|_{L^{3p^{-}}(\Omega)}^2
        +c_1^{-\frac{5(p^{-})^2-6p^{-}-12}{6p^-}}(t)\|\bar{\rho}\|_{L^{\infty}(\Omega)}^2\right)\nonumber\\
	&~~~+ C\|\rho^\gamma-\bar{\rho}^\gamma\|_{L^\frac{p^-}{p^--1}(\Omega)}^{2}
		\left(\|\nabla\uu\|_{L^\infty(\Omega)}+\|\nabla\bar{\uu}\|_{L^\infty(\Omega)}
		+\|\nabla\bar{\rho}^\gamma\|_{L^{3p^{-}}(\Omega)}^2\right)\nonumber\\
    &~~~+C\|\rho^\gamma-\bar{\rho}^\gamma\|_{L^\frac{p^-}{p^--1}(\Omega)}^{2}\|{\rho}\|_{L^{\infty}(\Omega)}^2 c_2^{-\frac{2-p^-}{p^-}}(t)\nonumber\\
    &\leqslant  C\left(\|\sqrt{\rho}(\uu-\bar{\bf u})\|_{L^2(\Omega)}^2+\|\rho-\bar{\rho}\|_{L^\frac{6p^{-}}{5p^{-}-6}(\Omega)}^{p^{-}}
		+\|\rho^\gamma-\bar{\rho}^\gamma\|_{L^\frac{p^-}{p^--1}(\Omega)}^{2}\right)\cdot\nonumber\\
	&~~~\left(1+\|\nabla\bar{\uu}\|_{L^\infty(\Omega)}+\|\nabla{\uu}\|_{L^\infty(\Omega)}+\|f\|_{L^2(\Omega)}^2
           +\|\bar{\uu}\|_{W^{1,\frac{12}{5}}(\Omega)}^4+\|\nabla\bar{\rho}^\gamma\|_{L^{3p^{-}}(\Omega)}^2\right.\nonumber\\
    &~~~\left.+\|\partial_t{\uu}\|_{L^2(\Omega)}^{p^-}\left(\|\uu\|_{L^{12}(\Omega)}^{2-p^{-}}+\|\bar{\uu}\|_{L^{12}(\Omega)}^{2-p^{-}}\right)
    +c_1^{-\frac{5(p^{-})^2-3p^{-}-16}{6p^-}}(t)\|\nabla\bar{\rho}\|_{L^{3p^{-}}(\Omega)}^2\right.\nonumber\\
	&~~~\left.
        +c_1^{-\frac{5(p^{-})^2-6p^{-}-12}{6p^-}}(t)\|\bar{\rho}\|_{L^{\infty}(\Omega)}^2
		+c_2^{-\frac{2-p^-}{p^-}}(t)\|{\rho}\|_{L^{\infty}(\Omega)}^2\right).
\end{align}
Note that
\begin{eqnarray*}
&&1-\frac{5(p^{-})^2-2p^{-}-16}{6p^-}=\frac{16+8p^{-}-5(p^{-})^2}{6p^-}>0,\\
&&1-\frac{5(p^{-})^2-6p^{-}-12}{6p^-}=\frac{12+12p^{-}-5(p^{-})^2}{6p^-}>0\mbox{ and }1-\frac{2-p^-}{p^-}=\frac{2(p^--1)}{p^-}>0,
\end{eqnarray*}
based on the fact that $\frac{7}{5}<p^-\leqslant p(t,x)\leqslant p^+\leqslant 2.$ Hence, one follows from the Gronwall's inequality that
$$\|\sqrt{\rho}(\uu-\bar{\bf u})\|_{L^2(\Omega)}^2+\|\rho-\bar{\rho}\|_{L^\frac{6p^{-}}{5p^{-}-6}(\Omega)}^{p^{-}}
	+\|\rho^\gamma-\bar{\rho}^\gamma\|_{L^\frac{p^-}{p^--1}(\Omega)}^{2}=0\ \mbox{ a.e. on } (0,T^*),$$
which contradicts to \eqref{E4-4-12}. The Claim is proved.

%%%%%%%%%%%%%%%%%%%%%%%%%%%%%%%%%%%%%%%%%%%%%%%%%%%%%%%%%%%%%%%%%%%%%%%%%%%%%%%%%%%%%%%%%%%%%%%%%%%%%%%%%%%%%%%%%%%%%%%%%%%%%%%%%%%%%%%%%%%%%%%%%
	
Further, one deduces from Claim and \eqref{E4-4-8} that
\begin{eqnarray*}
&&\frac{d}{dt}\|\sqrt{\rho}(\uu-\bar{\bf u})\|_{L^2(\Omega)}^2
   +\frac{1}{C}\left(\|\uu-\bar{\bf u}\|_{L^\frac{3p^-}{3-p^-}(\Omega)}^2+\|\nabla(\uu-\bar{\bf u})\|_{L^{p^-}(\Omega)}^2\right)\\
&&\leqslant \|\sqrt{\rho}(\uu-\bar{\bf u})\|_{L^2(\Omega)}^2\|\nabla\bar{\uu}\|_{L^\infty(\Omega)}.
\end{eqnarray*}
Obviously, it is easy to get
	$$\|\sqrt{\rho}(\uu-\bar{\bf u})\|_{L^2(\Omega)}^2=0\mbox{ a.e. on } (0,T^*)$$
and so $\uu=\bar{\uu}$ a.e. on $(0,T^*)\times\Omega.$
	
\section{Proof of Theorem 2 }\label{Sec-2}
	
This section establishes a blow-up criterion for strong solutions to the problem \eqref{XT-1}-\eqref{XT-2} under the compatible conditions \eqref{condition}. Let $(\rho,{\mathbf u})$ be a strong solution to the problem \eqref{XT-1}-\eqref{XT-2} on $[0,T^*)\times \Omega,$ which satisfies the regularity \eqref{1dlg-E6} on the time interval $[0,T^*).$ We proceed by contradiction to complete the proof of Theorem \ref{blow}. Suppose that
\begin{equation}\label{u3}
	\lim\limits_{T\rightarrow T^*}\|\rho\|_{L^\infty(0,T;L^\infty(\Omega))}+\|\nabla{\mathbf u}\|_{L^\infty(0,T;L^3(\Omega))}\leqslant M<\infty,
\end{equation}
then we are going to obtain a contradiction to the maximality of $T^*.$
	
To begin with, one deduces from energy inequality and \eqref{u3} that
\begin{align}\label{E-}
 \sup\limits_{0\leqslant t\leqslant T}\int_{\Omega}(\rho |{\mathbf u} |^2+(\rho )^\gamma)\:dx
		+\int_{0}^{T^*}\int_{\Omega}|{\mathbb D}{\mathbf u} |^{p(t,x)}\:dxdt\leqslant C
\end{align}
holds for any $T\in (0,T^*).$
	
Next, we drive bounds for time and second order derivatives of $\uu$, by taking similar argument in Section \ref{Sec-A}.
	
\begin{lemma}\label{uu3}
Under the condition of Theorem \ref{blow} and \eqref{u3}, there exists a positive constant $C$ such that for $r\in [\max \{1,\frac{5}{5p^--6}\},\frac{5p^--6}{2-p^-}),$
\begin{equation}\label{E-91}
  \sup\limits_{0\leqslant t\leqslant T}(\|\sqrt{\rho}\partial_t{\mathbf u}\|_{L^2(\Omega)}^2+\|\rho\|_{W^{1,3p^-}(\Omega)})
	+\int_{0}^{T}\left({\mathcal I}^r_\Phi({\mathbf u})+{\mathcal J}_\Phi({\mathbf u})+\|{\mathbf u}\|_{W^{2,3p^-}(\Omega)}^{p^-}\right)\:dt
   \leqslant C
\end{equation}
holds for any $T\in(0,\min\{\frac{1}{c_0(\alpha-1)}\left(f(0)+c_1(T^*+\|\partial_t{\bf f}\|_{L^2(0,T^*;L^2()\Omega)})\right)^{-(\alpha-1)},T^*\}),$ where $C$ depends on $T$ continuously, %$max\{1,\frac{5}{5p^--6}\}\leqslant r<\frac{5p^--6}{2-p^-},$
$$\alpha=\max \{\frac{p^-r}{p^--1},\frac{4r(p^--1)}{5p^--6-r(2-p^-)}\}\mbox{ and }f(0)\triangleq\|\sqrt{\rho_0}\partial_t({\mathbf u}_0)\|_{L^2(\Omega)}^2+\|\rho_0\|_{W^{1,3p^-}(\Omega)}.$$
\end{lemma}
	
\begin{proof}
Recalling the definition of \eqref{ID}, taking similar argument to derive \eqref{delu} and \eqref{L1}, and selecting appropriate coefficient, one gets that
\begin{align}\label{b1}
{\mathcal I}_\Phi({\mathbf u})
			&\leqslant \int_\Omega|\nabla\rho\|{\mathbf u}\|\nabla{\mathbf u}|^2~dx
			+\int_\Omega\rho|\nabla{\mathbf u}|^3~dx+\int_\Omega\rho|{\mathbf u}\|\nabla^2{\mathbf u}\|\nabla{\mathbf u}|~dx\nonumber\\
			&~~~+\int_\Omega\rho|\partial_t{\mathbf u}\|\Delta{\mathbf u}|~dx
			+\int_\Omega|\nabla \rho^\gamma\|\Delta{\mathbf u}|dx+\int_{0}^{T}\rho|{\mathbf f}\|\Delta{\mathbf u}|~dx
			+C\int_\Omega(\widetilde{D}{\mathbf u})^{p(t,x)+\sigma}~dx\nonumber\\
			&\leqslant \sum_{i=1}^{7}L_i,
\end{align}
where $\sigma\in (0,1).$ Now, one uses the definition \eqref{JD} and the condition \eqref{u3} to estimate each term on the right-hand side of \eqref{b1} as follows.
\begin{align*}
	L_{1}&\leqslant C\|\nabla\rho\|_{L^{3p^-}(\Omega)}\|{\mathbf u}\|_{L^\infty(\Omega)}\|\nabla{\mathbf u}\|_{L^3(\Omega)}^2
			\leqslant C \|\nabla\rho\|_{L^{3p^-}(\Omega)}\|\nabla{\mathbf u}\|_{L^3(\Omega)}^3\leqslant C\|\nabla\rho\|_{L^{3p^-}(\Omega)},\\
	L_{2}&\leqslant C\|\rho\|_{L^\infty(\Omega)}\|\nabla{\mathbf u}\|_{L^3(\Omega)}^3\leqslant C,\\
	L_{3}&\leqslant C\|\rho\|_{L^\infty(\Omega)}\|{\mathbf u}\|_{L^\infty(\Omega)}\|\nabla{\mathbf u}\|_{L^3(\Omega)}
			\|\nabla^2{\mathbf u}\|_{L^\frac{3p^-}{p^-+1}(\Omega)}\\
		 &\leqslant C\|\nabla{\mathbf u}\|_{L^3(\Omega)}^2({\mathcal I}_\Phi({\mathbf u})+1)^\frac{1}{p^-}\\
		 &\leqslant \varepsilon{\mathcal I}_\Phi({\mathbf u})+C,\\
	L_{4}&\leqslant C\|\rho\|_{L^\infty(\Omega)}^\frac{1}{2}\|\sqrt{\rho}\partial_t{\mathbf u}\|_{L^\frac{3p^-}{2p^--1}(\Omega)}
			\|\Delta{\mathbf u}\|_{L^\frac{3p^-}{p^-+1}(\Omega)}\\
			&\leqslant C\|\sqrt{\rho}\partial_t{\mathbf u}\|_{L^2(\Omega)}^\frac{4p^--4}{3p^--2}
			\|\sqrt{\rho}\partial_t{\mathbf u}\|_{L^{3p^-}(\Omega)}^\frac{2 - p^-}{3p^--2}\|{\mathbf u}\|_{W^{2,\frac{3p^-}{p^-+1}}(\Omega)}\\
			&\leqslant C\|\rho\|_{L^\infty(\Omega)}^\frac{p^-}{3p^--2}\|\sqrt{\rho}\partial_t{\mathbf u}\|_{L^2(\Omega)}^\frac{4p^--4}{3p^--2}
			\|\partial_t{\mathbf u}\|_{W^{1,\frac{3p^-}{p^-+1}}(\Omega)}^\frac{2 - p^-}{3p^--2}({\mathcal I}_\Phi({\mathbf u})+1)^\frac{1}{p^-} \\
			&\leqslant C\|\sqrt{\rho}\partial_t{\mathbf u}\|_{L^2(\Omega)}^\frac{4p^--4}{3p^--2}
			{\mathcal J}^\frac{1}{2}_\Phi({\mathbf u})({\mathcal I}_\Phi({\mathbf u})+1)^\frac{2 - p^-}{2p^-})^\frac{2 - p^-}{3p^--2}
			({\mathcal I}_\Phi({\mathbf u})+1)^\frac{1}{p^-} \\
			&\leqslant C\|\sqrt{\rho}\partial_t{\mathbf u}\|_{L^2(\Omega)}^\frac{4p^--4}{3p^--2}
			({\mathcal I}_\Phi({\mathbf u})+1)^\frac{p^- + 2}{2(3p^--2)}{\mathcal J}_\Phi({\mathbf u})^\frac{2 - p^-}{2(3p^--2)} \\
			&\leqslant \varepsilon{\mathcal I}_\Phi({\mathbf u})
			+C(1+\|\sqrt{\rho}\partial_t{\mathbf u}\|_{L^2(\Omega)}^\frac{8(p^--1)}{5p^--6}{\mathcal J}_\Phi({\mathbf u})^\frac{2-p^-}{5p^--6} ),\\
	L_5&\leqslant C\|\rho\|_{L^\infty(\Omega)}^{\gamma-1}\|\nabla\rho\|_{L^{3p^-}(\Omega)}\|\Delta{\mathbf u}\|_{L^\frac{3p^-}{p^-+1}(\Omega)}\\
			&\leqslant C\|\nabla\rho\|_{L^{3p^-}}({\mathcal I}_\Phi({\mathbf u})+1)^\frac{1}{p^-}\\
			&\leqslant \varepsilon{\mathcal I}_\Phi({\mathbf u})+C(1+\|\nabla\rho\|_{L^{3p^-}(\Omega)}^\frac{p^-}{p^--1}),\\
	L_6&\leqslant C\|\rho\|_{L^\infty(\Omega)}\|{\mathbf f}\|_{L^6(\Omega)}\|\Delta{\mathbf u}\|_{L^\frac{3p^-}{p^-+1}(\Omega)}
			\leqslant C({\mathcal I}_\Phi({\mathbf u})+1)^\frac{1}{p^-}\leqslant \varepsilon{\mathcal I}_\Phi({\mathbf u})+C,\\
	L_7&\leqslant C\|\nabla{\mathbf u}\|_{L^3(\Omega)}^3\leqslant  C,
\end{align*}
where $\varepsilon\in(0,1)$  is arbitrary but fixed. Substituting the above estimates into \eqref{b1}, one gets that
\begin{align*}
{\mathcal I}_\Phi({\mathbf u})
          \leqslant C(1+\|\nabla\rho\|_{L^{3p^-}(\Omega)}^\frac{p^-}{p^--1}
			+\|\sqrt{\rho}\partial_t{\mathbf u}\|_{L^2(\Omega)}^\frac{8(p^--1)}{5p^--6}{\mathcal J}_\Phi({\mathbf u})^\frac{2-p^-}{5p^--6}).
\end{align*}
So, for $r\in [\max \{1,\frac{5}{5p^--6}\},\frac{5p^--6}{2-p^-}),$
\begin{align}\label{31}
{\mathcal I}^r_\Phi({\mathbf u})
	  &\leqslant C(1+\|\nabla\rho\|_{L^{3p^-}(\Omega)}^\frac{p^-r}{p^--1}
       +\|\sqrt{\rho}\partial_t{\mathbf u}\|_{L^2(\Omega)}^\frac{8(p^--1)r}{5p^--6}{\mathcal J}_\Phi({\mathbf u})^\frac{(2-p^-)r}{5p^--6})\nonumber\\
      &\leqslant C(1+\eta{\mathcal J}_\Phi({\mathbf u})+\|\nabla\rho\|_{L^{3p^-}(\Omega)}^\frac{p^-r}{p^--1}
			+\|\sqrt{\rho}\partial_t{\mathbf u}\|_{L^2(\Omega)}^\frac{8r(p^--1)}{5p^--6-r(2-p^-)})
\end{align}
holds for any fixed $\eta\in(0,1).$ Moreover, one obtains from \eqref{pt} that
\begin{align}\label{b2}
   &\frac{d}{dt}\int_\Omega\rho|\partial_t{\mathbf u}|^2dx+{\mathcal J}_\Phi({\mathbf u})\nonumber\\
   &\leqslant C\int_\Omega\rho|\partial_t{\mathbf u}|^2|\nabla{\mathbf u}|dx+\int_\Omega\rho|{\mathbf u}\|\nabla{\mathbf f}|\partial_t{\mathbf u}|dx
		+\int_\Omega\rho|{\mathbf u}\|{\mathbf f}\|\nabla\partial_t{\mathbf u}|~dx
		+\int_\Omega\rho|\partial_t{\mathbf f}\|\partial_t{\mathbf u}|\nonumber~dx\\
   &~~~+\int_\Omega\gamma \rho^\gamma|\divv{\mathbf u}\|\divv\partial_t{\mathbf u}|~dx
			+\int_\Omega|\nabla P\|{\mathbf u}\|\divv\partial_t{\mathbf u}|~dx
			+\int_\Omega\rho|{\mathbf u}\|\partial_t{\mathbf u}\|\nabla\partial_t{\mathbf u}|~dx\nonumber\\
			&~~~+\int_\Omega\rho|{\mathbf u}\|\nabla{\mathbf u}|^2|\partial_t{\mathbf u}|~dx
            +\int_\Omega\rho|{\mathbf u}|^2|\nabla^2{\mathbf u}\|\partial_t{\mathbf u}|~dx\nonumber\\
			&~~~+\int_\Omega\rho|{\mathbf u}|^2|\nabla{\mathbf u}\|\nabla\partial_t{\mathbf u}|~dx
            +\int_\Omega|\widetilde{D}{\mathbf u}|^{p(t,x)+\sigma}\:dx\nonumber\\
            &\leqslant \sum_{i=1}^{11}A_i.
\end{align}
where $\sigma\in(0,1).$ Based on \eqref{I}-\eqref{J} and Sobolev embedding, one estimates $A_i~(i=1,2,\cdots,11)$ as follows
\begin{align*}
A_1&\leqslant C\|\rho\|_{L^\infty(\Omega)}^\frac{1}{2}\|\sqrt{\rho}\partial_t{\mathbf u}\|_{L^3(\Omega)}^2\|\nabla{\mathbf u}\|_{L^3(\Omega)}\\
   &\leqslant C\|\sqrt{\rho}\partial_t{\mathbf u}\|_{L^2(\Omega)}^\frac{4}{7}
           \|\sqrt{\rho}\partial_t{\mathbf u}\|_{L^\frac{15}{4}(\Omega)}^\frac{10}{7}\\
   &\leqslant C\|\sqrt{\rho}\partial_t{\mathbf u}\|_{L^2(\Omega)}^\frac{4}{7}\|\nabla\partial_t{\mathbf u}\|_{L^\frac{5}{3}(\Omega)}^\frac{10}{7}\\
   &\leqslant C\|\sqrt{\rho}\partial_t{\mathbf u}\|_{L^2(\Omega)}^\frac{4}{7}
           (\int_\Omega(|\nabla\partial_t{\mathbf u}|^2|\widetilde{D}{\mathbf u}|^{p^--2})^\frac{5}{6}
           |\widetilde{D}{\mathbf u}|^\frac{5(2-p^-)}{6}\:dx)^\frac{6}{7}\\
	&\leqslant C\|\sqrt{\rho}\partial_t{\mathbf u}\|_{L^2(\Omega)}^\frac{4}{7}
			(\int_\Omega|\nabla\partial_t{\mathbf u}|^2|\widetilde{D}{\mathbf u}|^{p^--2}\:dx)^\frac{5}{7}
               (\int_\Omega|\widetilde{D}{\mathbf u}|^{5(2-p^-)}\:dx)^\frac{1}{7}\\
			&\leqslant C\|\sqrt{\rho}\partial_t{\mathbf u}\|_{L^2(\Omega)}^\frac{4}{7}{\mathcal J}_\Phi({\mathbf u})^\frac{5}{7}
			(\int_\Omega|\widetilde{D}{\mathbf u}|^{5(2-p^-)}\:dx)^\frac{1}{7}\\
			&\leqslant \eta{\mathcal J}_\Phi({\mathbf u})+C\|\sqrt{\rho}\partial_t{\mathbf u}|_{L^2(\Omega)}^2
			(\int_\Omega|\widetilde{D}{\mathbf u}|^{5(2-p^-)}\:dx)^2\\
			&\leqslant \eta{\mathcal J}_\Phi({\mathbf u})+C\|\sqrt{\rho}\partial_t{\mathbf u}\|_{L^2(\Omega)}^2,\\
A_2&\leqslant C \|\rho\|_{L^\infty(\Omega)}^\frac{1}{2}\|\sqrt{\rho}\partial_t{\mathbf u}\|_{L^2(\Omega)}
			\|\nabla{\mathbf f}\|_{L^2(\Omega)}\|{\mathbf u}\|_{L^\infty(\Omega)}\\
			&\leqslant C\|\nabla{\mathbf u}\|_{L^3(\Omega)}\|\sqrt{\rho}\partial_t{\mathbf u}\|_{L^2(\Omega)}\\
			&\leqslant C\|\sqrt{\rho}\partial_t{\mathbf u}\|_{L^2(\Omega)},\\
A_3&\leqslant C \|\rho\|_{L^\infty(\Omega)}\|{\mathbf f}\|_{L^6(\Omega)}\|\nabla\partial_t{\mathbf u}\|_{L^\frac{6}{5}(\Omega)}
			\|{\mathbf u}\|_{L^\infty(\Omega)}\\
			&\leqslant C\|\nabla{\mathbf u}\|_{L^3(\Omega)}(\int_\Omega(|\nabla\partial_t{\mathbf u}|^2
			|\widetilde{D}{\mathbf u}|^{p^--2})^\frac{3}{5}|\widetilde{D}{\mathbf u}|^\frac{3(2-p^-)}{5}\:dx)^\frac{5}{6}\\
			&\leqslant C(\int_\Omega|\nabla\partial_t{\mathbf u}|^2|\widetilde{D}{\mathbf u}|^{p^--2}\:dx)^\frac{1}{2}
			(\int_\Omega|\widetilde{D}{\mathbf u}|^\frac{3(2-p^-)}{2}\:dx)^\frac{1}{3}\\
			&\leqslant C{\mathcal J}^\frac{1}{2}_\Phi({\mathbf u})\|\widetilde{D}{\mathbf u}\|_{L^3(\Omega)}\\
			&\leqslant \eta{\mathcal J}_\Phi({\mathbf u})+C,\\
A_4&\leqslant C\|\rho\|_{L^\infty(\Omega)}^\frac{1}{2}\|\sqrt{\rho}\partial_t{\mathbf u}\|_{L^2(\Omega)}\|\partial_t{\mathbf f}\|_{L^2(\Omega)}\\
			&\leqslant C\|\sqrt{\rho}\partial_t{\mathbf u}\|_{L^2(\Omega)}^2+\|\partial_t\mathbf{f}\|^2_{L^2(\Omega)},\\
A_5&\leqslant C\|\rho\|_{L^\infty(\Omega)}^\gamma\|\nabla{\mathbf u}\|_{L^3(\Omega)}\|\nabla\partial_t{\mathbf u}\|_{L^\frac{3}{2}(\Omega)}\\
			&\leqslant C\|\nabla{\mathbf u}\|_{L^3(\Omega)}(\int_\Omega(|\nabla\partial_t{\mathbf u}|^2
			|\widetilde{D}{\mathbf u}|^{p^--2})^\frac{3}{4}|\widetilde{D}{\mathbf u}|^\frac{3(2-p^-)}{4}\:dx)^\frac{2}{3}\\
			&\leqslant C(\int_\Omega|\nabla\partial_t{\mathbf u}|^2|\widetilde{D}{\mathbf u}|^{p^--2}\:dx)^\frac{1}{2}
			(\int_\Omega|\widetilde{D}{\mathbf u}|^{3(2-p^-)}\:dx)^\frac{1}{6}\\
			&\leqslant C{\mathcal J}^\frac{1}{2}_\Phi({\mathbf u})\|\widetilde{D}{\mathbf u}\|_{L^3(\Omega)}^\frac{1}{2}\\
			&\leqslant \eta{\mathcal J}_\Phi({\mathbf u})+C,\\
A_6&\leqslant C\|\rho\|_{L^\infty(\Omega)}^{\gamma-1}\|\nabla\rho\|_{L^{3p^-}(\Omega)}\|\nabla\partial_t{\mathbf u}\|_{L^\frac{3}{2}(\Omega)}
			\|{\mathbf u}\|_{L^\infty(\Omega)}\\
			&\leqslant C\|\nabla\rho\|_{L^{3p^-}(\Omega)}\|\nabla{\mathbf u}\|_{L^3(\Omega)}
			(\int_\Omega(|\nabla\partial_t{\mathbf u}|^2|\widetilde{D}{\mathbf u}|^{p^--2})^\frac{3}{4}
			|\widetilde{D}{\mathbf u}|^\frac{3(2-p^-)}{4}\:dx)^\frac{2}{3}\\
			&\leqslant C\|\nabla\rho\|_{L^{3p^-}(\Omega)}(\int_\Omega|\nabla\partial_t{\mathbf u}|^2
			|\widetilde{D}{\mathbf u}|^{p^--2}\:dx)^\frac{1}{2}(\int_\Omega|\widetilde{D}{\mathbf u}|^{3(2-p^-)}\:dx)^\frac{1}{6}\\
			&\leqslant C\|\nabla\rho\|_{L^{3p^-}(\Omega)}{\mathcal J}^\frac{1}{2}_\Phi({\mathbf u})
              \|\widetilde{D}{\mathbf u}\|_{L^3(\Omega)}^\frac{1}{2}\\
            &\leqslant \eta{\mathcal J}_\Phi({\mathbf u})+C\|\nabla\rho\|_{L^{3p^-}(\Omega)}^2,\\
A_7&\leqslant C\|{\mathbf u}\|_{L^\infty(\Omega)}\|\rho\|_{L^\infty(\Omega)}^\frac{1}{2}\|\sqrt{\rho}\partial_t{\mathbf u}\|_{L^3(\Omega)}
			\|\nabla\partial_t{\mathbf u}\|_{L^\frac{3}{2}(\Omega)}\\
			&\leqslant C\|\nabla{\mathbf u}\|_{L^3(\Omega)}\|\sqrt{\rho}\partial_t{\mathbf u}\|_{L^2(\Omega)}^\frac{2}{7}
			\|\nabla\partial_t{\mathbf u}\|_{L^\frac{5}{3}(\Omega)}^\frac{5}{7}\|\nabla\partial_t{\mathbf u}\|_{L^\frac{3}{2}(\Omega)}\\
			&\leqslant C\|\sqrt{\rho}\partial_t{\mathbf u}\|_{L^2(\Omega)}^\frac{2}{7}
			(\int_\Omega(|\nabla\partial_t{\mathbf u}|^2|\widetilde{D}{\mathbf u}|^{p^--2})^\frac{5}{6}
			|\widetilde{D}{\mathbf u}|^{5(2-p^-)}\:dx)^\frac{3}{7}\cdot\\
			&~~~~(\int_\Omega(|\nabla\partial_t{\mathbf u}|^2|\widetilde{D}{\mathbf u}|^{p^--2})^\frac{3}{4}
			|\widetilde{D}{\mathbf u}|^\frac{3(2-p^-)}{4}\:dx)^\frac{2}{3}\\
			&\leqslant C\|\sqrt{\rho}\partial_t{\mathbf u}\|_{L^2(\Omega)}^\frac{2}{7}(\int_\Omega|\nabla\partial_t{\mathbf u}|^2
			|\widetilde{D}{\mathbf u}|^{p^--2}\:dx)^\frac{6}{7}(\int_\Omega|\widetilde{D}{\mathbf u}|^{5(2-p^-)}\:dx)^\frac{3}{14}\cdot\\
			&~~~~(\int_\Omega|\widetilde{D}{\mathbf u}|^{3(2-p^-)}\:dx)^\frac{1}{6}\\
			&\leqslant C\|\sqrt{\rho}\partial_t{\mathbf u}\|_{L^2(\Omega)}^\frac{2}{7}{\mathcal J}_\Phi({\mathbf u})^\frac{6}{7}
			\|\widetilde{D}{\mathbf u}\|_{L^3(\Omega)}^\frac{5}{7}\\
			&\leqslant C{\mathcal J}_\Phi({\mathbf u})^\frac{6}{7}\|\sqrt{\rho}\partial_t{\mathbf u}\|_{L^2(\Omega)}^\frac{2}{7}\\
			&\leqslant \eta{\mathcal J}_\Phi({\mathbf u})+C\|\sqrt{\rho}\partial_t{\mathbf u}\|_{L^2(\Omega)}^2,\\
A_8&\leqslant C\|\rho\|_{L^\infty(\Omega)}^\frac{1}{2}\|{\mathbf u}\|_{L^\infty(\Omega)}\|\sqrt{\rho}\partial_t{\mathbf u}\|_{L^3(\Omega)}
			\|\nabla{\mathbf u}\|_{L^3(\Omega)}^2\\
			&\leqslant C\|\sqrt{\rho}\partial_t{\mathbf u}\|_{L^2(\Omega)}^\frac{2}{7}
			\|\sqrt{\rho}\partial_t{\mathbf u}\|_{L^\frac{15}{4}(\Omega)}^\frac{5}{7}\|\nabla{\mathbf u}\|_{L^3(\Omega)}^3\\
			&\leqslant C\|\sqrt{\rho}\partial_t{\mathbf u}\|_{L^2(\Omega)}^\frac{2}{7}
              \|\nabla\partial_t{\mathbf u}\|_{L^\frac{5}{3}(\Omega)}^\frac{5}{7}\\
			&\leqslant C\|\sqrt{\rho}\partial_t{\mathbf u}\|_{L^2(\Omega)}^\frac{2}{7}(\int_\Omega(|\nabla\partial_t{\mathbf u}|^2
			|\widetilde{D}{\mathbf u}|^{p^--2})^\frac{5}{6}|\widetilde{D}{\mathbf u}|^{\frac{5(2-p^-)}{6}}\:dx)^\frac{3}{7}\\
			&\leqslant C\|\sqrt{\rho}\partial_t{\mathbf u}\|_{L^2(\Omega)}^\frac{2}{7}(\int_\Omega|\nabla\partial_t{\mathbf u}|^2
			|\widetilde{D}{\mathbf u}|^{p^--2}\:dx)^\frac{5}{14}(\int_\Omega|\widetilde{D}{\mathbf u}|^{5(2-p^-)}\:dx)^\frac{1}{14}\\
			&\leqslant C\|\sqrt{\rho}\partial_t{\mathbf u}\|_{L^2(\Omega)}^\frac{2}{7}{\mathcal J}_\Phi({\mathbf u})^\frac{5}{14}
			\|\widetilde{D}{\mathbf u}\|_{L^3(\Omega)}^\frac{3}{14}\\
			&\leqslant \eta{\mathcal J}_\Phi({\mathbf u})+C(\|\sqrt{\rho}\partial_t{\mathbf u}\|_{L^2(\Omega)}^2+1),\\
A_9&\leqslant C\|\rho\|_{L^\infty(\Omega)}^\frac{1}{2}\|{\mathbf u}\|_{L^\infty(\Omega)}^2
			\|\sqrt{\rho}\partial_t{\mathbf u}\|_{L^3(\Omega)}\|\nabla^2{\mathbf u}\|_{L^\frac{3}{2}(\Omega)}\\
			&\leqslant C\|\nabla{\mathbf u}\|_{L^3(\Omega)}^2\|\sqrt{\rho}\partial_t{\mathbf u}\|_{L^2(\Omega)}^\frac{2}{7}
			\|\nabla\partial_t{\mathbf u}\|_{L^\frac{5}{3}(\Omega)}^\frac{5}{7}\|\nabla^2{\mathbf u}\|_{L^\frac{3}{2}(\Omega)}\\
			&\leqslant C\|\sqrt{\rho}\partial_t{\mathbf u}\|_{L^2(\Omega)}^\frac{2}{7}\left(\int_\Omega(|\nabla\partial_t{\mathbf u}|^2
			|\widetilde{D}{\mathbf u}|^{p^--2})^\frac{5}{6}|\widetilde{D}{\mathbf u}|^{\frac{5(2-p^-)}{6}}\:dx\right)^\frac{3}{7}\\
			&~~~~~~~\left(\int_\Omega(|\nabla^2{\mathbf u}|^2|\widetilde{D}{\mathbf u}|^{p^--2})^\frac{3}{4}
			|\widetilde{D}{\mathbf u}|^{\frac{3(2-p^-)}{4}}\:dx\right)^\frac{2}{3}\\
			&\leqslant C\|\sqrt{\rho}\partial_t{\mathbf u}\|_{L^2(\Omega)}^\frac{2}{7}(\int_\Omega|\nabla\partial_t{\mathbf u}|^2
			|\widetilde{D}{\mathbf u}|^{p^--2}\:dx)^\frac{5}{14}(\int_\Omega|\widetilde{D}{\mathbf u}|^{5(2-p^-)}\:dx)^\frac{1}{14}\\
			&~~~~(\int_\Omega(|\nabla^2{\mathbf u}|^2|\widetilde{D}{\mathbf u}|^{p^--2})^\frac{3}{4}\:dx)^\frac{1}{2}
			(\int_\Omega|\widetilde{D}{\mathbf u}|^{\frac{3(2-p^-)}{4}}\:dx)^\frac{1}{6}\\
			&\leqslant C\|\sqrt{\rho}\partial_t{\mathbf u}\|_{L^2(\Omega)}^\frac{2}{7}{\mathcal J}_\Phi({\mathbf u})^\frac{5}{14}
			{\mathcal I}_\Phi({\mathbf u})^\frac{1}{2}\|\widetilde{D}{\mathbf u}\|_{L^3(\Omega)}^\frac{5}{7}\\
			&\leqslant \eta{\mathcal J}_\Phi({\mathbf u})+\varepsilon{\mathcal I}_\Phi({\mathbf u})+C\|\sqrt{\rho}{\mathbf u}_t\|_{L^2(\Omega)}^2,\\
A_{10}&\leqslant C\|\rho\|_{L^\infty(\Omega)}\|{\mathbf u}\|_{L^\infty(\Omega)}^2\|\nabla{\mathbf u}\|_{L^3(\Omega)}
			\|\nabla\partial_t{\mathbf u}\|_{L^\frac{3}{2}(\Omega)}\\
			&\leqslant C\|\nabla{\mathbf u}\|_{L^3(\Omega)}^3(\int_\Omega(|\nabla\partial_t{\mathbf u}|^2
			|\widetilde{D}{\mathbf u}|^{p^--2})^\frac{3}{4}|\widetilde{D}{\mathbf u}|^\frac{3(2-p^-)}{4}\:dx)^\frac{2}{3}\\
			&\leqslant C(\int_\Omega|\nabla\partial_t{\mathbf u}|^2|\widetilde{D}{\mathbf u}|^{p^--2}\:dx)^\frac{1}{2}
			(\int_\Omega|\widetilde{D}{\mathbf u}|^{3(2-p^-)}\:dx)^\frac{1}{6}\\
			&\leqslant C{\mathcal J}^\frac{1}{2}_\Phi({\mathbf u})\|\widetilde{D}{\mathbf u}\|_{L^3(\Omega)}^\frac{1}{2}\\
			&\leqslant \eta{\mathcal J}_\Phi({\mathbf u})+C,\\
A_{11}&\leqslant C\|\nabla{\mathbf u}\|_{L^3(\Omega)}^3\leqslant C,
\end{align*}
where $\eta\in (0,1)$ is arbitrary but fixed. Adding above estimate into \eqref{b2} and choosing $\eta$ sufficiently small, one obtains that
\begin{align}\label{32}
  \frac{d}{dt}\int_\Omega\rho|\partial_t{\mathbf u}|^2\:dx+{\mathcal J}_\Phi({\mathbf u})
			\leqslant\varepsilon{\mathcal I}_\Phi({\mathbf u})+ C(\|\sqrt{\rho}\partial_t{\mathbf u}|_{L^2(\Omega)}^2
              +\|\nabla\rho\|_{L^{3p^-}(\Omega)}^2+\|\partial_t\mathbf{f}\|^2_{L^2(\Omega)}+1)
\end{align}
holds for any fixed $\varepsilon\in (0,1).$ Note that
\begin{equation*}\label{3dlg-E19}
  \|{\mathbf u} \|_{W^{2,s}(\Omega)}
   \leqslant C \|-\rho\partial_t {\mathbf u}-\rho ({\mathbf u}\cdot\nabla) {\mathbf u}-\nabla \rho^\gamma\|_{L^{s}(\Omega)}
			+\|\nabla{\mathbf u}\|_{L^{p^-}(\Omega)}^{p^-}.
\end{equation*}
It follows from Lemma \ref{ine} that
\begin{align}\label{qwe}
			&\|{\mathbf u}\|_{W^{2,3p^-}(\Omega)}\nonumber\\
			&\leqslant C(\|\rho\partial_t{\mathbf u}\|_{L^{3p^-}(\Omega)}+\|\rho({\mathbf u}\cdot\nabla){\mathbf u}\|_{L^{3p^-}(\Omega)}
			+\|\rho{\mathbf f}\|_{L^{3p^-}(\Omega)}+\|\nabla \rho^\gamma\|_{L^{3p^-}(\Omega)}
			+\|\mathbb{ D}{\mathbf u}\|_{L^{p^-}(\Omega)}^{p^-})\nonumber\\
			&\leqslant C\left(\|\rho\|_{L^\infty(\Omega)}\|\partial_t{\mathbf u}\|_{L^{3p^-}(\Omega)}
			+\|\rho\|_{L^\infty(\Omega)}\|{\mathbf u}\|_{L^\infty(\Omega)}\|\nabla{\mathbf u}\|_{L^{3p^-}(\Omega)}
			+\|\rho\|_{L^\infty(\Omega)}\|{\mathbf f}\|_{L^{3p^-}(\Omega)}\right.  \nonumber   \\
			&\left.~~~ +\|\rho\|_{L^\infty(\Omega)}^{\gamma-1}\|\nabla\rho\|_{L^{3p^-}(\Omega)}+1\right)\nonumber\\
			&\leqslant C({\mathcal I}_\Phi({\mathbf u})+{\mathcal J}_\Phi({\mathbf u})+1)^\frac{1}{p^-}
			+C{\mathcal I}_\Phi({\mathbf u})^\frac{1}{p^-}+C\|\nabla\rho\|_{L^{3p^-}(\Omega)}+C.
			%\nonumber\\  &\leqslant \varepsilon{\mathcal I}_\Phi({\mathbf u})+\eta{\mathcal J}_\Phi({\mathbf u})+C(\|\nabla\rho\|_{L^{3p^-}(\Omega)}+1)
\end{align}
Recall that
\begin{eqnarray*}
	&&\frac{d}{dt}\|\rho\|_{W^{1,{3p^-}}(\Omega)}\leqslant C\|{\mathbf u}\|_{W^{2,{3p^-}}(\Omega)}\|\rho\|_{W^{1,{3p^-}}(\Omega)}.
\end{eqnarray*}
One deduces from \eqref{qwe} and the Young's inequality that
\begin{align}\label{33}
	\frac{d}{dt}\|\rho\|_{W^{1,3p^-}(\Omega)}&\leqslant C\|{\mathbf u}\|_{W^{2,3p^-}(\Omega)}\|\rho\|_{W^{1,3p^-}(\Omega)}\nonumber\\
		&\leqslant C({\mathcal I}_\Phi({\mathbf u})+{\mathcal J}_\Phi({\mathbf u})+1)^\frac{1}{p^-}\|\rho\|_{W^{1,3p^-}(\Omega)} 			
			+C{\mathcal I}_\Phi({\mathbf u})^\frac{1}{p^-}\|\rho\|_{W^{1,3p^-}(\Omega)}\nonumber\\
		&~~+C\|\nabla\rho\|_{L^{3p^-}(\Omega)}\|\rho\|_{W^{1,3p^-}(\Omega)}+C\|\rho\|_{W^{1,3p^-}(\Omega)}\nonumber\\
		&\leqslant \varepsilon{\mathcal I}_\Phi({\mathbf u})+\eta{\mathcal J}_\Phi({\mathbf u})+C(\|\rho\|_{W^{1,3p^-}(\Omega)}^\frac{p^-}{p^--1}+1),
\end{align}
holds for any fixed $\varepsilon,\eta\in (0,1).$  Summary, one obtains from \eqref{31}, \eqref{32} and \eqref{33} that
\begin{align}\label{3}
  &\frac{d}{dt}(\|\sqrt{\rho}\partial_t{\mathbf u}\|_{L^2(\Omega)}^2+\|\rho\|_{W^{1,3p^-}(\Omega)})
	 +{\mathcal I}^r_\Phi({\mathbf u})+{\mathcal J}_\Phi({\mathbf u})\nonumber\\
  &\leqslant C(1+\|\rho\|_{W^{1,3p^-}(\Omega)}^\frac{p^-r}{p^--1}
    +\|\sqrt{\rho}\partial_t{\mathbf u}\|_{L^2(\Omega)}^\frac{8r(p^--1)}{5p^--6-r(2-p^-)}+\|\partial_t\mathbf{f}\|^2_{L^2(\Omega)}).
\end{align}
Hence,
\begin{align}\label{333}
	&\frac{d}{dt}(\|\sqrt{\rho}\partial_t{\mathbf u}\|_{L^2(\Omega)}^2+\|\rho\|_{W^{1,3p^-}(\Omega)})
		\leqslant c_0(\|\sqrt{\rho}{\mathbf u}_t\|_{L^2(\Omega)}^2+\|\rho\|_{W^{1,3p^-}(\Omega)})^{\alpha}+c_1(1+\|\partial_t\mathbf{f}\|^2_{L^2(\Omega)})
\end{align}
for some $c_0>0$ and $c_1>0,$ where $\alpha=\max \{\frac{p^-r}{p^--1},\frac{4r(p^--1)}{5p^--6-r(2-p^-)}\}>1.$ Setting
		$$f(t)\triangleq \|\sqrt{\rho}\partial_t{\mathbf u}\|_{L^2(\Omega)}^2+\|\rho\|_{W^{1,3p^-}(\Omega)}\mbox{ and }
h(t)\triangleq c_1(1+\|\partial_t\mathbf{f}\|^2_{L^2(\Omega)}),$$
one rewrites \eqref{333} as
\begin{align*}
	f^\prime(t)\leqslant c_0f^{\alpha}(t)+h(t).
\end{align*}
Thanks to the compatibility condition
\begin{eqnarray*}
\sqrt{\rho_0(x)}(\sqrt{\rho_0(x)}\partial_t{\mathbf u}(0,x)+\sqrt{\rho_0(x)}{\mathbf u}_0\cdot \nabla{\mathbf u}_0-\sqrt{\rho_0(x)}{\mathbf g} )=0,
\end{eqnarray*}
it holds that
\begin{eqnarray*}
\sqrt{\rho_0(x)}\sqrt{\rho_0(x)}\partial_t{\mathbf u}(0,x)
=-\sqrt{\rho_0(x)}{\mathbf u}_0\cdot \nabla{\mathbf u}_0+\sqrt{\rho_0(x)}{\mathbf g} \in L^2(\Omega).
\end{eqnarray*}
So, it follows from $\partial_t\mathbf{f}\in L^2(0,T^*;L^2(\Omega))$ and Lemma \ref{local-gronwall} that
\begin{align*}
	f(t)\leqslant H(t)\left(1-(\alpha-1)c_0tH^{\alpha-1}(t)\right)^{-\frac{1}{\alpha-1}}
\end{align*}
holds for any $t\in(0,\min\{\frac{1}{c_0(\alpha-1)}\left(f(0)+c_1(T^*+\|\partial_t{\bf f}\|_{L^2(0,T^*;L^2()\Omega)})\right)^{-(\alpha-1)},T^*\})$, where
$$H(t)=f(0)+\int_0^th(s)ds \mbox{ and }f(0)\triangleq \|\sqrt{\rho_0}\partial_t{\mathbf u}(0,\cdot)\|_{L^2(\Omega)}^2+\|\rho_0\|_{W^{1,3p^-}(\Omega)}.$$
		
Therefore,
\begin{align}\label{E-100}
	\sup\limits_{0\leqslant t\leqslant T}(\|\sqrt{\rho}\partial_t{\mathbf u}\|_{L^2(\Omega)}^2+\|\rho\|_{W^{1,3p^-}(\Omega)})
			+\int_{0}^{T}\left({\mathcal I}^r_\Phi({\mathbf u})+{\mathcal J}_\Phi({\mathbf u})\right)\:dt\leqslant C
\end{align}
holds for any $T\in(0,\min\{\frac{1}{c_0(\alpha-1)}\left(f(0)+c_1(T^*+\|\partial_t{\bf f}\|_{L^2(0,T^*;L^2()\Omega)})\right)^{-(\alpha-1)},T^*\}).$ In addition, integrating \eqref{qwe} with respect to $t$ yields that
\begin{align}\label{E-101}
	\int_{0}^{T}\|{\mathbf u}\|_{W^{2,3p^-}(\Omega)}^{p^-}\:dt\leqslant C,
\end{align}
holds for any $T\in(0,\min \{\frac{1}{(\alpha-1)Cf^{\alpha-1}(0)},T^*\}).$
\end{proof}
	
In view of Lemma \ref{uu3}, one takes $t_1\in(0,\min \{\frac{1}{(\alpha-1)Cf^{\alpha-1}(0)},T^*\})$ and finds that
	$$(\rho,{\mathbf u})(x,t_1)\triangleq\lim\limits_{t\rightarrow t_1}(\rho,{\mathbf u}),$$
which satisfies the conditions imposed on the initial data \eqref{CZ} at the time $t=t_1$. Furthermore,
\begin{align*}
	\left(-\divv((1+|{\mathbb D}{\mathbf u}|^2)^{\frac{p(t, x)-2}{2}}{\mathbb D}{\mathbf u})+\nabla \rho^\gamma\right)|_{t=t_1}
		=\lim\limits_{t\rightarrow t_1} \rho^\frac{1}{2}(x,t_1)g(x),
\end{align*}
where $g(x)\triangleq\lim\limits_{t\rightarrow t_1}(\rho^\frac{1}{2}(\partial_t{\mathbf u}+{\mathbf u}\cdot\nabla{\mathbf u}))(x,t_1)\in L^{3p^-}(\Omega).$ Thus, $(\rho,{\mathbf u})(x,t_1)$ also satisfies \eqref{condition}. So, one can take $(\rho,{\mathbf u})(x,t_1)$ as the initial data and apply Theorem \ref{theo1} to extend the local strong solution beyond $T^*$. This contradicts the assumption on $T^*.$

\bmhead{Conflict of Interest}
On behalf of all authors, the corresponding author states that there is no conflict of interest.
	
\bmhead{Acknowledgements}
	
The authors would like to be grateful to the anonymous referee for his (her) helpful comments and advice. Fang's research is supported by National Natural Science Foundation of China (No. 11501445). Guo's research is supported by National Natural Science Foundation of China (No. 11931013).

\end{document}